\documentclass[11pt,reqno]{amsproc}

\title[Scattering for focusing supercritical wave equations in odd dimensions]{Scattering for radial bounded solutions of focusing supercritical wave equations in odd dimensions}

\author[G.~Camliyurt]{Guher Camliyurt}
\address{Department of Mathematics, University of Chicago, Chicago, IL 60637}
\email{camliyurt@math.uchicago.edu}

\author[C. E.~Kenig]{Carlos E. Kenig}
\address{Department of Mathematics, University of Chicago, Chicago, IL 60637}
\email{cek@math.uchicago.edu}

\usepackage[margin=1in]{geometry}
\usepackage{amsmath, amsthm, amssymb}
\usepackage{mathrsfs}
\usepackage{graphicx}
\usepackage{times}
\allowdisplaybreaks

\usepackage[usenames,dvipsnames]{color}
\usepackage[colorlinks=true, pdfstartview=FitV, linkcolor=blue, citecolor=blue, urlcolor=blue]{hyperref}
\usepackage{amsmath}
\usepackage{amsfonts}
\usepackage{amssymb}
\usepackage{amsmath, amsthm, amssymb, enumerate}

\usepackage{fancyhdr}
\usepackage{mathtools}
\usepackage{bbm}

\DeclarePairedDelimiter\evaluat{.}{\rvert}

\numberwithin{equation}{section}

\definecolor{colorgggg}{rgb}{0.5,0.0,0.4}%%%out
\def\cole{}

\def\colb{\color{black}}
\def\inon#1{\,\,\,\,\,\,\hbox{#1}}                %in or on
\def\norm#1{\left\Vert #1\right\Vert} %norm
\def\llabel#1{\nonumber}
\newtheorem{Theorem}{Theorem}[section]
\newtheorem{Corollary}[Theorem]{Corollary}
\newtheorem{Proposition}[Theorem]{Proposition}
\newtheorem{Lemma}[Theorem]{Lemma}
\newtheorem{Remark}[Theorem]{Remark}

\newtheorem{Notation}[Theorem]{Notation}
\newtheorem{Claim}[Theorem]{Claim}
\newtheorem*{Claim*}{Claim}
\newtheorem*{Proposition*}{Proposition}
\newtheorem{definition}[Theorem]{Definition}

\def\comma{ {\rm ,\qquad{}} }            %comma in a formula
\def\les{\lesssim}

\def\indeq{\quad{}}                     %indentation in formulas
\def\indeqtimes{\indeq\indeq\times}
\def\period{.}                           %period in a formula
\def\Imax{\textmd{I}_{\scriptsize{\mbox{max}}}}

\begin{document}

\nocite{*} %checking unused references

\begin{abstract}
We consider the wave equation with an energy-supercritical focusing nonlinearity in general odd dimensions. 
We prove that any radial solution that remains bounded in the critical Sobolev space is global and 
scatters to a linear solution. 
\end{abstract}


%\noindent\thanks{\em{MSC Classification:  35Q35.}}

\thanks{{\em Date}. \today}

\maketitle

\tableofcontents

\section{Introduction}
\label{sec:intro}
In this paper we consider the Cauchy problem for the semilinear focusing wave equation 
\begin{align}
   \begin{split}
    & \partial^2_{t} u
    - \Delta u 
    -  |u|^{p-1} u
     = 0   \comma \mbox{in } \mathbb{R}^d \times I ,
    \\
   &  \vec{u} (0)
    = (u_0, u_1) \in \dot{H}^{s_p} \times \dot{H}^{s_p -1} (\mathbb{R}^d)
    \end{split}
   \label{wave equations}
  \end{align}
in the energy-supercritical radial setting. Here, the set $I$ is an interval around $0$ and  the dimension $d$ is an odd 
integer greater than or equal to seven.    
Also, we consider the energy-supercritical nonlinearities:
\begin{equation}
s_p = \frac{d}{2} - \frac{2}{p-1} \comma p > \frac{d+2}{d-2}.
\label{eq0.1}
\end{equation}
%\begin{equation}
%p \geq 1+ \frac{4}{d-2}  \comma s_p = \frac{d}{2} - \frac{2}{p-1}.
%\label{eq0.1}
%\end{equation}
The \textit{homogeneous Sobolev space }
$\dot{H}^{s_p}$ 
denotes the closure of the Schwartz functions under the norm 
 \begin{equation*}
 \norm{f}_{\dot{H}^{s_p}  (\mathbb{R}^d)} = \norm{D^{s_p} f}_{L^2 (\mathbb{R}^d)} .
 \end{equation*}
The class of solutions to the Cauchy problem \eqref{wave equations} is invariant under the scaling 
\begin{equation}
\vec{u} (t,x) \mapsto \left( {\lambda^{- a_p }} u \left( {t}/{\lambda}, {x}/{\lambda} \right), 
{\lambda^{-1- a_p }}  \partial_t u \left( {t}/{\lambda}, {x}/{\lambda}  \right)   \right)
\label{eq0.2}
\end{equation}
where $a_p = 2/ (p-1)$. The scaling in \eqref{eq0.2} also determines the critical regularity space for the initial data: we note that
the  $\dot{H}^{s_p} \times \dot{H}^{s_p -1}$
norm of $(u_0, u_1)$ stays invariant under \eqref{eq0.2}. 
%Since we let $p \geq 3$, 
Due to the lower bound on the $p$ exponent, 
the space for initial data is equipped with  $s_p > 1$, 
which places  the Cauchy problem \eqref{wave equations} in an energy-supercritical regime. 

In the present article, 
we concentrate on the cases where the nonlinearity is analytic, in particular we let
\begin{equation*}
 p = 2m+1,~  m \in \mathbb{N}
\end{equation*}
and
we prove that any radial solution of the Cauchy problem \eqref{wave equations} that is bounded in the critical regularity space 
$\dot{H}^{s_p} \times \dot{H}^{s_p -1}$ (throughout its maximal interval of existence) 
must be global and must scatter to a linear solution. 
We remark that the analogous assertions were established by Duyckaerts, Kenig, and Merle \cite{DKM14} in three dimensions and by  
Dodson and Lawrie \cite{DL15} in five dimensions.
% While we particularly work in seven dimensions, we expect that our approach 
%will generalize to all higher odd dimensions. 
Our main result is below. 
\begin{Theorem}
\label{main thm_1}
Let $p \geq 3$ be an odd integer and let $\vec {u} (t)$  be a radial solution to the equation \eqref{wave equations} with maximal interval of existence $\Imax (\vec{u}) = (T_{-} (\vec{u}), T_{+} (\vec{u}))$ such that 
\begin{align}
\sup_{t \in (0, T_+ (\vec{u}))}   \Vert (u(t), \partial_t u(t)) \Vert_{\dot{H}^{s_p} \times \dot{H}^{s_p -1} (\mathbb{R}^d)} 
< \infty
.
\end{align}
Then, $\Imax ( \vec{u}) \cap (0, \infty)  = (0, \infty)$  and $\vec{u} (t)$ scatters to a free wave as $t \to \infty$. 

\end{Theorem}

A direct consequence of Theorem~\ref{main thm_1} is that any finite time blow-up solution must admit a critical 
Sobolev norm diverging to infinity along a sequence of times. 
We point out that under slightly different boundedness conditions, there are blow-up constructions for the focusing 
nonlinear wave equation in the energy-supercritical regime; see Collot's result \cite{Co18} for a family of  blow-up solutions which become singular 
via concentration of a soliton profile. In \cite{Co18} solutions break down at finite time
even though the norms below critical scaling remain bounded, i.e., 
\begin{equation}
\limsup_{t \nearrow T} \norm{\vec{u} (t)}_{\dot{H}^s \times \dot{H}^{s-1} (\mathbb{R}^d)} < \infty
\label{eq0.5}
\end{equation}
for $s \in [1, s_p )$, where $T$ is the blow-up time 
(we note that the critical norms of these solutions are unbounded over $[0, T)$). 

Power-type nonlinear wave equations have received particular attention in the energy-critical setting
\begin{align}
 \begin{split}
    & \partial^2_{t} u
    - \Delta u 
    = \pm  |u|^{\frac{4}{d-2}} u
       \comma \mbox{in } \mathbb{R}^d \times I ,
    \\
   &  \vec{u} (0)
    = (u_0, u_1) \in \dot{H}^{1} \times L^2 (\mathbb{R}^d)
    \end{split}
   \label{wave equations_ec}
\end{align}
where $d \geq 3$ denotes the dimension and the signs $+$, $-$ correspond to the focusing and defocusing cases, respectively. 
For the defocusing problem, global existence and scattering results were first obtained in three dimensions by Struwe \cite{S89} in the radial setting, and then by 
Grillakis \cite{G90} in the general setting. The results were then generalized to higher dimensions by Grillakis \cite{G92} , Shatah-Struwe \cite{SS93, SS94, SS98}, 
Bahouri-Shatah \cite{BS98}, and Kapitanski \cite{K94}. 

In the energy-critical focusing case, the asymptotic dynamics of solutions with large initial data require a much closer look. 
In 1974, Levine \cite{L74} showed that if $(u_0, u_1) \in \dot{H}^1 \times L^2$ is a non-zero initial data where 
\begin{equation*}
E (u_0, u_1) = \int \frac{1}{2} ( |u_1|^2 + |\nabla u_0|^2 ) - \frac{d-2}{2 d} \int |u_0|^{\frac{2d}{d-2}} <0  
\end{equation*}
then the solution must break down in finite time.  Although this work does not provide an answer on the nature of the blow-up, 
it stimulated the search for subsequent blow-up constructions in the literature. 
%may be seen as a guide to the subsequent blow-up constructions  in the literature.  

Firstly, we observe that 
\begin{equation*}
\varphi (t,x) = \left( \frac{(d-2) d}{4} \right)^{\frac{d-2}{4}} (1-t)^{- \frac{(d-2)}{2}}
\end{equation*}
is a solution to the ODE, $\partial_{tt} \varphi = |\varphi |^{\frac{4}{d-2}} \varphi$, which fails to be in $\dot{H}^1 \times L^2$. Nevertheless, by truncating 
the data and using finite speed of propagation, we may find a solution $u(x,t)$ to the focusing problem \eqref{wave equations_ec} that has unbounded critical Sobolev norm, i.e.,  $\lim_{t \nearrow 1} \norm{u (x,t)}_{\dot{H}^1 \times L^2 (\mathbb{R}^d)} = \infty$. We refer to this behaviour as type-I blow up. 

Additionally, if a nonzero solution $\vec{u} (t)$ of \eqref{wave equations_ec} has critical norm that remains bounded on $(0, T_+ (\vec{u}))$, namely
\begin{equation*}
\sup_{0< t< T_+ (\vec{u})} \norm{ \vec{u} (t)}_{\dot{H}^1 \times L^2 (\mathbb{R}^{d})} < \infty ,
\end{equation*}
then we call $\vec{u} (t)$ a type-II solution.  
There are type-II solutions to the focusing problem that blow-up in finite time, i.e., type-II solutions with $T_+ (\vec{u}) < \infty$. Such behaviour is generally referred to as 
type-II blow-up. 
 In \cite{KST09} Krieger, Schlag, and Tataru  constructed a radial type-II 
blow-up solution for the energy critical focusing problem \eqref{wave equations_ec} 
in three dimensions using the unique radial ground state solution $W$ for the underlying elliptic equation.  
The blow-up occurs at $T_- =0$ and 
their blow-up solution
$u(t,x)$ has the form 
\begin{equation}
u (t,x) = \lambda (t)^{\frac{1}{2}} W (\lambda (t) x) + \eta (t, x)
\end{equation}
where as $t \to 0$, the scaling parameter $\lambda (t) = t^{-1-\nu}$ diverges to infinity  and the local energy of the term $\eta$ inside the light cone converges to $0$. 
The latter limiting behaviour is given by
\begin{equation}
\mathcal{E}_{loc} (\eta) = \int_{|x| <t} (\eta_t^2 + |\nabla \eta|^2 +\eta^6 )~dx \to 0 \quad \mbox{as~} t \to 0.  
\end{equation}
 Later on, the original condition $\nu > 1/2$ was relaxed to $\nu >0$ in \cite{KS14}. 
 Furthermore, Donninger, Huang, Krieger, and Schlag investigated  the question of existence of other rescaling functions 
 which could yield similar type of blow-up solutions. Indeed, their results in \cite{DHKS14} exhibit an uncountable family of admissible rates
  for $\lambda (t)$ that are not of  pure-power type. 
  
  In \cite{KM1, KM2}, Kenig and Merle developed a program to address the \textit{ground state conjecture} for critical focusing problems. 
  In particular, for the energy-critical focusing problem \eqref{wave equations_ec} they established that 
  the energy of the ground state solution W was a threshold for global existence and scattering. 
 The method behind these results, which has come to be known as the \textit{concentration-compactness/rigidity method} has found 
 numerous applications within nonlinear dispersive and wave equations.  We refer the reader to \cite[Introduction]{K15}
 for more details and further references. Moreover, in a series of articles \cite{DKM11, DKM12.2, DKM13}
 Duyckaerts, Kenig, and Merle gave a classification of  solutions that remain bounded  in the three dimensional radial case. 
  Particularly, in \cite{DKM13} the authors established the \textit{soliton resolution conjecture} in dimension three, 
 which yields that any type-II radial solution asymptotically resolves 
 into a sum of decoupled solitary waves and a radiation term in $\dot{H}^1 \times L^2$. 
 
 In the energy supercritical regime, global in time well-posedness and scattering results 
 accompanied by the boundedness of the critical Sobolev norm 
 were 
 obtained 
 %under boundedness of the critical Sobolev norm 
 firstly for the defocusing case. 
  In \cite{KM11},  Kenig and Merle  addressed these assertions in the radial setting for dimension three. 
  Killip and Visan generalized these claims
  to all  dimensions in  \cite{KV11} for a range of energy supercritical exponents. 
  Analogous results for the cubic nonlinear wave equation were studied by Bulut \cite{B14} in dimension five;  
  see also  \cite{B12, B15, KV11.2} for results addressing the non-radial setting. 
  
  Utilizing the channel of energy method, Duyckaerts, Kenig, and Merle extended the global well-posedness and scattering results of \cite{KM11} to the focusing case in \cite{DKM14}. Additionally, similar results were obtained in dimension three by Duyckaerts and Roy in \cite{DR17}, and by Duyckaerts and Yang in \cite{DY18}
  with an improvement on the uniform boundedness condition. 
  In \cite{DL15}, Dodson and Lawrie studied the focusing cubic wave equation in five dimensions as well as the one-equivariant wave maps equations in three dimensions. 
  We note that the methods in \cite{DL15} apply to all supercritical exponents, yielding analogous results to \cite{DKM14}.    
  For results addressing the nonradial setting, see \cite{DS16, BD21} and the references cited therein. Also, we refer the reader to 
 \cite{KV10} for a corresponding  result addressing the focusing nonlinear Schr\"{o}dinger equation. 

 Analogous to the results in \cite{DKM14, DL15}, in this article we are concerned with 
type-II solutions, namely solutions to the problem \eqref{wave equations} for which 
\begin{equation}
\sup_{t \in (0, T_+ (\vec{u})) } \norm{ \vec{u} (t)  }_{\dot{H}^{s_p} \times \dot{H}^{s_p -1} (\mathbb{R}^d)} < \infty 
. 
\label{eq0.3}
\end{equation}
Our main result Theorem~\ref{main thm_1} shows that radial solutions to \eqref{wave equations} with  \eqref{eq0.3} achieve
$T_+ = \infty$ and they scatter. 

Prior to the blow-up scenario constructed in \cite{Co18}, 
a related result for the focusing NLS  was shown in \cite{MRR15} by Merle, Raphael, and Rodnianski. Both of these blow-up scenarios 
are constructed in large dimensions $d \geq 11$, addressing sufficiently large energy-supercritical exponents $p$. 
 In \cite{DD21}, Dai and Duyckaerts proved the existence of a countable family of self-similar blow-up solutions 
to the focusing energy supercritical wave equations   under the assumption \eqref{eq0.5}
in dimension three and 
for a range of supercritical exponents in dimensions $d \geq 4$. Additionally, Merle, Raphael, Rodnianski, and Szeftel recently  
obtained a finite time blow-up construction for defocusing NLS. In \cite{MRRS22}, they
showed the existence of a smooth well-localized initial data that corresponds to a unique blow-up solution in dimensions $d \geq 5$ for a range of supercritical exponents. 
%This distinguishes the main result of the present paper and its analogues in \cite{DL15, DKM12} for dimensions three and five  
%due to the blow-up scenario constructed in \cite{Co18}. 

There are also several works in the energy-subcritical regime addressing the asymptotic dynamics of type-II solutions. 
For instance,  see \cite{MZ05, MZ15} for blow-up behaviour of solutions to the focusing nonlinear wave equations. 
In addition, conditional scattering results that are parallel 
to Theorem~\ref{main thm_1}  may be found in \cite{S13, S20, DL15.2, R17, DLMM20} for dimensions three, four, and five.

\subsection{Overview of the proof of Theorem~\ref{main thm_1}} 
The general framework for the proof of Theorem~\ref{main thm_1}  follows closely the concentration compactness/rigidity method introduced by Kenig and Merle 
in \cite{KM1, KM2}, and extended into the energy-supercritical regime in the works \cite{DKM14, DL15}. 

First, we establish the local wellposedness theory for the Cauchy problem \eqref{wave equations} 
by means of standard techniques based on the Strichartz estimates. 
These estimates guarantee that if the initial data is sufficiently small in $\dot{H}^{s_p} \times \dot{H}^{s_{p-1}}$, then the corresponding solution $\vec{u} (t)$ is a global solution and it scatters to free waves in both time directions as $t \to \pm \infty$. Nevertheless, the small data theory will not be sufficient to analyze global dynamics of solutions with large data. The goal of our main result is to address the asymptotic dynamics of such solutions in the energy-supercritical radial setting.

We observe that Theorem~\ref{main thm_1} is equivalent to the fact that the claim below holds for all $A >0$. 
\begin{Claim*}
Let $\vec{u}(t)$ be a radial solution of the Cauchy problem \eqref{wave equations} with $\Imax = (T_- , T_+)$ such that
\begin{equation}
\sup_{t \in [0, T_+) } \norm{(u (t), \partial_t u (t))}_{\dot{H}^{s_p} \times \dot{H}^{s_p -1} (\mathbb{R}^d)}  \leq A.
\label{eq0.4}
\end{equation}
Then, $T_+ = \infty$ and $\vec{u} (t)$ scatters to a free wave in the positive time direction.  
\end{Claim*}
The small data theory guarantees that the claim is true for sufficiently small $A>0$. 
If Theorem~\ref{main thm_1} failed, this would lead to a critical value $A_C >0$ so that the claim above holds 
for all $A < A_C$, and for $A \geq A_C$ it fails. The profile decomposition results for the wave equations then 
allow for an extraction of a minimal solution to \eqref{wave equations}, called  ``critical solution", which does not scatter. 
In this context, minimality refers to the size of the solution in the accompanying condition assumed under Theorem~\ref{main thm_1}.  
In the present application of the concentration compactness procedure, we appeal to a profile decomposition result \cite{B10} by Bulut, which extends 
the earlier work of Bahouri-Gerard \cite{BG99} from $\dot{H}^1 \times L^2$ initial data in three dimensions to $\dot{H}^{s} \times \dot{H}^{s -1}$ initial data in 
higher dimensions with $s \geq 1$. Such critical solutions are shown to have pre-compact trajectories, up to modulation, in the space $\dot{H}^{s_p} \times \dot{H}^{s_p -1}$, which is the main property that under further analysis produces a contradiction. 

As noted above, in order to prove Theorem~\ref{main thm_1} we need to show the non-existence of a non-zero critical element. To that end, we follow the 
rigidity argument developed for the energy-supercritical wave equations in \cite{KM11, DKM14, DL15}. 
As a first step in that direction, we define and study solutions that exhibit  a compactness property: a solution $\vec{u} (t)$ is said to have \textit{the compactness property}
if there exists $\lambda: \Imax (\vec{u}) \to (0, \infty)$ so that the set
\begin{align*}
K = \left\lbrace \left( \frac{1}{\lambda(t)^{\frac{2}{p-1}}} u \left(t,  \frac{x}{\lambda (t)}  \right), \frac{1}{\lambda(t)^{\frac{2}{p-1}+1}} \partial_t u\left( t, \frac{x}{\lambda (t)} \right) \right) 
: t \in \Imax (\vec{u}) \right\rbrace
\end{align*}
has compact closure in $\dot{H}^{s_p} \times \dot{H}^{s_p -1} (\mathbb{R}^d)$. Such solutions are obtained from critical solutions via convergence: 
if $\vec{u} (t)$ is a non-scattering solution to \eqref{wave equations} that satisfies 
\begin{align*}
\sup_{0 < t < T_+ (\vec{u})} \norm{ (u (t), \partial_t u (t))}_{\dot{H}^{s_p} \times \dot{H}^{s_p -1} (\mathbb{R}^d) } < \infty, 
\end{align*}
then there exists $t_n \to T_+ (u)$ such that, up to modulation, $({u} (t_n), \partial_t u (t_n) )$ converges to $(v_0, v_1)$ in $\dot{H}^{s_p} \times \dot{H}^{s_p -1}$, where the corresponding solution $\vec{v} (t)$ has the compactness property.
%We then appeal to a profile decomposition result by Bulut  from \cite{B10} 
%to extract a minimal solution to \eqref{wave equations}, called  ``critical solution", which does not scatter. 
%In this context, minimality refers to the size of the solution in the accompanying condition assumed under Theorem~\ref{main thm_1}.  

In \cite{DKM15}, the authors showed that a solution with the compactness property must be global. 
Their result addresses focusing wave equations both in the energy-critical and energy-supercritical cases (as well as energy-supercritical Schr\"{o}dinger equations) 
and it directly precludes the possibility of a self similar solution that remains bounded in critical Sobolev norm and  blows up in finite time.

Having eliminated a finite time blow up scenario, our main rigidity result takes the following form.
\begin{Proposition*}[Proposition~\ref{Proposition1}]
 Let $\vec{u} (t)$ be a radial solution of \eqref{wave equations} with $\Imax (\vec{u}) = \mathbb{R}$, which has the compactness property. Suppose that we also have $\lambda (t) > A_0 > 0 $  for all $t \in \mathbb{R}$. Then, $\vec{u} \equiv (0,0)$.  
\end{Proposition*} 

In order to implement the rigidity argument, we first show that solutions with the compactness property have better decay than we have assumed. More specifically, 
we prove that $\vec{u} (t) \in \dot{H}^1 \times L^2 (\mathbb{R}^d)$ for all $t \in \mathbb{R}$, and in fact the trajectory 
\begin{align*}
\{\vec{u}(t): t \in \mathbb{R} \}
\end{align*}
has compact closure in $\dot{H}^1 \times L^2 (\mathbb{R}^d)$. As a direct consequence, we obtain the following vanishing: For all $R >0$,
\begin{align*}
\limsup_{t \to -\infty} \norm{ \vec{u} (t)}_{\dot{H}^1 \times L^2 (r \geq R + |t|)} = 
\limsup_{t \to \infty} \norm{ \vec{u} (t)}_{\dot{H}^1 \times L^2 (r \geq R + |t|)} = 0.
\end{align*}
The additional decay that lands the solution trajectories in the energy space $\dot{H}^1 \times L^2 (\mathbb{R}^d)$ is achieved via the double-Duhamel trick. This technique
was introduced by Colliander, Keel, Staffilani, Takaoka, and Tao in \cite{CKSTT08} and has been extensively utilized (see for instance \cite{KV10, KV10.2,  KV11.2, KV13, B12}).
It was also employed in \cite{DL15}
for the analogous problem in dimension five. In this step,  we rely on the assumption that the exponent $p$ is an odd integer, which implies that the nonlinearity is an analytic type. Generalizing this step to non-integer  exponents that are greater than two introduces technical difficulties. We hope to return to this at a later time. 

An essential ingredient of the rigidity argument for supercritical focusing type equations in the radial setting  is the so-called \textit{channels of energy} method. 
These estimates were first considered for linear radial wave equation in three dimensions in \cite{DKM11}, and for the five dimensional case in \cite{KLS}. 
Both of these results were then utilized in the adaptation of rigidity arguments to supercritical focusing nonlinearities; see \cite{DKM14, DL15} and references therein. 
Here we rely on the general form of 
the channel of energy estimates, which were proven  in all odd dimensions  by Kenig, Lawrie, Liu, and Schlag in \cite{KLLS15}: More specifically, considering the solution 
$V(t, x)$ 
to the linear wave equation with radial initial data $(f, g) \in \dot{H}^1 \times L^2 (\mathbb{R}^d)$, the result in \cite{KLLS15} states that in any odd dimension $d$, 
the radial energy solution $V(t,r)$ satisfies
\begin{align*}
\max_{\pm} \lim_{t \to \pm \infty} 
 \int_{r \geq |t|+R} |\nabla_{x,t} V(t, r)|^2 r^{d-1} dr   
 \geq \frac{1}{2} \Vert \pi^{\perp}_{P(R)} (f, g) \Vert^2_{\dot{H}^1 \times L^2 (r \geq R,~r^{d-1} dr)}
\end{align*}
 for all $R>0$. 
 Similar to \cite{DKM14, DL15}, the estimates above
can be directly employed in our rigidity argument.  
 
The operator $ \pi^{\perp}_{P(R)}$  on the right hand side denotes the orthogonal projection onto the complement of the subspace $P(R)$ in
  $\dot{H}^1 \times L^2 (r \geq R,~ r^{d-1} dr)$. When $d=1$, we have $P(R) = \{ (0,0) \}$, and when $d=3$,  $P(R)$ is the line 
\begin{align*}
P(R) = \{ (k/r , 0): k \in \mathbb{R} \}.  
\end{align*}
The formula for the subspace in general odd dimensions is given by
\begin{align*}
P (R) = \mbox{span} \left\lbrace (r^{2 k_1 -d}, 0), (0, r^{2 k_2 -d}): k_1= 1,2, \cdots, \left[ \frac{d+2}{4}\right] ; k_2 = 1, 2, \cdots, \left[ \frac{d}{4} \right]  \right\rbrace .
\end{align*}
As a result, in order to adapt the rigidity arguments into our setting, where $d \geq7$, 
we need to project away from a ${(d-1)}/{2}$ dimensional subspace rather than a line as in \cite{DKM14} or a $2$-plane as in \cite{DL15}. The change in the level of complexity also manifests at every step of the rigidity argument. 

Another tool needed for the rigidity argument is a one-parameter family of solutions to the underlying elliptic equation, whose behaviour 
near infinity are similar to that of $u (t,r)$ given in the main rigidity result.  Similar to the focusing cubic wave equation in dimension five,  this can be done via 
phase portrait analysis after the equation is written as an autonomous ODE system. This way, we obtain a radial $C^2$ solution of the elliptic equation
\begin{align*}
- \partial_{rr} \varphi - \frac{d-1}{r} \partial_r \varphi = | \varphi |^{p-1} \varphi \comma r>0
\end{align*}
which fails to belong in the critical space $\dot{H}^{s_p} \times \dot{H}^{s_p -1} (\mathbb{R}^d)$. 

Finally, by applying the channel of energy method, we prove the main rigidity result: Let $\vec{u} (t)$ be as in Proposition~\ref{Proposition1}. Then, $u_0 (r)$
must coincide with a singular stationary solution, whose construction is outlined above. This produces a contradiction because such stationary solutions do not lie in 
$\dot{H}^{s_p} \times \dot{H}^{s_p -1} (\mathbb{R}^d)$. 

%n%\newpage
\section{The Cauchy problem}
\label{sec:A review of the Cauchy problem}
In order to study the global dynamics of solutions to the Cauchy problem \eqref{wave equations}, we must first establish 
a local well-posedness theory. To that end, we review the Strichartz estimates from \cite{GV95, T10} and develop 
the theory of the Cauchy problem for the nonlinear wave equation.   

In this section, we sketch the theory of the local Cauchy problem for the nonlinear wave equation 
\begin{align}
   \begin{split}
    & \partial^2_{t} u
    - \Delta u 
    -  |u|^{p-1} u
     = 0   \comma \mbox{in } \mathbb{R}^d \times I ,
    \\
   &  \vec{u} (0)
    = (u_0, u_1) \in \dot{H}^{s_p} \times \dot{H}^{s_p -1} (\mathbb{R}^d)
    \end{split}
   \label{nonlinear wave equations}
\end{align}
 where the exponent $p$ obeys the condition
 \begin{equation}
p >  \frac{d+1}{2}  \qquad \mbox{ or } \qquad  p = 2m+1,~  m \in \mathbb{N} 
\label{p-condition}
\end{equation}
and $s_p = d/2 - 2 / (p-1)$.

%First, we recall 
%the Strichartz estimates for the linear wave equation in $\mathbb{R}^d \times I$
We begin with the linear wave equation in $\mathbb{R}^d \times I$
\begin{align}
   \begin{split}
    & \partial^2_{t} \omega
    - \Delta \omega 
     = h   \comma % \mbox{in } \mathbb{R}^d \times I ,
    \\
   &  \vec{\omega} (0)
    = (\omega_0, \omega_1) \in \dot{H}^{s} \times \dot{H}^{s -1} (\mathbb{R}^d) .
    \end{split}
   \label{linear wave equations}
  \end{align}
The solution operator to \eqref{linear wave equations} is given by
\begin{align}
\omega (t) = S (t) (\omega_0, \omega_1) + \int_0^{t} \frac{\sin ((t-s) \sqrt{- \Delta})}{\sqrt{- \Delta}} h(s)~ds
\label{eq1.2}
\end{align} 
 where
\begin{align}
 S(t) (\omega_0, \omega_1) = \cos (t \sqrt{- \Delta}) ~\omega_0 + (- \Delta)^{-1/2} \sin (t \sqrt{- \Delta})~\omega_1 .
 \label{eq1.3}
\end{align}
The operators in  \eqref{eq1.2}--\eqref{eq1.3} are defined via the Fourier transform: namely, we have
\begin{equation}
 \mathcal{F} (\cos (t \sqrt{- \Delta}) ~ f ) (\xi )  = \cos ( t |\xi|) \mathcal{F}( f ) (\xi) 
 \label{eq1.4}
\end{equation}
and
\begin{equation}
\mathcal{F} ( (- \Delta)^{-1/2} \sin (t \sqrt{- \Delta})~ f ) (\xi ) = | \xi |^{-1/2} \sin (t |\xi| ) \mathcal{F}( f ) (\xi).
\label{eq1.5}
\end{equation}
Similarly, we define the fractional differentiation operators as
\begin{equation}
\mathcal{F} (D^{\alpha} f) (\xi) = |\xi |^{\alpha} \mathcal{F} (f) (\xi). 
\label{eq1.6}
\end{equation}

Next, we recall the definition of \textit{Littlewood-Paley multipliers}. Let $\varphi (\xi)$ be a radial bump function supported in the ball 
$\{ \xi \in \mathbb{R}^d: |\xi| \leq 2\}$
with 
$
\varphi (\xi) = 1$ on $\{\xi \in \mathbb{R}^d: |\xi| \leq 1\}$.

We say that $N$ is a \textit{dyadic number} if $N = 2^k$ where $k \in \mathbb{Z}$ is an integer. For each dyadic integer $N$, we define
\begin{align}
& \mathcal{F} (P_{\leq N} f) (\xi)  = \varphi (\xi/N) \mathcal{F} (f) (\xi) \\
& \mathcal{F} (P_{> N} f) (\xi)  = (1- \varphi (\xi/N)) \mathcal{F} (f) (\xi) \\
& \mathcal{F} (P_N f) (\xi)  = (\varphi (\xi/N) - \varphi (2 \xi /N)) \mathcal{F} (f) (\xi). \label{eq1.7p}
\end{align} 
It is common to use the alternate notation $P_k$ to denote the Litlewood-Paley multiplier corresponding to the frequency $2^k$ (cf. Section~\ref{sec:Decay}).

We now recall the definition of the homogeneous Besov spaces $\dot{B}^s_{r, q}$. Let $s \in \mathbb{R}$ and $1\leq r \leq \infty$, 
$1 \leq q \leq \infty$. We define
\begin{equation}
\norm{u}_{\dot{B}^s_{r, q}} = \left( \sum_{N \in 2^{\mathbb{Z}}} (N^s \norm{ P_N u}_{L^r (\mathbb{R}^d)} )^q   \right)^{\frac{1}{q}}
\label{EQ_Besov1}
\end{equation}
and 
\begin{equation*}
\dot{B}^s_{r, q} (\mathbb{R}^d ) = \{ u \in \mathcal{S}' (\mathbb{R}^d): \norm{u}_{\dot{B}^s_{r, q}} < \infty \} .
\end{equation*}

\subsection{Strichartz Estimates}
Below, we state the Strichartz estimates for the wave equation \eqref{linear wave equations}. The result we present below may be found in  \cite{GV95, T10}. 

In what follows we say that a triple $(q, r, \gamma)$ is admissible if 
\begin{align*}
q, r \geq 2 \comma \frac{1}{q}  \leq \frac{d-1}{2} \left( \frac{1}{2} - \frac{1}{r} \right)   \comma \frac{1}{q}+ \frac{d}{r} = \frac{d}{2} - \gamma .
\end{align*}  
Denote by $q'$ and $r'$ the conjugate indices to $q$ and $r$, respectively, i.e., we have
\begin{align*}
\frac{1}{q} + \frac{1}{q'} = 1 \comma \frac{1}{r} + \frac{1}{r'} =1. 
\end{align*}

\begin{Lemma}[Strichartz Estimates, \cite{T10}]
\label{Lemma2.1}
Let $(q, r, \gamma)$ and $(a,b, \rho)$ be admissible triples.
% with $r< \infty$ and  $b < \infty$. 
Then, any solution $\omega$ to the linear Cauchy problem \eqref{linear wave equations} 
with initial data $\vec{\omega} (0) \in \dot{H}^{s} \times \dot{H}^{s -1}$
satisfies
\begin{align}
\| \vec{\omega} (t)  \|_{L^{q}_t \left( I; \dot{B}^{s -\gamma}_{r,2} \times  \dot{B}^{s -\gamma -1}_{r,2} \right) }
 \les \| \vec{\omega} (0) \|_{\dot{H}^{s} \times \dot{H}^{s -1}}
 + \| h \|_{L^{a'}_t \left( I; \dot{B}^{s -1 + \rho}_{b',2} \right)}
 \period
 \label{eq2.2}
\end{align}
\end{Lemma}

\begin{Remark}
{ \rm
Lemma~\ref{Lemma2.1} may also be stated in homogeneous Sobolev spaces, in which case we need the additional  requirement $r, b < \infty$. 
We refer the reader to \cite[Corollary~8.3]{T10} for further details. }
\end{Remark}

\subsection{Function Spaces} We introduce the following spaces to establish the local well-posedness of the Cauchy problem \eqref{nonlinear wave equations}. 
Let $I \in \mathbb{R}$ denote any time interval. We fix an exponent $p$ that satisfies \eqref{p-condition}. For dimensions $d \geq 7$, we define
\begin{align}
\begin{split}
& \norm{u}_{S_p(I)} = \norm{u}_{L^{2(p-1)}_t L^{\frac{2d}{3} (p-1)}_x  (I \times \mathbb{R}^d )} \\
& \norm{u}_{\dot{S} (I)} = \sup \{ \norm{u}_{L^q_t \dot{B}^{s_p - \gamma}_{r,2} (I \times \mathbb{R}^d)}, \norm{ \partial_t u}_{L^q_t \dot{B}^{s_p -1 -\gamma}_{r,2} (
I\times \mathbb{R}^d)}: (q,r, \gamma)~ \mbox{is an admissible triple } \} 
\\
& \norm{u}_{W(I)} = \norm{u}_{L^2_t \dot{B}^{s_p -1}_{\frac{2d}{d-3},2} (I\times \mathbb{R}^d) } \\
& \norm{u}_{W' (I)} = \norm{u}_{L^1_t \dot{H}^{s_p -1} (I\times \mathbb{R}^d) } \\
& \norm{u}_{X(I)} = \norm{u}_{L^{q_x}_t \dot{B}^{s_p - \frac{1}{2}}_{r_x,2} (I\times \mathbb{R}^d)} \comma q_x = 2(p-1), r_x = \frac{2d (p-1)}{d(p-1) -p} \\
& \norm{u}_{X' (I)} = \norm{u}_{L^{\tilde{q}'}_t \dot{B}^{s_p -\frac{1}{2}}_{\tilde{r}', 2} (I\times \mathbb{R}^d) } \comma \tilde{q} = \frac{2(p-1)}{(p-2)}, \tilde{r}= \frac{2d (p-1)}{d (p-1) -2p+3} 
\period
\end{split}
\label{norms}
\end{align}

Applying \eqref{eq2.2} with $(a,b, \rho) = (\infty, 2, 0)$, $(q,r, \gamma)= (2 (p-1), \frac{2d}{3} (p-1), s_p )$, and any other admissible triple, which 
yields terms that may be absorbed into $ \Vert \omega \Vert_{ \dot{S} (I)}$,
we obtain the following Strichartz estimate
 for solutions to the linear Cauchy problem \eqref{linear wave equations}:
\begin{align}
\begin{split}
&
\sup_{t \in \mathbb{R}} \Vert  \vec{\omega} (t) \Vert_{\dot{H}^{s_p} \times \dot{H}^{s_p -1}}
+ 
\Vert \omega \Vert_{S_p (I)} 
+ \Vert \omega \Vert_{ \dot{S} (I)} 
\\
& \indeq \indeq \indeq \indeq  \les
\Vert \vec{\omega} (0) \Vert_{\dot{H}^{s_p} \times \dot{H}^{s_p -1}}
+
\Vert  h \Vert_{W' (I)}.
\end{split}
\label{eq2.3}
\end{align}

\begin{Remark}
{ \em 
Note that the right hand side of \eqref{eq2.3} also controls  $\Vert \omega \Vert_{X(I)}$ and $\Vert \omega \Vert_{W(I)}$ as these norms are also admissible.   
}
\end{Remark}

\begin{Lemma}
\label{Interpolations}
Let $d \geq 7$. Then, we have the following inequalities.
\begin{align}
& \norm{u}_{X(I)} \les \norm{u}_{\dot{S}(I)}  \label{embd1} \\
& \norm{u}_{W(I)} \les \norm{u}_{\dot{S}(I)} \label{embd2} \\
& \norm{u}_{S_p (I)} \les \norm{u}_{X(I)}   \label{embd3} \\
& \norm{u}_{X(I)} \les \norm{u}^{\theta}_{S_p (I)} \norm{u}_{\dot{S}(I)}^{1- \theta} \quad \quad \mbox{ for some } \theta \in (0,1). \label{embd4}
\end{align} 
\end{Lemma}

\begin{proof}
Inequalities \eqref{embd1} and \eqref{embd2} hold due to the admissability of the $X(I)$ and $W(I)$ norms. We combine the Sobolev embedding with the Besov-Sobolev embedding to obtain \eqref{embd3}. Finally, \eqref{embd4} holds by interpolation. 
\end{proof}

\begin{Lemma}[Inhomogeneous Strichartz Estimates]
Let $d \geq 7$ and $I \in \mathbb{R}$ denote any time interval. Then, we have
\begin{align}
\norm{ \int_0^t \frac{\sin ( (t- \tau) \sqrt{- \Delta}) }{\sqrt{- \Delta}} f(\tau) ~d \tau }_{X(I)} \les \norm{f}_{X'(I)} . \label{Inhom SE}
\end{align}
\end{Lemma}

\begin{proof}
Estimate \eqref{Inhom SE} is an application of Corollary~$8.7$ from \cite{T10}. We appeal to Corollary~$8.7$ with the exponent pairs  $(q_x, r_1)$ 
and $(\tilde{q}, \tilde{r}_1)$, where
\begin{align}
r_1 = \frac{2d}{d-1} \comma \tilde{r}_1 = \frac{2d (d-1)}{d^2 - 2d -1}.
\end{align}
Note that $r_1, \tilde{r}_1 \in [2, \infty)$ satisfying
\begin{align*}
\frac{1}{q_x}  < (d-1) \left( \frac{1}{2} - \frac{1}{r_1}  \right) \comma
\frac{1}{\tilde{q}}  < (d-1) \left(  \frac{1}{2} - \frac{1}{\tilde{r}_1} \right)
\end{align*}
as well as 
\begin{equation*}
\frac{1}{q_x} + \frac{1}{\tilde{q}} = \frac{d-1}{2} \left(  1- \frac{1}{r_1} - \frac{1}{\tilde{r}_1}  \right) 
\end{equation*}
and 
\begin{equation*}
\frac{d-3}{r_1} \leq \frac{d-1}{\tilde{r}_1} \comma \frac{d-3}{\tilde{r}_1} \leq \frac{d-1}{r_1}. 
\end{equation*}
Since $r_x \geq r_1$ and $\tilde{r} \geq \tilde{r}_1$, yielding
\begin{equation*}
s_p - \frac{1}{2} + d \left( \frac{1}{2} - \frac{1}{r_x} \right) - \frac{1}{q_x} = 1- \left( - s_p + \frac{1}{2} + d \left( \frac{1}{2}- \frac{1}{\tilde{r}} \right) - \frac{1}{\tilde{q}} \right) 
\end{equation*}
we may apply Corollary~$8.7$ to obtain \eqref{Inhom SE}.
\end{proof}
\subsection{Nonlinear Estimates}
Having defined the function spaces in which the Cauchy theory for the nonlinear wave equation will be established, we concentrate on the nonlinear estimates. 
First, we state a Leibniz rule and a version of $C^1$ chain rule, both of which are discussed in Besov spaces. 
\begin{Lemma}[Leibniz Rule] \label{Leibniz Rule}
Let $s \geq 0$ and $r \in [1, \infty)$. Then, 
\begin{equation}
\norm{fg}_{\dot{B}^{s}_{r,2}} \les \norm{f}_{\dot{B}^s_{p_1, 2}} \norm{g}_{L^{p_2}} + \norm{f}_{L^{q_1}} \norm{g}_{\dot{B}^s_{q_2, 2}}
\end{equation}
where $p_{1}, q_{2} \in [1, \infty)$, $p_2, q_1 \in [1, \infty]$, and
\begin{equation*}
\frac{1}{p_1} + \frac{1}{p_2} = \frac{1}{q_1} + \frac{1}{q_2} = \frac{1}{r}.
\end{equation*}
\end{Lemma}

\begin{Lemma}[$C^1$ Chain Rule] \label{Chain Rule}
Let $s \geq 0$ and $k= 1+ \left[  s \right]$. Assume that $G \in C^k (\mathbb{R})$ is a function of the form $G(x)= |x|^q$ or $G(x) = |x|^{q-1} x$, and let  $u \in \dot{B}^s_{{p}_2,2} \cap L^{p_1}$. Then, $G(u) \in \dot{B}^s_{r,2}$ with 
\begin{equation}
\norm{ G(u)}_{\dot{B}^s_{r, 2}} \les \norm{u}^{q-1}_{L^{p_1}} \norm{u}_{\dot{B}^s_{p_2, 2}}
\end{equation}
where
$1 / r = (q-1)/p_1 + 1/p_2$.
\end{Lemma}
% k= \left \lceil{s}\right \rceil
Note that if $G (x) = |x|^q$, then the hypothesis $G \in C^k (\mathbb{R})$ requires $q$ to be either an even integer or any number with $q>k$ . 
 In the case $G$ takes the latter form, we need to assume that $q$ is either an odd integer or that $q >k$. 
The proofs of  Lemma~\ref{Leibniz Rule}--\ref{Chain Rule} are provided in the appendix. 
We now apply these results to prove that the nonlinearity of \eqref{nonlinear wave equations} may be estimated in 
Strichartz spaces of Besov type. 
\begin{Lemma}\label{Nonlinear Estimate 1}
Let  $p -1 > 1+\left[ {s_p -1}\right]$ or let $p$ be an odd %even 
integer. Then, we have
\begin{align}
\norm{ F(u)}_{W' (I)} \les \norm{u}^{p-1}_{S_p (I)} \norm{u}_{W (I)}
\label{EQ_NE1}
\end{align}
and 
\begin{align}
\begin{split}
\norm{F(u) - F(v)}_{W' (I)} & \les \norm{u-v}_{W(I)} \left( \norm{u}^{p-1}_{S_p (I)} + \norm{v}^{p-1}_{S_p (I)} \right) \\
& \indeq + \norm{u-v}_{S_p (I)} \left( \norm{u}^{p-2}_{S_p (I)} + \norm{v}^{p-2}_{S_p (I)} \right) \\
& \indeq \indeq \indeqtimes \left( \norm{u}_{W(I)} + \norm{v}_{W(I)} \right) \period
\end{split}
\label{EQ_NE2}
\end{align}
\end{Lemma}

\begin{proof}
We express 
\begin{equation}
F(u) - F(v) = (u-v) \int_0^1 F' ( \lambda u + (1 - \lambda v) )~d \lambda
\end{equation}
by the fundamental theorem of calculus. Then, by Lemma~\ref{Leibniz Rule} and H\"older's inequality we have
\begin{align}
\begin{split}
\Vert   F(u) - F(v) \Vert_{W' (I)} & \les 
\Vert u-v \Vert_{ W(I)} \norm{ \int_0^1 F' ( \lambda u + (1 - \lambda v) )~d \lambda }_{L^2_t L^{\frac{2d}{3}}_x} \\
& \quad + \Vert u-v \Vert_{S_p (I)} 
\norm{   \int_0^1 F' ( \lambda u + (1 - \lambda v) )~d \lambda }_{ L^{q_{p}}_t  \dot{B}^{s_p -1}_{r_{p}, 2}}
\end{split}
\label{EQ2.25}
\end{align}
where $q_p = \frac{2 (p-1)}{p-3}$ and $r_p = \frac{2d (p-1)}{d (p-1)-3}$. Also, by Minkowski's inequality we get
\begin{align}
\begin{split}
\norm{ \int_0^1 F' ( \lambda u + (1 - \lambda v) )~d \lambda }_{L^2_t L^{\frac{2d}{3}}_x} 
& \les \sup_{\lambda \in [0,1]}  \Vert F' (\lambda u + (1- \lambda) v) \Vert_{L^2_t L^{2d/3}_x} \\
& \les \Vert u \Vert_{S_p (I)}^{p-1} + \Vert v \Vert_{S_p (I)}^{p-1} 
\period
\end{split}
\label{EQ2.22}
\end{align}
Similarly, we find 
\begin{equation}
\norm{   \int_0^1 F' ( \lambda u + (1 - \lambda v) )~d \lambda}_{ L^{q_{p}}_t  \dot{B}^{s_p -1}_{r_{p},2}}
\les \sup_{\lambda \in [0,1]} \norm{F' (\lambda u + (1-\lambda) v) }_{ L^{q_{p}}_t  \dot{B}^{s_p -1}_{r_{p},2}}
\period 
\label{EQ2.23}
\end{equation} 
Next, we apply Lemma~\ref{Chain Rule} to estimate the right hand side of \eqref{EQ2.23}. We get
\begin{equation*}
\norm{|u|^{p-1}}_{\dot{B}^{s_p -1}_{r_p, 2} } \les \norm{u}^{p-2}_{\frac{2d}{3} (p-1)} \norm{u}_{\dot{B}^{s_p -1}_{\frac{2d}{d-3}, 2}}
\end{equation*}
which results in
\begin{align}
\sup_{\lambda \in [0,1]} \norm{F' (\lambda u + (1-\lambda) v) }_{ L^{q_{p}}_t  \dot{B}^{s_p -1}_{r_{p},2}} 
\les \left( \Vert u \Vert^{p-2}_{S_p (I)} + \Vert v \Vert^{p-2}_{S_p (I)} \right)
\left( \Vert u \Vert_{W (I)} + \Vert v \Vert_{W(I)} \right) . 
\label{EQ2.24}
\end{align}
Collecting the estimates in \eqref{EQ2.25}, \eqref{EQ2.22}, and \eqref{EQ2.24}, we obtain
\eqref{EQ_NE2}. Also, setting $v=0$ in \eqref{EQ_NE2} we get \eqref{EQ_NE1}.
\end{proof}

Allowing a change in H\"older exponents, we may obtain the nonlinear estimates in $X'(I)$--$X(I)$
norms in the same manner. 
\begin{Lemma} \label{Nonlinear Estimate 2}
Let $p -1 > 1+ \left[ {s_p -\frac{1}{2}} \right]$ or let $p$ be an odd integer. Then, we have
\begin{align}
\norm{ F(u)}_{X' (I)} \les \norm{u}^{p-1}_{S_p (I)} \norm{u}_{X (I)}
\label{EQ_NE3}
\end{align}
and 
\begin{align}
\begin{split}
\norm{F(u) - F(v)}_{X' (I)} & \les \norm{u-v}_{X(I)} \left( \norm{u}^{p-1}_{S_p (I)} + \norm{v}^{p-1}_{S_p (I)} \right) \\
& \indeq + \norm{u-v}_{S_p (I)} \left( \norm{u}^{p-2}_{S_p (I)} + \norm{v}^{p-2}_{S_p (I)} \right) \\
& \indeq \indeq \indeqtimes \left( \norm{u}_{X (I)} + \norm{v}_{X (I)} \right) \period
\end{split}
\label{EQ_NE4}
\end{align}
\end{Lemma}
\begin{Remark} 
\label{R2.5}
{ \rm
Both Lemma~\ref{Nonlinear Estimate 1} and Lemma~\ref{Nonlinear Estimate 2} play essential roles in establishing the local well-posedness and perturbation results. Thus, 
we will need the exponent $p$ to obey the hypothesis of Lemma~\ref{Nonlinear Estimate 2}, which is slightly more restrictive than that of Lemma~\ref{Nonlinear Estimate 1}.
Namely, as $(d-1)/2$ is an integer, this hypothesis takes the form of \eqref{p-condition}.
For the remainder of Section~\ref{sec:A review of the Cauchy problem}, we assume that $p$ satisfies \eqref{p-condition}.
}
\end{Remark}

\subsection{Local Well-Posedness}
Using the estimates in Lemma~\ref{Nonlinear Estimate 1} and Lemma~\ref{Nonlinear Estimate 2}, we obtain the following local well-posedness result. 
First, we recall the definition of a solution to the Cauchy problem \eqref{nonlinear wave equations}. 

\begin{definition}\label{def1}
We say that $u$ is a \textit{strong solution} to \eqref{nonlinear wave equations} on a time interval $I$ with $t_0 \in I$ if
$(u, \partial_t u) \in C(I; \dot{H}^{s_p} \times \dot{H}^{s_p -1})$, $u \in {S_p}(I) \cap W(I)$, $(u, \partial_t u) = (u_0, u_1)$
and the integral equation 
\begin{equation*}
u (t) = S(t- t_0) ((u_0, u_1)) - \int_{t_0}^{t} \frac{\sin ((t-s)\sqrt{- \Delta})}{\sqrt{- \Delta}} F(u(s))~ds \comma t \in I
\end{equation*}
holds with $F(u)= |u|^{p-1} u$. 
\end{definition}

\begin{Theorem}
\label{Thm2.2}
Let $d \geq 7$, $(u_0, u_1) \in \dot{H}^{s_p} \times \dot{H}^{s_p -1}$, and $p$ satisfy the condition \eqref{p-condition}.
Let $I \in \mathbb{R}$  be an interval containing $0 \in I^{\circ}$. 
 Assume that  
\begin{equation}
\Vert (u_0, u_1) \Vert_{ \dot{H}^{s_p} \times \dot{H}^{s_p -1}} \leq A.
\end{equation}
Then, there exists $\delta_0 = \delta_0 (d, p, A) > 0$ such that if 
\begin{align*}
\Vert S(t) (u_0, u_1) \Vert_{S_p (I)} = \delta < \delta_0,
\end{align*}
the Cauchy problem \eqref{nonlinear wave equations} admits a unique solution $u$  in $\mathbb{R}^d \times I$ 
as defined in Definition~\ref{def1}. Moreover, we have
\begin{align*}
\Vert u \Vert_{S_p (I)} & < 2 \delta \\
\Vert u \Vert_{\dot{S} (I)} &< \infty \\
\Vert u \Vert_{X(I)}&  \leq C_0 \delta^{\theta} A^{1 - \theta}
\end{align*}
for some $\theta \in (0,1)$. 
Furthermore, if $(u_{0,k}, u_{1,k}) \to (u_0, u_1)$  as $k \to \infty$ in $\dot{H}^{s_p} \times \dot{H}^{s_p -1}$, 
then 
\begin{align*}
(u_k, \partial_t u_k) \to (u, \partial_t u) \quad \quad \mbox{in~} C(I; \dot{H}^{s_p} \times \dot{H}^{s_p -1}) 
\end{align*}
where $u_k$ is the solution corresponding to $(u_{0,k}, u_{1,k})$.
%maybe, move this continuity claim into another Remark, and cite Remark~2.17 in KM1. 
\end{Theorem}

We note that the proof of Theorem~\ref{Thm2.2} is very similar to the proof of local well-posedness results 
presented in \cite{BCLPZ13} and  \cite[Theorem~2.7]{KM2} (see also \cite{KM1, DL15}). Below, we give an outline of the ideas. 
\begin{proof}
We begin with the size of the linear part $S(t) (u_0, u_1) $ under the hypothesis: by the Strichartz inequality
\begin{align}
\norm{ S(t) (u_0, u_1)}_{\dot{S}(I)} \leq C A. 
\label{EQ2.1}
\end{align}
Also, we obtain 
\begin{align}
\norm{ S(t) (u_0, u_1)}_{X(I)} & \leq C \norm{ S(t) (u_0, u_1)}^{\theta}_{S_p (I)} \norm{S(t) (u_0, u_1)}^{1 - \theta}_{\dot{S} (I)} \\
& \leq C \delta^{\theta} A^{1 - \theta}
\label{EQ2.3}
\end{align}
by \eqref{embd4}. 

Next, we set $\eta := C \delta^{\theta} A^{1 - \theta}$, and define
\begin{align*}
u^0 & = S(t) (u_0, u_1), \\
u^{n+1} & = S(t) (u_0, u_1) + \int_0^{t} \frac{\sin ((t-s) \sqrt{- \Delta})}{\sqrt{- \Delta}} F(u^n) (s)~ds.  
\end{align*}
We then show that the sequence $\{ u^n\}$ is bounded in $\dot{S} (I)$ and in $X(I)$.  More precisely, setting 
\begin{align*}
a := 2 \eta ~\mbox{ and }~ b:= 2 CA
\end{align*}
we will show by induction that 
\begin{align}
\norm{u^n}_{S_p (I)} \leq 2 \delta \comma
\norm{u^n}_{X(I)} \leq a \comma  \norm{ u^n}_{\dot{S}(I)} \leq b. 
\label{EQ2.2}
\end{align}
Note that \eqref{EQ2.2} is satisfied for $n=0$ by \eqref{EQ2.1} and \eqref{EQ2.3}. Furthermore, by 
Strichartz inequality and \eqref{EQ_NE1}, we obtain
\begin{align*}
\norm{u^{n+1}}_{\dot{S} (I)} &\leq CA + C \norm{F(u^n)}_{W' (I)} \\
& \leq C A + C \norm{u^n}^{p-1}_{S_p (I)} \norm{u^n}_{W(I)} \\
& \leq C A + C \delta^{p-1} b. 
\end{align*}
In the same fashion, our hypothesis on $\norm{S(t) (u_0, u_1)}_{S_p (I)}$ along with \eqref{embd3} implies that 
\begin{align*}
\norm{u^{n+1}}_{{S}_p (I)} &\leq \delta + C \norm{F(u^n)}_{W' (I)} \\
& \leq \delta + C \norm{u^n}^{p-1}_{S_p (I)} \norm{u^n}_{W(I)} \\
& \leq \delta + C \delta^{p-1} b. 
\end{align*}
Since $p > 2$, letting $\delta_0$ to be sufficiently small we guarantee that
\begin{equation}
\norm{u^{n+1}}_{S_p (I)} \leq 2 \delta \comma
\norm{u^{n+1}}_{\dot{S} (I)} \leq b. 
\end{equation}
Next, we use \eqref{EQ_NE3} and observe that 
\begin{align*}
\norm{u^{n+1}}_{X (I)} &\leq \eta + C   \norm{u^n}_{X(I)}  \left( \norm{u^n}^{p-1}_{S_p (I)} + \norm{u^n}^{p-2}_{S_p (I)} \norm{u^n}_{X(I)} \right) \\
& \leq \frac{a}{2} + \frac{C a^{p}}{2^p}  \\
& \leq  a
\end{align*}
assuming that $a$ is sufficiently small. Thus, by induction \eqref{EQ2.2} holds. 

Next, we show that the sequence $\{u^n\}$ is Cauchy in $X(I)$. Applying Lemma~\ref{Inhom SE} and  the nonlinear difference estimate \eqref{EQ_NE4} we get
\begin{align*}
\norm{u^{n+1} - u^n}_{X(I)} & \les \norm{F (u^{n}) - F(u^{n-1}) }_{X' (I)} \\
& \les \norm{u^n - u^{n-1}}_{X(I)} \left( \norm{u^n}^{p-1}_{S_p (I)} + \norm{u^{n-1}}^{p-1}_{S_p (I)} \right) \\
& \indeq + \norm{u^n - u^{n-1}}_{S_p (I)} \left( \norm{u^n}^{p-2}_{S_p (I)} + \norm{u^{n-1}}^{p-2}_{S_p (I)} \right) \\
& \indeq \indeq \indeqtimes \left( \norm{u^n}_{X (I)} + \norm{u^{n-1}}_{X (I)} \right) \period
\end{align*}
By \eqref{embd3} and \eqref{EQ2.2}, we then note that the right hand side above may be estimated by
\begin{align*}
C_0 \norm{u^n - u^{n-1}}_{X(I)} \delta^{p-2} a . 
\end{align*}
Allowing $\delta_0$ to be sufficiently small guarantees that the sequence converges to some $u$ in $X(I)$. 
The rest of the arguments follow from the ideas outlined in \cite[Theorem~3.3]{BCLPZ13}. Namely, 
$\{ u^n\}$ is a bounded sequence in $\dot{S} (I)$, which in particular  implies that it is bounded in $W(I)$. Since 
$W(I)$ is a reflexive Banach space, $u$ is the weak limit of $\{ u^n\}$ in $W(I)$. By the Strichartz inequality, we then deduce 
that $u \in \dot{S}(I)$. The fact that $u$ is a unique solution of \eqref{nonlinear wave equations} in $\dot{S}(I)$ is also a consequence of the 
Strichartz estimates, coupled with the nonlinear estimates above. 

The last continuity statement also follows from an application of the same arguments as done above. Since $(u_{0,k}, u_{1,k}) \to (u_0, u_1)$   in $\dot{H}^{s_p} \times \dot{H}^{s_p -1}$ as $k \to \infty$, we have for large $k$, $\norm{ S(t) (u_{0,k}, u_{1,k})}_{S_p (I)} \leq \delta_0$. We may then repeat the estimates from above to show that
\begin{align*}
(u_k, \partial_t u_k) \to (u, \partial_t u) \quad \quad \mbox{in~} C(I; \dot{H}^{s_p} \times \dot{H}^{s_p -1}).
\end{align*}
\end{proof}

\begin{Remark}
\label{R2.3}
{\rm
As noted in analogous results in \cite{KM11, KM2, KM1}, the proof of Theorem~\ref{Thm2.2} implies that there exists $\tilde{\delta} >0$ such that 
if $\Vert (u_0, u_1) \Vert_{\dot{H}^{s_p} \times \dot{H}^{s_p -1}} < \tilde{\delta}$, the above result holds with $I = \mathbb{R}$.  
Moreover, given  $(u_0, u_1) \in \dot{H}^{s_p} \times \dot{H}^{s_p -1}$, there exists a time interval $I \in \mathbb{R}$ containing $0$ such that the hypothesis of Theorem~\ref{Thm2.2} is justified. In fact, following the arguments as in  \cite[Definition~2.10]{KM1}, we deduce that there exists a maximal time interval $I ((u_0, u_1))$ with $0 \in I$, where the solution $\vec{u} (t) = (u (t), \partial u (t))$ is defined in the sense of Definition~\ref{def1}. We will denote the maximal time interval of a solution $\vec{u}$ by
\begin{equation}
\Imax (\vec{u}) = (T_- (\vec{u}), T_+ (\vec{u}))
\label{Imax}
\end{equation}
}
\end{Remark}

As a consequence of the proof of Theorem~\ref{Thm2.2}, we may obtain a finite time blow-up criterion and a scattering result in a standard manner.
We state  these results below. For analogous proofs, we refer the reader to see \cite[Lemma~2.11 and Remark~2.15, resp.]{KM1}. 
\begin{Remark}[Global Existence and Scattering]
\label{R02}
{\rm
Let $\vec{u} (t)$ be a solution of \eqref{nonlinear wave equations} in $(T_- (\vec{u}), T_+ (\vec{u}))$. 
If $T_+  (\vec{u})< \infty$, then we have
\begin{align*}
\Vert u \Vert_{S_p ([0, T_+ ))} = \infty . 
\end{align*}
Along with the statement above, we recall the equivalence between scattering and boundedness of $S_p$ norms. More precisely, we have   
$\Vert u \Vert_{S_p ([0, T_+ (\vec{u})))} < \infty$ if and only if $\vec{u} (t)$
scatters to a free wave as $t \to \infty$, i.e.,  
there exists $(u_0^+, u_1^+) \in \dot{H}^{s_p} \times \dot{H}^{s_p -1}$  so that  
\begin{align*}
\lim_{t \to \infty} \Vert (\vec{u} (t) - \vec{S} (t) (u_0^+, u_1^+)) \Vert_{\dot{H}^{s_p} \times \dot{H}^{s_p -1}} =0 .
\end{align*}
The same equivalence also holds for solutions $\vec{u} (t)$ on $(T_- (\vec{u}), 0]$. 
A finite time blow-up criterion may be stated for $T_{-} (\vec{u}) > - \infty$ as well.  }
\end{Remark}

Next, we prove a long-time perturbation theorem for approximate solutions to \eqref{nonlinear wave equations}. Our first long-time perturbation result is established in Besov spaces. 

\begin{Theorem}[Long Time Perturbation in Besov Spaces]
\label{perturbation_thm}
Let  $(u_0, u_1) \in \dot{H}^{s_p} \times \dot{H}^{s_p -1}$ and $I \subset \mathbb{R}$ be an open interval containing $t_0$. 
Assume that $v$ solves the equation
\begin{align*}
\partial^2_{t} v
    - \Delta v 
     =    |v|^{p-1} v + e  
\end{align*}
in the sense of the corresponding  integral equation, and it satisfies
\begin{itemize}
\item[(i)]
$\sup_{t \in I} \Vert v \Vert_{\dot{H}^{s_p} \times \dot{H}^{s_p -1}} \leq A $
\item[(ii)]
$\Vert v \Vert_{S_{p} (I)} \leq A$ and $\Vert  v \Vert_{W(I)} < \infty$.  
\end{itemize}
Additionally, assume that we have
\begin{align}
\Vert (u_0 - v(t_0), u_1 - \partial_t v (t_0) ) \Vert_{\dot{H}^{s_p} \times \dot{H}^{s_p -1}} \leq A'
\label{EQ2.4p}
\end{align}
as well as the smallness condition 
\begin{align}
\Vert D^{s_p -1} e \Vert_{L^{1}_t L^{2}_x} + \Vert S(t- t_0) (u_0 - v(t_0), u_1 - \partial_t v (t_0) ) \Vert_{X (I)} 
\leq \epsilon .
\label{EQ2.5p}
\end{align}
Then there exists $\epsilon_0 = \epsilon_0 (d, p, A, A') >0$  such that if  $0 < \epsilon < \epsilon_0$,  there is a unique solution $u$ of \eqref{nonlinear wave equations}
 in $I$ with $(u(t_0), \partial_t u (t_0))= (u_0, u_1)$ 
satisfying
 \begin{align*}
\Vert u-v  \Vert_{\dot{S}(I)} & \leq C(d, p, A, A') ,
\\
 \Vert u - v \Vert_{X (I)} & \leq  C(d, p, A, A') \cdot \epsilon .
 \end{align*}
 \end{Theorem}
 
 First, we prove the following short-time perturbation result proof of which follows from the arguments given in the analogous version in  \cite{BCLPZ13}. We point out the estimates that are different from those in \cite{BCLPZ13}.  
 \begin{Theorem}[Short Time Perturbation]
 \label{short time perturbation}
 Let  $(u_0, u_1) \in \dot{H}^{s_p} \times \dot{H}^{s_p -1}$ and $I \subset \mathbb{R}$ be an open interval containing $t_0$. 
Assume that $v$ solves the equation
\begin{align*}
\partial^2_{t} v
    - \Delta v 
     =    |v|^{p-1} v + e  
\end{align*}
in the sense of the corresponding  integral equation such that
\begin{equation}
\Vert (u_0 - v(t_0), u_1 - \partial_t v (t_0) ) \Vert_{\dot{H}^{s_p} \times \dot{H}^{s_p -1}} \leq A' .
\label{EQ2.6}
\end{equation}
Furthermore, we assume that 
 \begin{align*}
 \Vert D^{s_p -1} e \Vert_{L^{1}_t L^{2}_x}  & \leq \epsilon , \\
  \Vert S(t- t_0) (u_0 - v(t_0), u_1 - \partial_t v (t_0) ) \Vert_{X (I)}  & \leq \epsilon , \\
  \Vert v \Vert_{X(I)} + \Vert v \Vert_{W(I)} & \leq \delta.
 \end{align*}
 Then, there exist $\delta_0 = \delta_0 (d,p) >0$, $\epsilon_0 = \epsilon_0 (A') >0$ such that if $0 < \delta < \delta_0$ and $0< \epsilon < \epsilon_0$,  there is a unique  solution $u$ of \eqref{nonlinear wave equations}
 in $I$ with $(u(t_0), \partial_t u (t_0))= (u_0, u_1)$, 
 which obeys
 \begin{align}
\Vert u-v  \Vert_{\dot{S}(I)} & \leq C(d, p, A')  , \label{EQ2.7}
\\
 \Vert u - v \Vert_{X (I)} & \leq  C(d, p) \cdot \epsilon , \label{EQ2.8}
 \\ 
 \Vert F(u) - F(v) \Vert_{X' (I)} & \leq C(d, p) \cdot \epsilon.  \label{EQ2.9}
 \end{align}
 \end{Theorem}
 
 \begin{proof}
 Without loss of generality, we take $t_0 =0$, and set $(v(t_0), \partial_t v (0)) = (v_0, v_1)$. Letting $\epsilon_0 >0$ to be sufficiently small, we may assume that $u$ exists on $I$, and that it satisfies 
\begin{equation*}
\partial^2_{t} (u - v)
    - \Delta (u- v)
     =   e + |u|^{p-1} u -  |v|^{p-1} v . 
\end{equation*} 
 Letting $F(x)= |x|^{p-1} x$, we observe that 
\begin{align*}
\begin{split}
\norm{u-v}_{X(I)} & \les \norm{S(t) (u_0 - v_0, u_1 - v_1  )}_{X(I)} 
+ 
\norm{\int_0^t \frac{\sin ((t - \tau) \sqrt{- \Delta})}{\sqrt{- \Delta}} e(\tau) ~ d\tau }_{X(I)} \\
& \indeq
+ \norm{\int_0^t \frac{\sin ((t - \tau) \sqrt{- \Delta})}{\sqrt{- \Delta}} (F(u) - F(v)) ~d\tau }_{X(I)} \\
& = I_1 + I_2 + I_3
\period
\end{split}
\end{align*} 
Our assumption on the inital data $(u_0 - v_0, u_1 - v_1)$ implies that 
$
I_1 \les \epsilon $.
For $I_2$, we use the Strichartz inequality and Lemma~\ref{Interpolations}, 
\begin{align*}
I_2 & \les 
\norm{\int_0^t \frac{\sin ((t - \tau) \sqrt{- \Delta})}{\sqrt{- \Delta}} e(\tau) ~ d\tau }_{\dot{S}(I)} \\
& \les
\norm{D^{s_p -1} e}_{L^1_t L^2_x} \\
& \les \epsilon. 
\end{align*}
Next, we estimate $I_3$. By Lemma~\ref{Inhom SE}, Lemma~\ref{Nonlinear Estimate 2}, and the embedding \eqref{embd3} we get
\begin{align*}
I_3 & \les \norm{ F(u) - F(v)}_{X' (I)} \\
& \les \norm{u-v}_{X(I)} \left( \norm{u-v}^{p-1}_{S_p (I)} + \norm{v}^{p-1}_{S_p (I)} \right) \\
& \indeq + \norm{u-v}_{S_p (I)} \left( \norm{u-v}^{p-2}_{S_p (I)} + \norm{v}^{p-2}_{S_p (I)} \right) \\
& \indeq \indeq \indeqtimes \left( \norm{u-v}_{X (I)} + \norm{v}_{X (I)} \right) \\
& \les \norm{u-v}_{X(I)} \left( \norm{u- v}^{p-1}_{X (I)} + \delta^{p-1}  \right) \period
\end{align*}
In the last step above we also used the Young's inequality to find a common upper bound 
\begin{equation}
 ( \norm{u- v}_{X(I)}^{p-2} + \delta^{p-2}) (  \norm{u- v}_{X(I)} + \delta ) \les  \norm{u- v}_{X(I)}^{p-1} + \delta^{p-1} .
\end{equation}
Next, we add up the estimates above to find
\begin{align*}
\norm{u-v}_{X(I)} \les \epsilon + \norm{u-v}_{X(I)} \left( \norm{u- v}^{p-1}_{X (I)} + \delta^{p-1}  \right). 
\end{align*}
Since the estimate for $\norm{u-v}_{X(I)}$ does not depend on the norm difference in the $\dot{S} (I)$ norm,  we may close this estimate 
independently 
by selecting $\delta \leq \delta_0 (d, p)$ to be sufficiently small. Note that the choice for $\delta_0$ only depends on the dimension $d$ and the exponent $p$.
As a result, we obtain \eqref{EQ2.8} and \eqref{EQ2.9} for small enough $\delta_0$ and $\epsilon_0$. 

Similarly, we may estimate $\norm{u-v}_{\dot{S} (I)}$. By the Strichartz inequality and Lemma~\ref{Nonlinear Estimate 1}, we obtain
\begin{align*}
\norm{u-v}_{\dot{S} (I)} & \les A' + \epsilon + \norm{F(u) - F(v)}_{W' (I)} \\
& \les A' + \epsilon + \norm{u-v}_{W(I)}  \left( \norm{u-v}^{p-1}_{S_p (I)} + \norm{v}^{p-1}_{S_p (I)} \right) \\
& \indeq + \norm{u-v}_{S_p (I)} \left( \norm{u-v}^{p-2}_{S_p (I)} + \norm{v}^{p-2}_{S_p (I)} \right) \\
& \indeq \indeq \indeqtimes \left( \norm{u-v}_{W (I)} + \norm{v}_{W (I)} \right)
\period
\end{align*}
We then combine our assumptions with the embeddings in Lemma~\ref{Interpolations} to have
\begin{align*}
\begin{split}
\norm{u-v}_{\dot{S} (I)} & \les A' + \epsilon + \norm{u-v}_{\dot{S} (I)} \left( \norm{u-v}_{X (I)}^{p-1} + \delta^{p-1} \right) \\
& \indeq + \norm{u-v}_{X (I)} \left( \norm{u-v}^{p-2}_{X (I)} + \delta^{p-2} \right) \\
& \indeq \indeq \indeqtimes \left( \norm{u-v}_{\dot{S} (I)} + \delta \right) \\
& \les A' + \epsilon + \norm{u-v}_{\dot{S} (I)} \left( \epsilon^{p-1} + \delta^{p-1} \right)
\period
\end{split}
\end{align*}
Absorbing the last term above on the left hand side and noting that $\epsilon < \epsilon_0 = \epsilon_0 (A')$, we obtain \eqref{EQ2.7}.
 \end{proof}
 
 \begin{proof}[Proof of Theorem~\ref{perturbation_thm}] 
 The same arguments used in the proof Theorem~$5.3$  in \cite{BCLPZ13} yield the 
 proof of the Theorem~\ref{perturbation_thm}. 
\end{proof}  

Next, we state the Sobolev version of our long time perturbation theorem. 
\begin{Theorem}[Long Time Perturbation in Sobolev Spaces]
\label{perturbation_thm_Sobolev}
Let  $(u_0, u_1) \in \dot{H}^{s_p} \times \dot{H}^{s_p -1}$ and $I \subset \mathbb{R}$ be an open interval containing $t_0$. 
Assume that $v$ solves the equation
\begin{align*}
\partial^2_{t} v
    - \Delta v 
     =    |v|^{p-1} v + e  
\end{align*} 
in the sense of the corresponding  integral equation, and it satisfies
\begin{equation}
\sup_{t \in I} \Vert v \Vert_{\dot{H}^{s_p} \times \dot{H}^{s_p -1}} 
+\Vert D^{s_p -1}   v \Vert_{L^2_t L^{\frac{2d}{d-3}}_x (I \times \mathbb{R}^d)} \leq A
\label{EQ2.19}
\end{equation}
Additionally, assume that we have
\begin{align}
\Vert (u_0 - v(t_0), u_1 - \partial_t v (t_0) ) \Vert_{\dot{H}^{s_p} \times \dot{H}^{s_p -1}} \leq A'
\label{EQ2.4}
\end{align}
as well as the smallness condition 
\begin{align}
\Vert D^{s_p -1} e \Vert_{L^{1}_t L^{2}_x} + \Vert S(t- t_0) (u_0 - v(t_0), u_1 - \partial_t v (t_0) ) \Vert_{S_p (I)} 
\leq \epsilon .
\label{EQ2.5}
\end{align}
Then there exists $\epsilon_0 = \epsilon_0 (d, p,  A, A') >0$  such that if  $0 < \epsilon < \epsilon_0$,  there is a unique solution $u$ of \eqref{nonlinear wave equations}
 in $I$ with $(u(t_0), \partial_t u (t_0))= (u_0, u_1)$ 
satisfying
 \begin{align}
\Vert u-v  \Vert_{\dot{S}(I)} & \leq C(d, p, A, A')  , \label{EQ2.20}
\\
 \Vert u - v \Vert_{S_p (I)} & \leq  C(d, p, A, A') \cdot \epsilon^{\theta} \label{EQ2.21}
 \end{align}
 where $\theta \in (0,1)$. 
 \end{Theorem}
Before we address the proof of  Theorem~\ref{perturbation_thm_Sobolev}, we will estimate the $\dot{S}(I)$ norm of the perturbed solution. 
\begin{Lemma} 
\label{Lem2.3}
Assume that $v$ solves the equation
\begin{align*}
\partial^2_{t} v
    - \Delta v 
     =    |v|^{p-1} v + e  
\end{align*}
in the sense of the corresponding  integral equation, and it satisfies 
\begin{align}
\sup_{t \in I} \Vert v \Vert_{\dot{H}^{s_p} \times \dot{H}^{s_p -1}} 
+\Vert D^{s_p -1}   v \Vert_{L^2_t L^{\frac{2d}{d-3}}_x (I \times \mathbb{R}^d)} & \leq A \label{EQ2.10}
 \\
\Vert D^{s_p -1} e \Vert_{L^{1}_t L^{2}_x}  & \leq A.  \label{EQ2.11}
\end{align}
Then we get
\begin{equation}
\Vert v \Vert_{\dot{S} (I)} \leq C(d,p, A). \label{EQ2.12}
\end{equation}
\end{Lemma}

\begin{proof}
First, we note that the hypothesis \eqref{EQ2.10} implies that the $S_p$ norm of $v$ stays bounded. By interpolation, 
\begin{equation}
\Vert v \Vert_{S_p (I)} \les \Vert v \Vert^{1- \frac{1}{p-1}}_{L^{\infty}_t \dot{H}^{s_p}_x} \Vert D^{s_p -1} v \Vert^{\frac{1}{p-1}}_{L^2_t L^{\frac{2d}{d-3}}_x }
\les A.
\label{EQ2.18}
\end{equation}
We  assume that $I = [0, T_+)$, and partition the interval $[0, T_+)$ as
\begin{equation*}
I = \cup_{j=1}^{K} I_j = \cup_{j=1}^{K} [t_{j-1}, t_j )
\end{equation*}
such that for each $j = 1, \ldots, K$, we get
\begin{align}
\Vert D^{s_p -1}   v \Vert_{L^2_t L^{\frac{2d}{d-3}}_x (I_j \times \mathbb{R}^d)} & \leq \eta \label{EQ2.13}
 \\
\Vert D^{s_p -1} e \Vert_{L^{1}_t L^{2}_x (I_j \times \mathbb{R}^d)}  & \leq \eta  \label{EQ2.14}
\\
\Vert S(t - t_{j-1}) (v (t_{j-1}), \partial_t v (t_{j-1})) \Vert_{S_p (I_j)} & <  \delta_0  \label{EQ2.15}
\end{align}
where $\delta_0$ is given by Theorem~\ref{Thm2.2}, and $\eta  >0$ is to be determined. 
%We also note that \eqref{EQ2.13} combined with the Strichartz inequality implies that 
%\begin{equation}
%\Vert S(t - t_{j-1}) (v (t_{j-1}), \partial_t v (t_{j-1})) \Vert \leq C \eta.  \label{EQ2.15}
%\end{equation}
Similar to \eqref{EQ2.18}, we get
\begin{align}
\Vert v \Vert_{S_p (I_j)}^{p-1}  \les A^{p-2} \eta 
\label{EQ2.18.2}
\end{align}
 on each $I_j$.
%Next, we set $\eta:= \min \left( \frac{1}{2C}, \frac{\delta_0}{C}  \right)$, where $\delta_0$ is given by Theorem~\ref{Thm2.2}, 
%where the constant $C$ is due to the implicit constants  in the  Strichartz inequality and Lemma~\ref{Nonlinear Estimate 1}, and depends only on $d, p$. 
We may then obtain 
\begin{equation}
\Vert v \Vert_{\dot{S} (I_j)} \leq C (d, p, A)
\label{EQ2.17}
\end{equation}
for each $j$ by following  the proof of Theorem~\ref{Thm2.2} and selecting $\eta$ to be sufficiently small. 
More precisely, using the Strichartz inequality, \eqref{EQ2.14}, and Lemma~\ref{Nonlinear Estimate 1} we obtain
\begin{align}
\begin{split}
\Vert v \Vert_{\dot{S} (I_j)}  &  \leq C \left( \sup_{t \in I_j} \Vert v \Vert_{\dot{H}^{s_p} \times \dot{H}^{s_p -1}} +  \Vert D^{s_p -1} e \Vert_{L^{1}_t L^{2}_x (I_j \times \mathbb{R}^d)} 
+ \Vert F(v) \Vert_{W' (I_j)}  \right)
\\ & 
\leq C \left( A + \eta + A^{p-2} \eta \Vert v \Vert_{\dot{S} (I_j)}  \right)
\end{split}
\label{EQ2.16}
\end{align}
On the second line above, we combined Lemma~\ref{Nonlinear Estimate 1} with \eqref{EQ2.18.2} and \eqref{embd2}. Also, the constant $C$ is due to the implicit constants in the  Strichartz inequality, Lemma~\ref{Nonlinear Estimate 1}, and \eqref{EQ2.18.2}. Letting $\eta= \eta (d, p,  A) >0$ to be sufficiently small, we obtain \eqref{EQ2.17}. 
Since the total number of intervals, $K$, must be finite due to the hypotheses \eqref{EQ2.10}--\eqref{EQ2.11}, which is given by
$K = K(A, d, p)$, we observe that  \eqref{EQ2.12} holds. 
\end{proof}

\begin{proof}[Proof of Theorem~\ref{perturbation_thm_Sobolev}]
Combining Lemma~\ref{Lem2.3} with \eqref{EQ2.19} and \eqref{EQ2.5}, we get
\begin{align*}
\Vert v \Vert_{\dot{S} (I)} \leq C (d,p, A). 
\end{align*}
Also, we note that by \eqref{embd4} and the Strichartz inequality,
\begin{align*}
\begin{split}
& \Vert S(t) (u_0 - v (t_0), u_1 - \partial v (t_0)) \Vert_{X (I)} \\
& \indeq \les  \Vert S(t) (u_0 - v (t_0), u_1 - \partial v (t_0)) \Vert_{S_p (I)}^{\theta} \Vert S(t) (u_0 - v (t_0), u_1 - \partial v (t_0)) \Vert_{\dot{S} (I)} \\
& \indeq \les \epsilon^{\theta} A'^{1- \theta} .
\end{split}
\end{align*}
Letting $\epsilon= \epsilon (d, p, A, A') >0$ to be sufficiently small, we apply Theorem~\ref{perturbation_thm}, which gives us
\begin{align*}
\Vert u- v \Vert_{\dot{S} (I)} & \leq C (d, p, A, A') \\
\Vert u-v \Vert_{X(I)} & \leq C (d, p, A, A') \epsilon^{\theta} .
\end{align*}
Finally, we obtain \eqref{EQ2.21} by applying \eqref{embd3}. 
\end{proof}

\begin{Remark}[Continuity property] 
\label{R2.4}
{\rm
As an application to the perturbation results stated above, we deduce that the flow associated to \eqref{nonlinear wave equations} has a continuity property. 
More precisely, let
 $(u_0, u_1) \in \dot{H}^{s_p} \times \dot{H}^{s_{p-1}}$ and let $\vec{u} (t)$ be the solution of \eqref{nonlinear wave equations}
 with maximal interval of existence
\begin{equation*}
\Imax (\vec{u}) = (T_- (\vec{u}), T_+ (\vec{u})).
\end{equation*} 
 Assume that  $(u_{0,n}, u_{1,n}) \to (u_0, u_1)$ in $\dot{H}^{s_p} \times \dot{H}^{s_{p-1}}$ and denote by $\vec{u_n} (t)$
 the corresponding solution with
 \begin{equation*}
 \Imax (\vec{u_n }) = (T_- (\vec{u_n }), T_+ (\vec{u_n })).
\end{equation*}  
Then, 
\begin{equation}
(T_- (\vec{u}), T_+ (\vec{u})) \subset (\limsup_n T_- (\vec{u_n}), \liminf_n T_+ (\vec{u_n})). 
\label{eq2.14}
\end{equation}
Moreover, for each $t \in (T_- (\vec{u}), T_+ (\vec{u}))$ we have
\begin{equation}
(u_n (t), \partial_t u_n (t)) \to (u(t), \partial_t u (t)) \quad \mbox{in~} \dot{H}^{s_p} \times \dot{H}^{s_{p-1}}. 
\label{eq2.15}
\end{equation}
}
\end{Remark}

%n%\newpage
\section{Concentration compactness procedure}
\label{sec:reduction}

The first component in establishing Theorem~\ref{main thm_1} is a concentration compactness argument. 
The approach we follow here was introduced by Kenig and Merle \cite{KM1, KM2} and studied further in several works \cite{KM11, DKM14, DL15, KLLS}. 

\subsection{Existence and compactness of a critical solution}
In order to highlight the essential tools in the proof of Theorem~\ref{main thm_1}, we begin with some notation.
\begin{definition}
For $A>0$, we define
\begin{align}
\mathcal{B} (A) := \{(u_0, u_1) \in \dot{H}^{s_p} \times \dot{H}^{s_p -1} :~ \sup_{t \in [0, T_+ (\vec{u})) } \Vert \vec{u}(t)\Vert_{\dot{H}^{s_p} \times \dot{H}^{s_p -1}} \leq A \}
\end{align}
where $\vec{u} (t)$ denotes the unique solution to \eqref{wave equations} in $\dot{H}^{s_p} \times \dot{H}^{s_p -1} $ with initial data $(u_0, u_1)$ and the maximal interval of existence $\Imax (\vec{u})= (T_- (\vec{u}), T_+ (\vec{u}))$. 
\end{definition}

\begin{definition}
We say that $\mathcal{SC} (A)$ holds if for each $(u_0, u_1) \in \mathcal{B} (A)$ we have $T_+ (\vec{u}) = \infty$ and $\Vert \vec{u} \Vert_{S_p ([0, \infty))} < \infty$. We also say that $\mathcal{SC} (A; (u_0, u_1))$ holds if $(u_0, u_1) \in \mathcal{B} (A)$ implies  $T_+ (\vec{u}) = \infty$ and $\Vert \vec{u} \Vert_{S_p ([0, \infty))} < \infty$.
\end{definition}

Analogously, we may define $\mathcal{B}_{rad} (A)$ and $\mathcal{SC}_{rad} (A)$ by restricting $(u_0, u_1)$ to radial initial data. 

\begin{Remark}
{\rm
By Theorem~\ref{Thm2.2} and Remarks \ref{R2.3}--\ref{R02} there exists $\delta>0$ sufficiently small so that 
%if $\Vert (u_0, u_1) \Vert_{\dot{H}^{s_p} \times \dot{H}^{s_p -1}} < \delta$, then 
$\mathcal{SC} ( \delta )$ holds.  Combined with Remark~\ref{R02} we deduce that Theorem~\ref{main thm_1} is equivalent to the statement that $\mathcal{SC}_{rad} (A)$ holds for all $A>0$. Hence, in the event that Theorem~\ref{main thm_1} fails, there exists a critical value $A_C >0$ so that for $A < A_C$, $\mathcal{SC} (A)$ holds, and for $A > A_C$, $\mathcal{SC} (A)$ fails. }
\end{Remark}

The next  result states that in the failure of Theorem~\ref{main thm_1}, there exists radial initial data $(u_{0,C}, u_{1,C}) \in \dot{H}^{s_p} \times \dot{H}^{s_p -1}$  so that
$\mathcal{SC} (A_C, (u_{0,C}, u_{1,C}))$ fails. 
%The solution $\vec{u}_C (t) $ to \eqref{wave equations} that corresponds to $(u_{0,C}, u_{1,C})$ satisfies a compactness property and is referred to as the critical solution. 
Furthermore, the solution to \eqref{wave equations} that  corresponds to $(u_{0,C}, u_{1,C})$  satisfies a compactness property and plays a crucial role in our discussion. 

\begin{Proposition} 
\label{Proposition3}
Suppose that Theorem~\ref{main thm_1} is false. Then, there exists $(u_{0,C}, u_{1,C})$ radial such that $\mathcal{SC} (A_C, (u_{0,C}, u_{1,C}))$ fails. 
Additionally, there exists a continuous function $\lambda: [0, T_+ (\vec{u}_C)) \to (0, \infty)$  so that the set 
\begin{equation}
\left\lbrace \left(\frac{1}{\lambda(t)^{\frac{2}{p-1}} } u_C \left( t, \frac{x}{\lambda (t)} \right), \frac{1}{\lambda(t)^{\frac{2}{p-1}+1 }} \partial_t u_C \left(t, \frac{x}{\lambda (t)}  \right) \right):~ t \in [0, T_+ (\vec{u}_C) )  \right\rbrace
\label{eq3.1}
\end{equation}
has compact closure in $\dot{H}^{s_p} \times \dot{H}^{s_p -1}$. 
\end{Proposition}

\begin{definition}
Let $\vec{u}_C (t)  \in \dot{H}^{s_p} \times \dot{H}^{s_p -1}$ be a radial solution to \eqref{wave equations}.  We say that $\vec{u}_C (t)$ is a critical solution if 
it satisfies the conclusions of Proposition~\ref{Proposition3}. More precisely, we have 
\begin{equation*}
\sup_{t \in [0, T_+ (\vec{u}_C))} \norm{\vec{u}_C (t)}_{\dot{H}^{s_p} \times \dot{H}^{s_{p-1}}} = A_C \comma \norm{\vec{u}_C }_{S_p ([0, T_+ (\vec{u}_C ) )} = \infty,
\end{equation*}
and there exists a continuous function $\lambda: [0, T_+ (\vec{u}_C)) \to (0, \infty)$ so that the set given in \eqref{eq3.1} is pre-compact in $\dot{H}^{s_p} \times \dot{H}^{s_{p-1}}$. 
\end{definition}
Once a critical solution $\vec{u}_C (t)$ is given, it is possible to construct another critical solution with  a corresponding scaling function that is bounded away from zero. 
We state this result next and refer the reader to \cite[Lemma~3.10]{KM10} for an analogous proof.
\begin{Lemma} 
\label{Lem20}
 There is a critical solution $\vec{\omega} (t)$ 
with a corresponding $\lambda_{\omega}$ continuously defined on $[0, T_+ (\vec{\omega}))$  
such that 
\begin{equation*}
\inf_{t \in [0, T_+ (\vec{\omega}))} \lambda_{\omega} (t) \geq A_0 >0. 
\end{equation*}
\end{Lemma} 

Going back to Proposition~\ref{Proposition3}, the main ingredient in extracting a critical solution is a profile decomposition theorem for linear solutions. The profile decomposition for the wave equation is introduced 
by Bahouri--Gerard \cite{BG99} for initial data belonging to $\dot{H}^1 \times L^2$ in three dimensions, and extended to higher dimensions by Bulut \cite{B10}.  Below, we state a higher dimensional version of the profile decomposition for initial data in $\dot{H}^{s_p} \times \dot{H}^{s_p-1}$. 

\begin{Theorem}[\text{\cite[Theorem~1.3]{B10}}]
Let $s \geq 1$ be given and let $(u_{0,n}, u_{1,n})_{n \in \mathbb{N}}$ be a bounded sequence in $\dot{H}^{s} \times \dot{H}^{s -1} (\mathbb{R}^d)$ with $d \geq 3$. 
Then there exists a subsequence of $(u_{0,n}, u_{1,n})$ (still denoted $(u_{0,n}, u_{1,n})$), a sequence $(V_0^j , V_1^j )_{j \in \mathbb{N}} \subset \dot{H}^{s} \times \dot{H}^{s -1} (\mathbb{R}^d)$, and a sequence of triples $(\epsilon_n^j, x_n^j, t_n^j) \in \mathbb{R}^+ \times \mathbb{R}^d \times \mathbb{R}$ such that, for every $j \neq j'$, 
\begin{align*}
\frac{\epsilon_n^j}{\epsilon_n^{j'}} + \frac{\epsilon_n^{j'}}{\epsilon_n^j} + \frac{ |t_n^j - t_n^{j'}|}{\epsilon_n^j} + \frac{ |x_n^j - x_n^{j'}|}{\epsilon_n^j} 
\to \infty \comma  n \to \infty
\end{align*}  
and for every $l \geq 1$, if 
\begin{align*}
V^j = S(t) (V_0^j , V_1^j ) \mbox{~and~} V_n^j (t,x) = \frac{1}{(\epsilon_n^j)^{\frac{d}{2} -s} } V^j \left( \frac{t- t_n^j}{\epsilon_n^j}, \frac{x- x_n^j}{\epsilon_n^j} \right)
\end{align*}
then
\begin{align*}
u_{0,n} (x) & = \sum_{j=1}^{l} V_n^j (0,x) +\omega_{0,n}^{l} (x)  \\
u_{1,n} (x) & = \sum_{j=1}^{l} \partial_t V_n^j (0,x) + \omega_{1,n}^{l} (x) 
\end{align*}
with 
\begin{align*}
\limsup_{n} \Vert S(t) (\omega_{0,n}^l , \omega_{1,n}^l) \Vert_{L^q_t L^r_x} \to 0 \comma  l \to \infty
\end{align*}
for every pair $(q,r)$  with $q,r \in (2, \infty)$ which satisfies
\begin{align*}
 \frac{1}{q} + \frac{d-1}{2r} \leq \frac{d-1}{4} \comma
 \frac{1}{q} + \frac{d}{r} = \frac{d}{2} -s. 
\end{align*}
For every $l \geq 1$, we also have
\begin{align*}
\Vert u_{0,n} \Vert^2_{\dot{H}^s} + \Vert u_{1,n} \Vert^2_{\dot{H}^{s-1}} 
 = \sum_{j=1}^l \left( \Vert V_0^j \Vert^2_{\dot{H}^s} + \Vert V_1^j \Vert^2_{\dot{H}^{s-1}} \right) 
+ \Vert \omega_{0,n}^{l} \Vert^2_{\dot{H}^s} +  \Vert \omega_{1,n}^{l} \Vert^2_{\dot{H}^{s-1}} + o (1)
\comma n \to \infty.
\end{align*}
\end{Theorem}

\begin{Remark}
{ \rm When the sequence $(u_{0,n}, u_{1,n})$ is radial, we may select $(V_0^j, V_1^j)$ radial with $x_n^j \equiv 0$. 
}
\end{Remark}

Invoking the profile decomposition theorem as stated above, the proof of Proposition~\ref{Proposition3} follows from the procedure 
developed in \cite[Section~3]{KM10}.  In a broad manner, the failure of Theorem~\ref{main thm_1} along with the profile decomposition result above 
leads to a minimizing sequence of non-scattering solutions to \eqref{wave equations} in $L^{\infty}_t (\dot{H}^{s_p} \times \dot{H}^{s_p -1})$ norm. 
Through further analysis, a critical solution $\vec{u}_C$ which possesses additional compactness properties may be extracted. The continuity  property of
$\lambda (t)$ on $[0, T_+ (\vec{u}_C))$ follows from the continuity of $\vec{u}_C (t)$ on $[0, T_+ (\vec{u}_C))$ in $\dot{H}^{s_p} \times \dot{H}^{s_p -1}$.  For a detailed treatment, please see Remark~5.4 of \cite{KM1}.

\subsection{The compactness property}
In view of the properties of a critical solution from Lemma~\ref{Lem20}, we deduce that Theorem~\ref{main thm_1} follows from 
the next result.
\cole
\begin{Theorem} 
\label{Theorem3.1}
Let $\vec{u} (t)$ be a radial solution of \eqref{wave equations}.
% with $p \geq 3$. 
Assume that there exists a continuous function $\lambda: [0, T_+ (\vec{u}) ) \to (0, \infty)$ 
so that
\begin{equation}
K_+ := \left\lbrace \left(\frac{1}{\lambda(t)^{\frac{2}{p-1}} } u \left( t, \frac{x}{\lambda (t)} \right), \frac{1}{\lambda(t)^{\frac{2}{p-1}+1 }} \partial_t u \left(t,  \frac{x}{\lambda (t)} \right) \right):~t \in [0, T_+ (\vec{u}))  \right\rbrace
\end{equation}
has compact closure in $\dot{H}^{s_p} \times \dot{H}^{s_p -1}$ and we have
\begin{equation}
\inf_{t \in [0, T_+ (\vec{u}) )} \lambda (t) >0 . 
\end{equation}
Then, $\vec{u} \equiv (0,0)$.
\end{Theorem}
\colb

Before we approach the proof of Theorem~\ref{Theorem3.1}, we consider two separate scenarios. Let $\vec{u} (t)$ be a solution as in Theorem~\ref{Theorem3.1}. 
First, we eliminate the case with a finite time blow-up, i.e., we cannot have  $T_+ (\vec{u})< \infty$. Secondly, we consider the case where $\Imax (\vec{u}) = \mathbb{R}$ and 
\begin{equation*}
\inf_{t \in \Imax (\vec{u})} \lambda (t) >0 . 
\end{equation*}
We argue that in this case $\vec{u} (t)$ must be the zero solution. The proof of Theorem~\ref{Theorem3.1} then follows from studying the properties of the scaling parameter $\lambda (t)$. 
%Since the hypothesis of  the theorem is based on the pre-compact trajectory of the critical solution, we define this property 
%next. 
To begin with, we introduce the following definition.  
 
\begin{definition}\label{compactness property}
Let $\vec{u} (t)$ be a solution of  \eqref{wave equations} defined on its maximal interval of existence $\Imax (\vec{u}) = (T_- (\vec{u}), T_+ (\vec{u}))$. 
We say that $\vec{u} (t)$ has the compactness property if there exists  $\lambda: \Imax (\vec{u}) \to (0, \infty)$ 
so that
% with \begin{align}
%\inf_{I} \lambda (t) \geq A_0 >0 
%\label{eq_comp_1}
%\end{align}
 the set 
\begin{align}
K = \left\lbrace \left( \frac{1}{\lambda(t)^{\frac{2}{p-1}}} u \left( t, \frac{x}{\lambda (t)}  \right), \frac{1}{\lambda(t)^{\frac{2}{p-1}+1}} \partial_t u\left( t, \frac{x}{\lambda (t)}  \right) \right) 
: t \in \Imax (\vec{u}) \right\rbrace
\label{eq_comp_2}
\end{align}
has compact closure in ${\dot{H}^{s_p} \times \dot{H}^{s_p -1}}$. 
\end{definition}

Note that in the definition above the scaling function $\lambda (t)$ is defined on the entire interval $\Imax$ as opposed to the half-open interval $[0, T_+ (\vec{u}))$.  
The fact that  pre-compactness is preserved when we pass from $K_+$ to $K$ is depicted in the next lemma. 

\begin{Lemma}
\label{Lem23}
Let $\vec{u} (t)$ be a solution of \eqref{wave equations} as in Theorem~\ref{Theorem3.1}. Let $\{ t_n \}_n$ be a sequence of times 
in $[0, T_+ (\vec{u}))$ such that $\lim_n t_n = T_+ (\vec{u})$. Assume that there exists $(v_0, v_1) \in \dot{H}^{s_p} \times \dot{H}^{s_p -1}$
such that 
\begin{equation}
\left( \frac{1}{\lambda (t_n)^{\frac{2}{p-1}}} u \left( t_n, \frac{x}{\lambda (t_n )} \right), \frac{1}{\lambda (t_n)^{\frac{2}{p-1} +1 }} \partial_t u \left( t_n,  \frac{x}{\lambda (t_n )} \right)    \right)
\to 
(v_0, v_1)
\quad 
\mbox{as~} n \to \infty
\label{eq3.2}
\end{equation}
in $\dot{H}^{s_p} \times \dot{H}^{s_p -1}$. Let $\vec{v} (t)$ be the solution of \eqref{wave equations} with initial data $(v_0, v_1)$ at $t=0$. 
Then, $\vec{v} \not\equiv (0,0)$  provided that $\vec{u}  \not\equiv (0,0)$. 
Additionally, $\vec{v} (t)$ has the compactness property. 
\end{Lemma}

By the hypothesis of Theorem~\ref{Theorem3.1}, the sequence on the left hand side of \eqref{eq3.2} belongs to the pre-compact set $K_+$. As a result, 
after passing to a subsequence, we deduce that there exists $(v_0, v_1) \in \overline{K}_+$ so that the limit in \eqref{eq3.2} holds. 
\begin{proof}
The proof of Lemma~\ref{Lem23} is very similar to those in \cite[Lemma~6.1]{KM11} and \cite[Claim~C.1]{DKM15.2}. We will show the changes below. 
First, we may assume that $(u_0, u_1) \neq (0,0)$ since the result becomes trivial with $(0,0)$ initial data. Below, we denote by $a = 2/ (p-1)$. 
%Under this hypothesis, we note that $(0,0) \notin \overline{K}_+$. 

Since $(u_0, u_1) \neq (0,0)$,  we have
\begin{equation}
\inf_{\Imax (\vec{u})} \norm{\vec{u} (t)}_{\dot{H}^{s_p} \times \dot{H}^{s_{p-1}}} > 0
\label{eq3.12} 
\end{equation}
by the local Cauchy theory. As a result,  for any sequence $\{ t_n \}_n \subset \Imax$ we deduce 
\begin{equation*}
\inf_n \norm{ \left(  \frac{1}{\lambda (t_n)^{a}} u \left( t_n, \frac{x}{\lambda (t_n)} \right), 
\frac{1}{\lambda (t_n)^{a+1}} \partial_t u \left( t_n,  \frac{x}{\lambda (t_n)}  \right) \right) }_{\dot{H}^{s_p} \times \dot{H}^{s_{p-1}}} >0. 
\end{equation*}
Hence, the limit  $(0,0) \notin \overline{K}_+$, which further implies that $(v_0, v_1) \neq (0,0)$. 
This proves the first claim.

\textit{Step~1.} We claim that for each $s \in (T_- (\vec{v}), T_+ (\vec{v}))$ we have
\begin{equation}
t_n + s/ \lambda (t_n) \geq 0 
\label{eq3.4}
\end{equation}
for large $n$. As $\{ t_n\}$ and $\{ \lambda (t_n) \}$  are non-negative sequences, the inequality above holds for $s \in [0, T_+ (\vec{v}))$.   
 We assume for a contradiction that \eqref{eq3.4} fails for some $s \in (T_- (\vec{v}), 0)$. We may then extract a subsequence so that 
\begin{equation}
t_n \lambda (t_n) + s <0 
\label{eq3.4c}
\end{equation}
for every $n$. Set $s_n = - t_n \lambda (t_n)$. Note that $s_n \in [s,0]$. By passing into a subsequence if necessary, we have
$\lim_n s_n = \theta \in [s,0] \subset (T_- (\vec{v}), T_+ (\vec{v}))$. Noting Remark~\ref{R2.4}  and the fact that $s_n + t_n \lambda (t_n) =0$, we get  
\begin{equation}
\left( \frac{1}{\lambda (t_n)^{a}} u \left(0, \frac{x}{\lambda (t_n)} \right), 
\frac{1}{\lambda (t_n)^{a+1}} \partial_t u \left(0, \frac{x}{\lambda (t_n)}  \right)  \right)
\to \left( v (\theta, x), \partial_t v ( \theta, x) \right)
\end{equation}
in $\dot{H}^{s_p} \times \dot{H}^{s_{p-1}}$. 
Since $(v_0, v_1) \neq (0,0)$, we also get 
$\vec{v} (\theta) \neq (0,0)$. 
Therefore, we obtain
\begin{equation}
\frac{1}{C}\leq \lambda (t_n) \leq C \quad \mbox{for every }n
\label{eq3.5}
\end{equation}
for some constant $C>0$, which will then yield a contradiction. If $T_+ (\vec{u}) = \infty$, \eqref{eq3.5} contradicts with \eqref{eq3.4c}. 
If $T_+ (\vec{u}) < \infty$, by Proposition~5.3 in \cite{KM11} ( see also  \cite[Prop.~5.3]{KM1}) we obtain
\begin{equation}
\lambda (t_n) \geq \frac{C_0}{T_+ (\vec{u}) - t_n} 
\end{equation} 
which implies that $\lambda (t_n) \to \infty$, contradicting \eqref{eq3.5}.

\textit{Step~2.} We aim to show that for every $s \in (T_-(\vec{v}), T_+ (\vec{v}))$ there exists  $\tilde{\lambda} (s)>0$ so that 
\begin{equation}
\left( \frac{1}{\tilde{\lambda} (s)^{a}} v \left(s, \frac{x}{\tilde{\lambda} (s)} \right), \frac{1}{\tilde{\lambda} (s)^{a+1}} \partial_t v \left(s, \frac{x}{\tilde{\lambda} (s)} \right) \right) \in \overline{K}_+
\period
\label{eq3.6}
\end{equation}
Setting $\tau_n =  t_n + s/ \lambda (t_n) $, we note that 
\begin{equation}
\left( \frac{1}{\lambda (\tau_n )^{a}} u \left(\tau_n, \frac{x}{\lambda (\tau_n) } \right), 
\frac{1}{\lambda (\tau_n )^{a+1}} \partial_t u \left( \tau_n, \frac{x}{\lambda (\tau_n)}  \right)  \right)
\in K_+
\end{equation}
as $t_n + s / \lambda (t_n) \geq 0$ for $n$ sufficiently large. By passing to a subsequence, we find $(\omega_0 (s), \omega_1 (s)) \in \overline{K}_+$
such that 
\begin{equation}
\left( \frac{1}{\lambda (\tau_n )^{a}} u \left(\tau_n, \frac{x}{\lambda (\tau_n)} \right), 
\frac{1}{\lambda (\tau_n )^{a+1}} \partial_t u \left(\tau_n,  \frac{x}{\lambda (\tau_n)} \right)  \right)
\to (\omega_0 (s, x), \omega_1 (s,x)) 
\label{eq3.14}
\end{equation}
in $\dot{H}^{s_p} \times \dot{H}^{s_{p-1}}$. 
At the same time, combining \eqref{eq3.2} and the continuity property of the flow as stated in Remark~\ref{R2.4} we also get
\begin{equation}
\left( \frac{1}{\lambda (t_n)^{a}} u \left( t_n + \frac{s}{\lambda (t_n)} , \frac{x}{\lambda (t_n )}\right), \frac{1}{\lambda (t_n)^{a +1 }} \partial_t 
u \left( t_n + \frac{s}{\lambda (t_n)} , \frac{x}{\lambda (t_n )} \right)    \right)
\to 
(v(s,x), \partial_s v(s,x))
\label{eq3.13}
\end{equation} 
in  $\dot{H}^{s_p} \times \dot{H}^{s_{p-1}}$. 
We  rescale \eqref{eq3.13} in $x$ by $\lambda (t_n) / \lambda (\tau_n)$ so that the convergence in \eqref{eq3.14} 
may be utilized. We then find that 
\begin{equation}
\left( \frac{\lambda (t_n)^{a}}{\lambda (\tau_n )^{a}} v \left(s, \frac{x \lambda (t_n)}{\lambda (\tau_n) } \right), 
\frac{\lambda (t_n)^{a+1}}{\lambda (\tau_n )^{a+1}} \partial_s v \left( s, \frac{x \lambda (t_n)}{\lambda (\tau_n)} \right)  \right)
\to (\omega_0 (s,x), \omega_1 (s,x))
\end{equation}
in  $\dot{H}^{s_p} \times \dot{H}^{s_{p-1}}$. 
Since $(\omega_0 (s), \omega_1 (s)) \neq (0,0)$ as it belongs to the compact set $\overline{K}_+$, we deduce that 
\begin{equation}
0 < \frac{1}{\tilde{C} (s)} \leq \frac{\lambda (t_n)}{\lambda (t_n + s/ \lambda (t_n))} \leq \tilde{C} (s) < \infty
\end{equation}
for every $n$. 
Therefore, we may find a further subsequence  so that 
\begin{equation}
\lim_{n \to \infty} \frac{\lambda (t_n)}{\lambda (t_n + s / \lambda (t_n))} =: \tilde{\lambda} (s) \in (0, \infty)
\end{equation}
and 
\begin{equation*}
\left( \frac{1}{\tilde{\lambda} (s)^{a}} v \left(s, \frac{x}{\tilde{\lambda} (s)} \right), \frac{1}{\tilde{\lambda} (s)^{a+1}} \partial_s v \left( s, \frac{x}{\tilde{\lambda} (s)} \right) \right) \in \overline{K}_+
\end{equation*}
for every $s \in (T_- (\vec{v}) , T_+ (\vec{v}))$, which completes the proof. 
\end{proof}

%Following the earlier works by Duyckaerts, Kenig, and Merle \cite{DKM14}, and Dodson and Lawrie \cite{DL15}  in dimensions three and five,  respectively,  the proof of Theorem~\ref{main thm_1} reduces to the following rigidity result. 

%\cole
%\begin{Proposition}
%\label{Proposition1}
%Let $\vec{u} (t)$ be a solution of the \eqref{wave equations} which has the compactness property on its maximal interval of existence 
%$\Imax= ( T_- (\vec{u} (0)), T_+ (\vec{u} (0)) )$. 
%Then, $\Imax = \mathbb{R}$ and $\vec{u} (t) \equiv (0,0)$. 
%\end{Proposition}
%n%\newpage
%\colb

Next, we show that there is no solution $\vec{u}(t)$ to \eqref{wave equations} as in Theorem~\ref{Theorem3.1} with $T_+ (\vec{u}) < \infty$. 
In \cite[Section~3]{DKM15}, the authors consider the equation \eqref{wave equations} under the hypothesis that $p$ is an odd integer or large enough so that the local well-posedness theory holds. 
Proposition~3.1 in \cite{DKM15} shows that a solution of the equation  \eqref{wave equations} which has the compactness property on its maximal interval of existence is global. For exponents $p$ that are not odd integers, the range  $p> N/2$  is provided as a sufficient condition in which $N$ denotes the dimension. The local well-posedness theory in Section~\ref{sec:A review of the Cauchy problem} lets us carry through the proof of Proposition~3.1 in \cite{DKM15} and  eliminate the possibility of a self similar  solution that  blows up in finite time.  For convenience, we will provide the details below. 
%in Section~\ref{sec: rigidity}. 

\cole
\begin{Proposition}[\text{\cite[Proposition~3.1]{DKM15}}]
Let 
%$p \geq 3$ and  let
$\vec{u}(t)$ be a solution of \eqref{wave equations} with the compactness property. Then, $\vec{u} (t)$ is global. 
\end{Proposition}
\colb

\begin{proof}[Sketch of the proof] Let $\vec{u}(t)$ be a solution of \eqref{wave equations} on its maximal interval of existence $\Imax (\vec{u}) = (T_- (\vec{u}), T_+ (\vec{u}) )$
which has the compactness property as defined in Definition~\ref{compactness property}. Since we are concerned with the radial case, we will assume that $\vec{u} (t)$
is a radial solution.  Fixing a non-zero radial solution $\vec{u} (t)$, we simplify the notation and write $\Imax= (T_-, T_+)$. 
For a contradiction, suppose that $T_-$ is finite. 

\textit{Step~1.} We claim that for every $t \in \Imax$
\begin{equation}
\mbox{supp}~\vec{u} (t) \subset \left\lbrace |x| \leq |T_- - t|\right\rbrace.
\label{eq3.7}
\end{equation}
The proof of this step may be completed by following the same strategy in \cite[Section~4]{KM11}, which relies on finite speed of propagation, 
the small data theory, and the perturbation result we
established in Section~\ref{sec:A review of the Cauchy problem}. In particular, since Theorem~\ref{Thm2.2} and Theorem~\ref{perturbation_thm} hold for $p$ satisfying \eqref{eq0.1}, we may proceed with the rest of the proof as outlined below. 
Firstly, we obtain that for all $t \in \Imax$
 \begin{equation}
 \lambda (t) \geq \frac{C(\bar{K})}{t - T_-} >0
 \end{equation}
 as done in \cite[Lemma~4.14]{KM11} (cf. \cite[Lemma~4.7]{KM2} for the details). In particular, we deduce 
 \begin{equation}
 \lim_{t \to T_-} \lambda (t) = \infty
 .
 \end{equation}
 Then, by following the same arguments in \cite[Lemma~4.15]{KM11} we prove the statement \eqref{eq3.7}. 
 
 \textit{Step~2.} We conclude the proof by finding a monotone function in time in terms of the solution $\vec{u} (t)$. We simply provide the proof of Step.~2 in \cite[Proposition~3.1]{DKM15} for the sake of completeness. Let
 \begin{equation}
 y (t) = \int u^2 (t,x)~dx.
 \end{equation}
 By Step~1,  $y (t)$ is well-defined for all $t \in \Imax$, and furthermore  $\vec{u} (t) \in \dot{H}^1 \times L^2$. 
 Using the equation \eqref{wave equations}, we obtain
 \begin{equation}
  y' (t) = 2 \int u(t,x) \partial_t u (t,x)~dx
 \label{eq3.9}
 \end{equation}
 and 
 \begin{align*}
 y'' (t) = 2 \int (\partial_t u (t,x))^2 dx - 2 \int |\nabla u (t,x)|^2 dx + 2 \int |u (t,x)|^{p+1} dx. 
 \end{align*}
 
 Next, we recall that the conserved energy for the flow is given by
 \begin{equation*}
 E (\vec{u} (t)) = \int \left(  \frac{1}{2}|\partial_t u(t,x)|^2 + \frac{1}{2} |\nabla u (t,x)|^2 - \frac{1}{p+1} |u (t,x)|^{p+1} \right) dx. 
 \end{equation*}
 Noting that $\vec{u} (t)$ is uniformly bounded in $\dot{H}^{s_p} \times \dot{H}^{s_{p-1}}$ by the compactness property with $s_p >1$, 
 we deduce that the condition \eqref{eq3.7} on the support of $\vec{u}(t)$ leads to 
 \begin{equation}
 \lim_{t \searrow T_-} E (\vec{u} (t)) =0 \quad \mbox{and} \quad \lim_{t \searrow T_-} y(t) = \lim_{t \searrow T_-} y' (t)=0.
 \label{eq3.8}
 \end{equation}
 We then obtain by conservation of the energy that
 \begin{equation}
 E(\vec{u} (t)) =0
 \end{equation}
 for all $t \in \Imax$, and we rewrite 
 \begin{equation}
 y'' (t) = (p+3) \int (\partial_t u (t, x))^2 dx + (p-1) \int |\nabla u (t,x)|^2 dx >0. 
 \label{eq3.10}
 \end{equation}
 Thus, we must have $y' (t)>0$ for all $t \in \Imax$.  We note that in the case $T_+ < \infty$ we also obtain
 \begin{equation}
 \lim_{t \nearrow T_+} y' (t) =0
 \end{equation}
 which contradicts with the fact that $y (t)$ is a strictly convex function with the limit \eqref{eq3.8}. We then deduce that $T_+ = \infty$. 
 Combining \eqref{eq3.9} with \eqref{eq3.10}  we obtain
 \begin{align}
 y' (t)^2 \leq 4 \left( \int u^{2} (t,x)~dx \right) \left( \int (\partial_t u(t,x) )^2 dx\right)
  \leq \frac{4}{p+3} y (t) ~ y''(t).
  \label{eq3.11}
 \end{align}
 Using \eqref{eq3.11} and the fact that $y' (t)>0$ for all $t \in \Imax$, we claim that $y^{- (p-1)/4}$ is strictly decreasing and concave down. To see this, 
 note that
 \begin{equation*}
 \frac{d}{dt} \left(y^{- (p-1)/4} (t) \right) = - \frac{(p-1)}{4} y^{- (p+3)/4} (t) ~y' (t) <0
 \end{equation*}
 and 
 \begin{equation*}
 \frac{d^2}{dt^2} \left( y^{- (p-1)/4} (t) \right) =  \frac{(p-1)}{4} y^{- (p+7)/4} (t) \left( \frac{p+3}{4} (y' (t))^2 - y(t)~ y'' (t) \right) \leq 0. 
 \end{equation*}
 This however contradicts the fact that $T_+ = \infty$. 
\end{proof}

We apply the above proposition for the solution $\vec{v}(t)$ constructed in Lemma~\ref{Lem23}, which satisfies the compactness property on its maximal interval of existence
$(T_- (\vec{v}), T_+ (\vec{v}))$. Note that Remark~\ref{R2.4} implies that if $\vec{v} (t)$ is a global, then $\vec{u} (t)$ as in Theorem~\ref{Theorem3.1} must be a global solution as well. 
Having eliminated the case $T_+ (\vec{u}) < \infty$, we focus on the following result. 
\cole
\begin{Proposition}
\label{Proposition1}
Let $\vec{u} (t)$ be a  radial solution of  \eqref{wave equations} with $\Imax (\vec{u}) = \mathbb{R}$, which has the compactness property on $\mathbb{R}$. Suppose that
we have
\begin{equation}
\inf_{t \in (- \infty, \infty) } \lambda (t) > 0. 
\label{eq3.3}
\end{equation}
Then,  $\vec{u}  \equiv (0,0)$. 
\end{Proposition}
%n%\newpage
\colb

The proof of Theorem~\ref{Theorem3.1} may be completed by proceeding as in \cite[Section~6]{DKM14}. 
More specifically, letting $\vec{u}$ be as in Theorem~\ref{Theorem3.1}, we follow the arguments in Lemma~6.3--6.6 and employ Proposition~\ref{Proposition1}
to  examine the further properties of the corresponding scaling function $\lambda (t)$, and we arrive at the conclusion that $\vec{u} \equiv (0,0)$. 

The remainder of the article deals with the proof of Proposition~\ref{Proposition1}. Firstly, we focus on showing that solutions of \eqref{wave equations} with the compactness property enjoy additional spacial decay, which yields the fact that the trajectory of $\vec{u} (t) \in \dot{H}^1 \times L^2$. Next, we highlight a family of  singular stationary solutions with asymptotic properties similar to those of solutions given in the hypothesis of the Proposition~\ref{Proposition1}
yet these singular stationary solutions fail to belong to the critical space $\dot{H}^{s_p}$. Finally, using the exterior energy estimates from \cite{KLLS15} we may go through the rigidity method in three main steps to show that a non-zero radial solution of \eqref{wave equations} with the compactness property has to coincide with a singular stationary solution. 
As a result, we obtain the desired conclusion of Proposition~\ref{Proposition1} that $\vec{u}  \equiv (0,0)$. 

\section{Decay results for solutions with the compactness property}
\label{sec:Decay}
In this section, 
we apply 
the double Duhamel method to study the decay rates
of solutions to the Cauchy problem \eqref{wave equations} which has the compactness property. 
  The methods in this section are analogous to the discussion 
in \cite[Section~4]{DL15} for the focusing cubic wave equation in $\mathbb{R}^{5}$.

First, we recall some preliminary facts from harmonic analysis which will be frequently used throughout the section.
We begin with a radial Sobolev inequality, quoted verbatim from \cite[Corollary~A.3]{TVZ07} for convenience of readers. 
\begin{Lemma}[Radial Sobolev inequality]
\label{Lem21}
Let $1 \leq p,q \leq \infty$, $0< s < d$, and $\beta \in \mathbb{R}$ obey the conditions 
\begin{equation*}
\beta > - \frac{d}{q'} \comma 1 \leq \frac{1}{p} + \frac{1}{q} \leq 1+s 
\end{equation*}
and the scaling condition
\begin{equation*}
 d - \beta -s =  \frac{d}{p'} + \frac{d}{q'}
\end{equation*}
with at most one of the equalities 
\begin{equation*}
p = 1, ~ p = \infty,~ q =1,~ q = \infty,~ \frac{1}{p} + \frac{1}{q} = 1+s
\end{equation*}
holding. Then, for any radial function $f \in \dot{W}^{s, p} (\mathbb{R}^d)$, we have
\begin{equation*}
\Vert |x|^{\beta} f   \Vert_{L^{q'} (\mathbb{R}^d)} \leq C \Vert  D^s f  \Vert_{L^{p} (\mathbb{R}^d)}
\period
\end{equation*}
\end{Lemma}

\begin{Remark} 
\label{R4.1}
{\rm
For the remainder of this section, we may assume that 
%$p \geq 3$ and 
$\vec{u} (t)$ is a solution of \eqref{wave equations}
as in Proposition~\ref{Proposition1}. Namely, $\vec{u} (t)$ is a radial solution of \eqref{wave equations} with 
$\Imax (\vec{u}) = \mathbb{R}$. Additionally, $\vec{u} (t)$ has  the compactness property on $\mathbb{R}$
and the corresponding scaling parameter $\lambda$  satisfies
\begin{equation}
\inf_{t \in (- \infty, \infty)} \lambda (t) >0.
\label{eq4.58}
\end{equation}
%Nevertheless, the results in this section continue to hold when we waive our assumption that $\vec{u} (t)$ is global. 
%Also, we may allow $p >2$. 
}
\end{Remark}
Next, we state a quantitative result for solutions with the compactness property. 
By the Arzela-Ascoli theorem, we may simply obtain the following uniform estimates on the  $\dot{H}^{s_p} \times \dot{H}^{s_{p-1}}$ norm 
of a solution that has the compactness property. For similar estimates, see \cite[Remark~4]{DL15}.
\begin{Lemma}[Uniformly Small Tails] 
\label{Lemma_tails}
Let $\vec{u} (t)$ be a solution of the equation \eqref{wave equations} as in Remark~\ref{R4.1}. 
Then for any $\eta >0$ there are  $0< c(\eta) < C(\eta) < \infty$ such that for all $t \in \mathbb{R}$
we have
\begin{align}
\begin{split}
& \int_{|x| \geq \frac{C(\eta)}{\lambda (t)}} \left|  D^{s_p} u (t,x)  \right|^2 dx
+
\int_{|\xi| \geq C(\eta) \lambda (t)} |\xi|^{2 s_p} |\hat{u} (t, \xi)|^2 d\xi 
 \leq \eta^2 \\
& \int_{|x| \leq \frac{c(\eta)}{\lambda (t)}} \left|  D^{s_p} u (t,x)  \right|^2 dx
+
\int_{|\xi| \leq c(\eta) \lambda (t)} |\xi|^{2 s_p} |\hat{u} ( t, \xi )|^2 d\xi 
 \leq \eta^2 \\
& \int_{|x| \geq \frac{C(\eta)}{\lambda (t)}} \left|  D^{s_p -1} u_t (t,x)  \right|^2 dx
+
\int_{|\xi| \geq C(\eta) \lambda (t)} |\xi|^{2 s_p -2} |\hat{u}_t (t, \xi)|^2 d\xi 
 \leq \eta^2 \\
& \int_{|x| \leq \frac{c(\eta)}{\lambda (t)}} \left|  D^{s_p -1} u_t (t, x)  \right|^2 dx
+
\int_{|\xi| \leq c(\eta) \lambda (t)} |\xi|^{2 s_p -2} |\hat{u}_t ( t, \xi )|^2 d\xi 
 \leq \eta^2 .
\end{split}
\label{eq4.uniform_tails}
\end{align}
\end{Lemma}

We will also utilize the following version of Duhamel's formula for solutions to \eqref{wave equations} with the compactness property. 
The standard Duhamel formula combined with 
the fact that the linear part of the evolution vanishes weakly in $\dot{H}^{s_p} \times \dot{H}^{s_p -1}$ yields the following lemma. 
Analogous results on weak limits are proved in \cite[Section~6]{TVZ08} and \cite[Proposition~3.8]{S13}.  
\begin{Lemma}[Weak Limits]
\label{Lemma_weak_limits}
Let $\vec{u} (t)$ be a solution of the equation \eqref{wave equations} 
as in Remark~\ref{R4.1}. Then, for any $t_0 \in \mathbb{R}$ we have
\begin{align*}
u (t_0) &  =  - \lim_{T \to \infty } \int_{t_0}^{T} \frac{\sin ((t_0 - \tau)  \sqrt{- \Delta})}{\sqrt{- \Delta}} ~  |u|^{p-1} u ~d \tau
\qquad   \mbox{ weakly in }~\dot{H}^{s_p} (\mathbb{R}^d)
\\
u_{t} (t_0) & = - \lim_{T \to \infty} \int_{t_0}^{T} {\cos  ((t_0 - \tau)  \sqrt{- \Delta})} ~  |u|^{p-1} u ~d \tau
\qquad \mbox{ weakly in }~\dot{H}^{s_p -1} (\mathbb{R}^d)
\\ 
u (t_0) &  =   \lim_{T \to -\infty } \int_{T}^{t_0} \frac{\sin ((t_0 - \tau)  \sqrt{- \Delta})}{\sqrt{- \Delta}} ~  |u|^{p-1} u ~d \tau
\qquad   \mbox{ weakly in }~\dot{H}^{s_p} (\mathbb{R}^d)
\\
u_{t} (t_0) & =  \lim_{T \to -\infty } \int_{T}^{ t_0} {\cos  ((t_0 - \tau)  \sqrt{- \Delta})} ~ |u|^{p-1} u ~d \tau
\qquad \mbox{ weakly in }~\dot{H}^{s_p -1} (\mathbb{R}^d).
\end{align*}
%Then $e^{-it \Delta } u(t)$ weakly converges to zero in 
 %$\dot{H}^{s_p} $ as $t \to \sup(\Imax)$ or $t \to \inf(\Imax)$. 
\end{Lemma}

%%%%% introduce the proposition

The following is the main result of this section on the decay of compact solutions to the equation \eqref{wave equations}. 
\begin{Proposition}
\label{decay_proposition}
Let $\vec{u} (t)$ be a solution to \eqref{wave equations} as in Remark~\ref{R4.1}. 
Then, for all $t \in \mathbb{R}$ 
\begin{align}
\Vert\vec{u} (t) \Vert_{\dot{H}^{3/4} \times \dot{H}^{-1/4} (\mathbb{R}^d)} \leq C_p
\label{eq_decay1}
\end{align}
where the constant $C$ is uniform in time. 
\end{Proposition}
The proof of Proposition~\ref{decay_proposition} follows from a double Duhamel technique as shown in \cite{DL15} and 
\cite{TVZ08}. 
Following the procedure introduced in \cite[Section~4.2]{DL15}, we define 
\begin{align}
v = u + \frac{i}{\sqrt{- \Delta}} u_t 
\period
\label{eq4.1}
\end{align}
As $\vec{u}$ solves \eqref{wave equations}, we get 
\begin{align}
v_t & = u_t + \frac{i}{\sqrt{- \Delta}} \left( \Delta u +  |u|^{p-1} u \right) \\
& = -i \sqrt{- \Delta} v + \frac{i}{\sqrt{- \Delta}}  |u|^{p-1} u \period
\label{eq4.2}
\end{align} 
Duhamel's formula then gives us
\begin{align}
v(t) = e^{-it \sqrt{- \Delta}} v(0) + i  \int_0^{t} \frac{e^{-i (t - \tau) \sqrt{- \Delta}}}{\sqrt{-\Delta}}
 |u|^{p-1} u (\tau)~d\tau
\period
\label{eq4.3}
\end{align} 
Assuming that $\vec{u} (t)$ has the compactness property, 
%which then has to be a global solution as discussed in Section~\ref{sec:reduction}, 
we deduce by Lemma~\ref{Lemma_weak_limits} that for any $t_0 \in \mathbb{R}$
\begin{align}
\int_{t_0}^{T} \frac{e^{-i (t_0 - \tau) \sqrt{- \Delta}}}{\sqrt{-\Delta}}
u |u|^{p-1} (\tau)~d\tau
\rightharpoonup 
i v (t_0)
\comma \mbox{as } ~T \to \pm \infty
\label{eq4.3p}
\end{align}
weakly in $\dot{H}^{s_p}$.
 Moreover, 
 \begin{align}
 \Vert \vec{u} (t) \Vert_{\dot{H}^{s_p} \times \dot{H}^{s_p -1}} \cong
 \Vert v(t) \Vert_{\dot{H}^{s_p}}
 \period
 \label{eq4.4}
 \end{align}
 
 We may now begin the proof of Proposition~\ref{decay_proposition}.
\begin{proof}[Proof of Proposition~\ref{decay_proposition}]
Our goal  is to find a sequence $\beta = \{ \beta_k\}$ of positive numbers such that
\begin{align}
\sup_{t \in \mathbb{R}} \Vert \left( P_k u (t), P_k u_t (t) \right) \Vert_{\dot{H}^{3/4} \times \dot{H}^{-1/4}} \les 2^{-\frac{3k}{4}} \beta_k
 \label{eq4.5}
\end{align}
 for all $k \in \mathbb{Z}$, and
 \begin{align}
 \Vert \{ 2^{-\frac{3k}{4}} \beta_k \}_k\Vert_{\ell^{2}} \les 1 \period
 \label{eq4.6}
 \end{align}
A sequence $\beta $ that satisfies the above properties is called a frequency envelope. In this section, 
 $P_k$ denotes the Littlewood-Paley projection corresponding to the dyadic number $2^k$. The definition of $P_k$ is recalled in \eqref{eq1.7p}. Equivalently, $P_k f$ is given 
 by convolution 
 \begin{align}
 P_k f = 2^{d k} \check{\phi} (2^k \cdot) \ast f
 \label{convolution_pk}
 \end{align}
 where $\check{\phi}$ belongs to the Schwartz class. 
 
  \begin{Claim}
 \label{claim_1}
 A frequency envelope that satisfies \eqref{eq4.5}--\eqref{eq4.6} may be defined as below: we take
 \begin{align}
 \beta_k := 1 \inon{ for }~ k \geq 0
 \label{eq4.7}
 \end{align}
 and for $k <0$, we set
 \begin{align}
 \beta_k := \sum_{j} 2^{- |j-k|} a_j
 \label{eq4.8}
 \end{align}
 where 
 \begin{align}
 a_j = 2^{s_p j} \Vert P_j u \Vert_{L_t^{\infty} L^2} + 2^{( s_{p} -1) j} \Vert P_j u_t \Vert_{L_t^{\infty} L^2} 
 \inon{ for }~ j \in \mathbb{Z}.
 \label{eq4.9}
 \end{align}
 \end{Claim}
 
 % \begin{Remark}
 %\label{remark1}
 Recalling the definition of $v$ in \eqref{eq4.3} we observe that
 \begin{align}
 \Vert P_j  v \Vert_{L_t^{\infty} \dot{H}^{s_p}} \cong a_j
 \period
 \label{eq4.10}
 \end{align}
 In order to verify \eqref{eq4.5}--\eqref{eq4.6} we first estimate $\Vert P_j  v \Vert_{L_t^{\infty} \dot{H}^{s_p}}$ for $j \in \mathbb{Z}^{-}$.  
 Throughout the proof of Claim~\ref{claim_1}, the estimates will be uniform in $t$. For that reason, 
it suffices to estimate the term  $\Vert P_j  v (0) \Vert_{\dot{H}^{s_p}}$, i.e. we set $t=0$.  
We also note that the implicit constants carried through the computations below are
allowed to depend on the norm $ \Vert v \Vert_{L^{\infty}_t \dot{H}^{s_p}}$.
%\end{Remark}

 Let $M= 2^j$, for $j \in \mathbb{Z}^{-}$ be an arbitrary dyadic number in $(0,1)$. By \eqref{eq4.3}, 
 \begin{align*}
 \left\langle P_M v(0), P_M v(0) \right\rangle_{\dot{H}^{s_p}} = 
 \left\langle P_M \left(e^{i T_1 \sqrt{-\Delta}} v(T_1) -i  \int_0^{T_1} \frac{e^{i \tau \sqrt{-\Delta}}}{\sqrt{-\Delta}}  |u|^{p-1} u (\tau)~d\tau \right) , P_M v(0)  \right\rangle_{\dot{H}^{s_p}}
\end{align*}  
 for any $T_1>0$.  We then take the limit $T_1 \to \infty$, which yields
 \begin{align}
 \begin{split}
&  \left\langle P_M v(0), P_M v(0) \right\rangle_{\dot{H}^{s_p}}  \\
 &  \qquad  = \lim_{T_1 \to \infty} 
  \left\langle P_M \left(e^{i T_1 \sqrt{-\Delta}} v(T_1) -i  \int_0^{T_1} \frac{e^{i \tau \sqrt{-\Delta}}}{\sqrt{-\Delta}} |u|^{p-1} u (\tau)~d\tau \right) , P_M v(0)  \right\rangle_{\dot{H}^{s_p}} \\
  & \qquad = 
  -i
   \left\langle P_M \int_0^{\infty} \frac{e^{i \tau \sqrt{-\Delta}}}{\sqrt{-\Delta}}  |u|^{p-1} u (\tau)~d\tau , P_M v(0)  \right\rangle_{\dot{H}^{s_p}} 
   \period
   \end{split}
   \label{eq4.11}
 \end{align}
  On the last line in \eqref{eq4.11}, we used Lemma~\ref{Lemma_weak_limits}  to have the weak limit
 \begin{align*}
 \lim_{t \to \infty} e^{it \sqrt{- \Delta}} v(t) = 0
 \end{align*}
 in ${\dot{H}^{s_p}}$. We also observe that
  \begin{align*}
 \lim_{t \to -\infty} e^{it \sqrt{- \Delta}} v(t) = 0
 \period
 \end{align*}
weakly in ${\dot{H}^{s_p}}$. Similarly, we use the formula \eqref{eq4.3} on the second term in \eqref{eq4.11}, and take the weak limit 
$ T \to -\infty$ to obtain the reduction 
 \begin{align}
 \begin{split}
&  \left\langle P_M v(0), P_M v(0) \right\rangle_{\dot{H}^{s_p}}  \\
 & \qquad = 
   \left\langle P_M \int_0^{\infty} \frac{e^{i \tau \sqrt{-\Delta}}}{\sqrt{-\Delta}}  |u|^{p-1} u (\tau)~d\tau , 
   P_M \int^0_{-\infty} \frac{e^{i \tau \sqrt{-\Delta}}}{\sqrt{-\Delta}}  |u|^{p-1} u (\tau)~d\tau \right\rangle_{\dot{H}^{s_p}} 
   \period
   \end{split}
   \label{eq4.12}
\end{align}
Next, we take a radial non-increasing bump function $\chi \in C_0^{\infty} (\mathbb{R}^d)$, which satisfies
 \begin{align*}
\chi (x) = 
\begin{cases}
1 & \inon{if  }~  |x| \leq 1, \\
0 & \inon{if  }~  |x| \geq 2 \period
\end{cases}
\end{align*}
 Also, let $c>0$ be a small fixed constant, say $c= 1/4$. We then express the $\dot{H}^{s_p}$ inner product 
  \eqref{eq4.12} as
\begin{align}
\begin{split}
\langle A+B, \tilde{A}+ \tilde{B} \rangle_{\dot{H}^{s_p}} = 
\langle A, \tilde{A}+ \tilde{B} \rangle_{\dot{H}^{s_p}} + 
\langle A+B,  \tilde{A} \rangle_{\dot{H}^{s_p}} 
- \langle A, \tilde{A} \rangle_{\dot{H}^{s_p}} 
+ \langle B, \tilde{B} \rangle_{\dot{H}^{s_p}} 
\end{split}
\label{EQ4.2}
\end{align} 
 where 
\begin{align}
 \begin{split}
 A &:= 
  \int_0^{\Lambda/M} \frac{e^{i \tau \sqrt{-\Delta}}}{\sqrt{-\Delta}} P_M  |u|^{p-1} u (\tau)~d\tau
  +  \int_{\Lambda/M}^{\infty} \frac{e^{i \tau \sqrt{-\Delta}}}{\sqrt{-\Delta}}  (1- \chi) \left(x/{c \tau} \right) P_M  |u|^{p-1} u (\tau)~d\tau \\
  & = A_1 + A_2
 \end{split}
 \label{eq4.13}
 \end{align}
 and 
\begin{equation}
B := \int_{\Lambda/M}^{\infty} \frac{e^{i \tau \sqrt{-\Delta}}}{\sqrt{-\Delta}}  \chi \left(x/{c \tau} \right) P_M  |u|^{p-1} u (\tau)~d\tau \period \label{eq4.14} 
\end{equation} 
The terms $\tilde{A}$ and $\tilde{B}$ are defined analogously in the negative time direction. 
 The constant $\Lambda>0$ will be determined below.  
 
 We remark that the left side of \eqref{EQ4.2} is a non-negative real number. As a result, it suffices to estimate the real part of each term on the right side of \eqref{EQ4.2}; this will be particularly important when we estimate $\langle B, \tilde{B} \rangle_{\dot{H^{s_p}}}$.
 
  First, we treat the term 
 $\langle A, \tilde{A} \rangle$ by estimating $\Vert A \Vert_{\dot{H}^{s_p}}$ and $\Vert \tilde{A} \Vert_{\dot{H}^{s_p}}$.
 
 \begin{Claim}
\label{claim_2}
Let $\eta >0$ be arbitrary. There is $N_0 >0 $ so that
\begin{equation}
\Vert A_1 \Vert_{\dot{H}^{s_p}} \les \Lambda M^{s_p} \eta^{p-1}  \Vert P_{> M/2p} u \Vert_{L_{t}^{\infty} L^{2}} + \Lambda M^{s_p} N_0^{-s_p} \period
\label{eq4.28}
\end{equation}
\end{Claim} 

%check the rest of the proof. 

First, we note that
\begin{align}
\begin{split}
& \norm{ P_M \left(   \int_0^{\Lambda/M} \frac{e^{i \tau \sqrt{-\Delta}}}{\sqrt{-\Delta}}   |u|^{p-1} u (\tau)~d\tau \right)  }_{\dot{H}^{s_p}} \\
& \qquad \les 
  M^{s_p -1} \norm{ \int_0^{\Lambda/M} {e^{i \tau \sqrt{-\Delta}}}  P_M (  |u|^{p-1} u ) (\tau)~d\tau }_{L^2} \\
& \qquad \les \Lambda M^{s_p -2} \norm{ P_M  (u |u|^{p-1}) }_{L_t^{\infty} L^2} \period  
\end{split}
\label{eq4.15}
\end{align}
Let $\eta >0$ be a small positive number. Since $u$ has the compactness property, we have $\inf_{t \in \mathbb{R}} \lambda (t) > A_0 $
(cf. \eqref{eq4.58}) for some positive constant $A_0$, and so Lemma~\ref{Lemma_tails} yields that there is a positive number $N_0 = N_0 (\eta)$ such that
\begin{equation}
\Vert P_{\leq N_0} u \Vert_{\dot{H}^{s_p}}^2  \leq \sum_{N \leq N_0} N^{2 s_p} \Vert P_N u \Vert_{L^2}^2 \les   \eta^2
\label{EQ4.15p}
\end{equation}
which then leads to
\begin{equation}
\Vert P_{\leq N_0} u \Vert_{L^{{d (p-1)}/{2}}} \les \eta
\label{eq4.16}
\end{equation}
 by the Sobolev embedding. Note that we also deduce from \eqref{EQ4.15p} that
 \begin{equation}
 \Vert P_{\leq K} u \Vert_{L^{d (p-1)/2 }} \les \eta 
 \label{EQ4.16p}
 \end{equation}
 whenever $K \leq  N_0$.
 
  In order to estimate the term $\Vert P_M (u |u|^{p-1}) \Vert_{L_t^{\infty} L^{2}}$ in \eqref{eq4.15}, we start with the following decomposition
  \begin{align}
  \begin{split}
  \Vert P_M (u |u|^{p-1}) \Vert_{L^{2}}& = \Vert P_M ( u |P_{\leq N_0} u|^{p-1} - (u |u|^{p-1} - u |P_{\leq N_0} u |^{p-1}) ) \Vert_{L^2}
   \\  & \leq \Vert  P_M (u | P_{\leq N_0} u |^{p-1}) \Vert_{L^2}  +   \Vert P_M (u |u|^{p-1} - u |P_{\leq N_0} u |^{p-1})  \Vert_{L^2}
   \\ & = I +\tilde{I}.
	\end{split}
	\label{eq4.17p} 
  \end{align}
 We then  
   write  $I$
 as a product of two factors decomposed into high-low frequencies around ${M}/{2p}$ and $N_0$. 
 In other words, we get
 \begin{align}
 \begin{split}
I &
= \Vert P_M (( P_{\leq M/2p} u + P_{>M/2p} u ) | P_{\leq N_0} u |^{p-1} ) \Vert_{L^{2}}
\\ & \leq \Vert P_M ( P_{\leq M/2p} u  | P_{\leq N_0} u |^{p-1} ) \Vert_{L^{2}} 
+ \Vert P_M (  P_{> M/2p} u  | P_{\leq N_0} u |^{p-1} )  \Vert_{L^{2}}
\\ &  = I_1 + I_2.
\end{split}
\label{eq4.17}
\end{align} 

We begin with
 \begin{align}
 I_1 = \Vert P_M (( P_{\leq M/2p} u )  | P_{\leq N_0} u |^{p-1} ) \Vert_{L^{2}} \period
 \label{eq4.18}
 \end{align}
 Since $p-1$ is an even integer, we deduce that 
 \begin{equation*}
 \mbox{ supp } \left( \mathcal{F} \left( ( P_{\leq N_0 } u )^{p-1} \right)  \right) \subset \{ \xi: ~|\xi| \leq 2 (p-1) N_0\} .
 \end{equation*}
 As a result, in the case $N_0 \leq M/2p$, we get $I_1 =0$. 
 
 We simply assume that 
 $N_0 > M/2p$, and split $I_1$ in two parts
 \begin{align}
 \begin{split}
 & I_1 \les 
 \Vert 
 P_M (( P_{\leq M/2p} u )  | P_{ \leq {M}/{2p}} u |^{p-1} )
 \Vert_{L^2}
 \\ & \qquad + 
 \Vert
 P_M ( ( P_{\leq M/2p} u ) (   | P_{ \leq N_0} u |^{p-1} - | P_{\leq M/2p} u |^{p-1} ))
 \Vert_{L^2}
 \\ & \qquad = I_{11} + I_{12} \period
 \end{split}
 \label{eq4.19}
 \end{align}
 As noted above, we deduce that $I_{11} =0$. 
 For $I_{12}$, we express the difference by the fundamental theorem of calculus, and we get
 \begin{align}
\left|  | P_{ \leq N_0} u (x) |^{p-1} - | P_{\leq M/2p} u (x) |^{p-1} \right|
  \leq  \left|	P_{{M}/{2p} < \cdot < N_0} u (x) \right| 
  \int_0^1   \left| \lambda P_{\leq N_0} u (x) +  (1- \lambda) P_{\leq M/2p} u (x) \right|^{p-2} d\lambda
 .
 \end{align}
 Next, 
 we recall \eqref{convolution_pk}, and  apply Young's and H\"older's inequality
 %\begin{equation}
%\frac{1}{q } = \frac{1}{2}	+ \frac{2 (p-2)}{d (p-1)} 
%\label{EQ4.54} 
 %\end{equation}
 \begin{align}
 \begin{split}
I_{12}  & \les 
 \Vert    M^d   \check{\phi} (M \cdot)    \Vert_{L^{d/ (d-2)}}
 ~
 \Vert \left|  P_{\leq M/2p} u \right| ( | P_{\leq N_0} u |^{p-2}+ | P_{\leq M/2p} u|^{p-2} ) \Vert_{L^{d/2}}
 ~
  \Vert     P_{{M}/{2p} < \cdot < N_0} u    \Vert_{L^2}
\\ 
 & \les
 M^2 \eta^{p-1}  \Vert     P_{> {M}/{2p}} u    \Vert_{L^2} 
 \end{split}
 \label{eq4.20}
 \end{align}
 where we used \eqref{EQ4.16p} on the last line.

By \eqref{eq4.19}, we then obtain
\begin{equation}
I_{1} \les M^2 \eta^{p-1} \Vert P_{> M/2p} u \Vert_{L^2}.
 \label{EQ4.52}
\end{equation}

 Similarly, we estimate
  \begin{align}
 \begin{split}
 I_2  & \les \Vert    M^d   \check{\phi} (M \cdot)    \Vert_{L^{d/(d-2)}}
 ~ 
 \Vert     P_{> M/2p} u    \Vert_{L^2}
 ~\Vert | P_{\leq  N_0}    u  |^{p-1}   \Vert_{L^{d/2}}
  \\
 & \les
 M^2 \eta^{p-1}  \Vert     P_{> {M}/{2p}} u    \Vert_{L^2} .
 \end{split}
 \label{eq4.21}
 \end{align}
 
 Combining \eqref{EQ4.52} and \eqref{eq4.21}, we obtain
 \begin{align}
 I \les  M^2 \eta^{p-1}  \Vert     P_{> {M}/{2p}} u    \Vert_{L^2} 
 \period
 \label{eq4.21p}
 \end{align}
 
Next, we estimate $\tilde{I}$ in \eqref{eq4.17p}. Similarly, we may control the difference 
 \begin{align*}
\left|  | u |^{p-1} - | P_{\leq N_0 } u |^{p-1} \right|
\leq 	\left| P_{>  N_0} u \right|
  \int_0^1 |  \lambda u + (1- \lambda )  P_{\leq N_0} u |^{p-2}  d\lambda
  .
 \end{align*}
 Using the Young's and H\"older's inequality followed by the Sobolev embedding and Bernstein's inequality 
 at the last step, we obtain
 \begin{align}
 \begin{split}
 \tilde{I}  & \les 
 \Vert    M^d  \check{\phi} (M \cdot)    \Vert_{L^{d/ (d-2)}}
 ~
 \Vert u~ ( | u|^{p-2} +  |P_{\leq N_0} u|^{p-2} ) \Vert_{L^{d/2}}
 ~
  \Vert     P_{ > N_0} u    \Vert_{L^2}
\\ 
 & \les M^2 
 \Vert u \Vert_{L^{d (p-1)/2}}^{p-1} 
 ~ 
 \Vert   P_{>  N_0} u  \Vert_{L^2}
 \\ 
 & \les M^2 
 N_0^{- s_p} 
 ~ 
 \Vert u \Vert_{L^{\infty}_t \dot{H}^{s_p}}^{p} . 
 \end{split}
 \label{eq4.22}
 \end{align}

 Collecting the estimates in \eqref{eq4.17p}, \eqref{eq4.21p}, and \eqref{eq4.22} we then have
 \begin{align}
 \| P_M (u |u|^{p-1}) \|_{L^2} \les M^2 \eta^{p-1} \| P_{> M/2p} u \|_{L^2} + 
 M^2 N_0^{- s_p} 
 \label{eq4.22p}
 \end{align}
which yields by \eqref{eq4.15} 
 \begin{align}
\begin{split}
 \norm{ P_M \left(   \int_0^{\Lambda/M} \frac{e^{i \tau \sqrt{-\Delta}}}{\sqrt{-\Delta}}   |u|^{p-1} u (\tau)~d\tau \right)  }_{\dot{H}^{s_p}} 
& \les \Lambda M^{s_p -2} \Vert P_M  |u|^{p-1} u \Vert_{L_t^{\infty} L^2} \\
&  \les 
\Lambda M^{s_p} \eta^{p-1} \Vert P_{> M/2p} u \Vert_{L_t^{\infty} L^{2}} 
+ \Lambda M^{s_p} N_0^{- s_p}
\period  
\end{split}
\label{eq4.27}
\end{align}
 This completes the proof of Claim~\ref{claim_2}. 
 
 Next, we consider $A_2$ in \eqref{eq4.13}. Recall that 
 \begin{align}
 A_2 = \norm{
 \int_{\Lambda / M}^{\infty}
 \frac{e^{-it \sqrt{- \Delta}}}{\sqrt{- \Delta}}
( 1-\chi ) \left( x/{ct}  \right) P_M \left(  |u|^{p-1} u (t)  \right)~dt
 }_{\dot{H}^{s_p}}
 \period
 \end{align}
First, we move the spacial norm inside the integral and obtain
\begin{align}
A_2 \les \int_{\Lambda / M}^{\infty}
\norm{
( 1-\chi ) \left( x/{ct}  \right) P_M \left(  |u|^{p-1} u (t)  \right)
 }_{\dot{H}^{s_p -1}} dt
 \period
 \label{eq4.29}
\end{align} 
Denote by 
\begin{align}
\tilde{A_2} (t, x) := 
( 1-\chi ) \left( x/{ct}  \right) P_M \left(  |u|^{p-1} u (t)  \right)
\period
 \label{eq4.29p}
\end{align}
Noting the condition \eqref{eq0.1} on $p$ values, we recall that $s_p -1  \in [{(d-4)}/{2}, {(d-2)}/{2})$. We express $s_p -1$ as 
\begin{equation}
s_p -1 = m+s \label{EQ_sp}
\end{equation}
where $m = \left[ s_p -1 \right]$ and $s = s_p - m$. 
By the Leibniz rule, we then have 
\begin{align*}
\norm{  \tilde{A_2} (t,x) }_{\dot{H}^{s_p -1}} \les 
\sum_{i=0}^m J_i (t)
\end{align*}
where each $J_i (t)$ is defined as
\begin{equation}
J_i (t )= \sup_{\mathcal{I}_i} \norm{ \partial_x^{\beta_i}  ( 1-\chi ) \left( x/{ct}  \right)  \partial_x^{\beta_{m_i}} P_M \left(  |u|^{p-1} u (t)  \right) }_{\dot{H}^s}
\label{EQ4.30p}
\end{equation}
over the index set 
\begin{equation*}
\mathcal{I}_i = \{ (\beta_i, \beta_{m_i}): \beta_i \mbox{ and }\beta_{m_i} \mbox{ are multi indices with } |\beta_i|=i, |\beta_{m_i}| = m-i   \}.
\end{equation*}
Firstly, we treat $J_0 (t)$. Letting $\beta_m$ and $\beta_{m+1}$ denote two  multi-indices of order $m$ and $m+1$, respectively,  we bound $J_0 (t)$ from above using interpolation by
\begin{align}
\begin{split}
 &  \sup_{ \beta_{m} } 
 \norm{ ( 1-\chi ) \left( x/{ct}  \right)  \partial_x^{ \beta_{m+1}} P_M \left(  |u|^{p-1} u (t)  \right)}^s_{L^2}
  \norm{  ( 1-\chi ) \left( x/{ct}  \right)  \partial_x^{\beta_{m}} P_M \left(  |u|^{p-1} u (t)  \right)}^{1-s}_{L^2} 
  \\
  &  \quad + \sup_{ \beta_{m} } 
 \norm{  \frac{1}{ct} \chi'  \left( x/{ct}  \right)  \partial_x^{\beta_{m}} P_M \left(  |u|^{p-1} u (t)  \right)}_{L^2}^s 
 \norm{  ( 1-\chi ) \left( x/{ct}  \right)  \partial_x^{\beta_{m}} P_M \left(  |u|^{p-1} u (t)  \right)}^{1-s}_{L^2}  .
\end{split}
\label{EQ4.31p}
\end{align}
Recall that
\begin{equation}
\mbox{supp} (1- \chi ) ( x/ {ct} )  \subset \{ x:  |x| \geq  |ct| \}.
\label{EQ4.32p}
\end{equation}
and 
\begin{equation}
\mbox{supp} ( \chi' ) ( x/ {ct} )  \subset \{ x: |ct| \leq |x| \leq 2 |ct| \}.
\label{EQ4.32.2p}
\end{equation}
Using the radial Sobolev inequality and Bernstein's inequality we bound each factor in \eqref{EQ4.31p} 
as below. We obtain
\begin{align}
\begin{split}
\norm{ ( 1-\chi ) \left( x/{ct}  \right)  \partial_x^{ \beta_{m+1}} P_M \left(  |u|^{p-1} u (t)  \right) }_{L^2} 
& \les
  \frac{1}{|ct|^{2- \frac{2}{d}}} \norm{  |x|^{2- \frac{2}{d}}   \partial_x^{ \beta_{m+1}} P_M \left(  |u|^{p-1} u (t)  \right)}_{L^2} \\
  & \les
  \frac{1}{|ct|^{2- \frac{2}{d}}} \norm{ D^{\frac{2}{d}}  \partial_x^{ \beta_{m+1}} P_M \left(  |u|^{p-1} u (t)  \right)}_{L^{\frac{2d}{d+4}}} \\
  & \les 
 \frac{M^{m+1+ \frac{2}{d}}}{|ct|^{2 - \frac{2}{d}}} \norm{ P_M \left(   |u|^{p-1} u (t)  \right) }_{L^{\frac{2d}{d+4}}} 
 \end{split} \label{EQ4.35p}
\end{align}
and
\begin{align}
\begin{split}
\norm{ ( 1-\chi ) \left( x/{ct}  \right)    \partial_x^{ \beta_{m}} P_M \left(  |u|^{p-1} u (t)  \right) }_{L^2} 
& \les
 \frac{1}{|ct|^{2 - \frac{2}{d}}} \norm{ {|x|^{2 - \frac{2}{d}}} \partial_x^{ \beta_{m}} P_M \left(  |u|^{p-1} u (t)  \right)}_{L^2} \\
  & \les 
  \frac{M^{m+ \frac{2}{d}}}{|ct|^{2- \frac{2}{d}}} \norm{ P_M \left(   |u|^{p-1} u (t)  \right) }_{L^{\frac{2d}{d+4}}} 
 \end{split} 
  \label{EQ4.36p}
\end{align}
 Similarly, we get
 \begin{align}
\begin{split}
\norm{  \chi'  \left( x/{ct}  \right)    \partial_x^{ \beta_{m}} P_M \left(  |u|^{p-1} u (t)  \right) }_{L^2} 
& \les
 \frac{1}{|ct|^{1 - \frac{2}{d}}} \norm{ {|x|^{1 - \frac{2}{d}}} \partial_x^{ \beta_{m}} P_M \left(  |u|^{p-1} u (t)  \right)}_{L^2} \\
  & \les 
  \frac{M^{m+1+ \frac{2}{d}}}{|ct|^{1 - \frac{2}{d}}} \norm{ P_M \left(   |u|^{p-1} u (t)  \right) }_{L^{\frac{2d}{d+4}}} 
 \end{split} 
  \label{EQ4.37p}
\end{align}
 As a result, the right hand side of  \eqref{EQ4.31p} may be estimated from above by
 \begin{align}
 J_0 (t) \les \frac{M^{m+s+ \frac{2}{d}}}{|ct|^{2- \frac{2}{d}}} \norm{ P_M  \left(  |u|^{p-1} u (t)  \right) }_{L^{\frac{2d}{d+4}}}
 . 
 \label{EQ4.34p}
 \end{align}
Integrating $J_0  (t)$ in time and using \eqref{EQ_sp}, we then find 
\begin{equation}
\int_{\Lambda / M}^{\infty} J_0 (t) ~dt \les \frac{M^{s_p}}{\Lambda^{1- \frac{2}{d}}} \norm{ P_M  \left(  |u|^{p-1} u (t)  \right) }_{L^{\infty}_t L^{\frac{2d}{d+4}}} .
\label{EQ4.38p}
\end{equation} 
 
 Next, we estimate $J_i (t)$ by following the same line of arguments. We find that $J_i (t)$ is bounded from above by
 \begin{align}
\begin{split}
 &  \sup_{ \beta_{m_i} } 
\frac{1}{|ct|^{i}}
 \norm{ \chi^{(i)}  \left( x/{ct}  \right) \nabla  \partial_x^{ \beta_{m_i}} P_M \left(  |u|^{p-1} u (t)  \right)}^s_{L^2}
  \norm{   \chi^{(i)} \left( x/{ct}  \right)  \partial_x^{\beta_{m_i}} P_M \left(  |u|^{p-1} u (t)  \right)}^{1-s}_{L^2} 
  \\
  & \quad   + \sup_{ \beta_{m_i} } \frac{1}{|ct|^i}
 \norm{  \frac{1}{ct} \chi^{(i+1)}  \left( x/{ct}  \right)  \partial_x^{\beta_{m_i}} P_M \left(  |u|^{p-1} u (t)  \right)}_{L^2}^s 
 \norm{   \chi^{(i)}  \left( x/{ct}  \right)  \partial_x^{\beta_{m_i}} P_M \left(  |u|^{p-1} u (t)  \right)}^{1-s}_{L^2}  
\end{split}
\label{EQ4.39p}
\end{align}
by interpolation. 
Once again we apply  the radial Sobolev inequality followed by Bernstein's inequality on each factor, and we find
\begin{align}
\begin{split}
\norm{ \chi^{(i)} \left( x/{ct}  \right) \nabla  \partial_x^{ \beta_{m_i}} P_M \left(  |u|^{p-1} u (t)  \right) }_{L^2} 
& \les
 \frac{1}{|ct|} \norm{ |x| \nabla  \partial_x^{ \beta_{m_i}} P_M \left(  |u|^{p-1} u (t)  \right)}_{L^2} \\
  & \les 
 \frac{M^{m-i + 2}}{|ct|} \norm{ P_M \left(   |u|^{p-1} u (t)  \right) }_{L^{\frac{2d}{d+4}}} 
 \end{split} 
 \label{EQ4.40p}
\end{align}
and
\begin{align}
\begin{split}
\norm{ \chi^{(i)} \left( x/{ct}  \right)   \partial_x^{ \beta_{m_i}} P_M \left(  |u|^{p-1} u (t)  \right) }_{L^2} 
& \les
  \frac{1}{|ct|} \norm{  |x|  \partial_x^{ \beta_{m_i}} P_M \left(  |u|^{p-1} u (t)  \right)}_{L^2} \\
  & \les 
 \frac{M^{m-i +1}}{|ct|} \norm{ P_M \left(   |u|^{p-1} u (t)  \right) }_{L^{\frac{2d}{d+4}}} .
 \end{split} 
 \label{EQ4.41p}
\end{align}
Furthermore, we have
\begin{align}
\begin{split}
\norm{ \frac{1}{ct} \chi^{(i+1)} \left( x/{ct}  \right)   \partial_x^{ \beta_{m_i}} P_M \left(  |u|^{p-1} u (t)  \right) }_{L^2} 
& \les
  \frac{1}{|ct|} \norm{    \partial_x^{ \beta_{m_i}} P_M \left(  |u|^{p-1} u (t)  \right)}_{L^2} \\
  & \les 
 \frac{M^{m-i+ 2}}{|ct|} \norm{ P_M \left(   |u|^{p-1} u (t)  \right) }_{L^{\frac{2d}{d+4}}} .
 \end{split} 
 \label{EQ4.42p}
\end{align}
Then, we plug the estimates \eqref{EQ4.40p}--\eqref{EQ4.42p} in \eqref{EQ4.39p}, which yields
\begin{align*}
J_i  (t) \les \frac{ M^{m-i + 1+s} }{|ct|^{i+1}}   \norm{P_{M} \left(   |u|^{p-1} u (t)  \right) }_{L^{\frac{2d}{d+4}}} .
\end{align*}
Next, we integrate the above estimate in time  and use that $\Lambda \geq 1$.  We obtain
\begin{align}
\begin{split}
\int_{\Lambda / M}^{\infty} J_i (t)~dt & \les  \frac{M^{s_p}}{\Lambda^i} 
 \norm{ P_M  \left(  |u|^{p-1} u (t)  \right) }_{L^{\infty}_t L^{\frac{2d}{d+4}}} \\
 & \les   \frac{M^{s_p}}{\Lambda}   \norm{ P_M  \left(  |u|^{p-1} u (t)  \right) }_{L^{\infty}_t L^{\frac{2d}{d+4}}} 
\end{split} 
 \label{EQ4.43p}
\end{align}
Moreover, we may
estimate the term 
$\norm{P_M ( |u|^{p-1} u)}_{L^{2d/{(d+4)}}}$  as done in the proof of Claim~\ref{claim_2}.  We
decompose $u = P_{\leq M/2p} u + P_{> M/2p} u$, and we get
\begin{align*}
\begin{split}
\norm{P_M ( |u|^{p-1} u)}_{L^{ \frac{2d}{d+4}}}
& \leq
 \norm{P_M (( P_{\leq M/2p} u + P_{> M/2p} u ) |u|^{p-1} ) }_{L^{\frac{2d}{d+4}} }
\\  & \leq
 \norm{P_M( ( P_{\leq M/2p} u ) |P_{\leq M/2p} u|^{p-1} )}_{L^{\frac{2d}{d+4}} }
 \\ & \indeq \indeq  + 
\norm{P_M( ( P_{\leq M/2p} u ) ( |u|^{p-1} - |P_{\leq M/2p} u|^{p-1} ) )}_{L^{\frac{2d}{d+4}} }
\\ & \indeq \indeq 
+ \norm{P_M ( (  P_{> M/2p} u ) |u|^{p-1} ) }_{L^{\frac{2d}{d+4}} }
\end{split}
\end{align*}
Note that the first term $ \norm{P_M \left( P_{\leq M/2p} u \right) |P_{\leq M/2p} u|^{p-1} }_{L^{2d/{(d+4)}}}=0$. 
 Meanwhile,  we may control the last two terms by applying Young's inequality and the Sobolev embedding, which leads to the estimate
 \begin{align}
 \begin{split}
 \norm{P_M  \left(  |u|^{p-1} u (t)  \right) }_{L^{\frac{2d}{d+4}}} 
&  \les \norm{M^d  \check{\phi} (M \cdot) }_{L^1} 
 \norm{ P_{> M/2p} u}_{L^2}
 \norm{u}^{p-1}_{\frac{d(p-1)}{2}}
 \\ 
 & \les \norm{P_{> M/2p} u}_{L^2} 
 \end{split}
 \label{EQ4.44p}
 \end{align}
 as done in \eqref{eq4.22}. On the last line above we once again used the Sobolev embedding to get
\begin{align*}
\norm{u}_{L^{d (p-1)/2}} \les 
\norm{u}_{\dot{H}^{s_p}}
\end{align*}
and absorbed the $\dot{H}^{s_p}$ norm of $u$ in the implicit constant. 
 Then,  
 the estimate for $A_2$ follows from collecting the upper bounds in \eqref{EQ4.38p}, \eqref{EQ4.43p}, and \eqref{EQ4.44p}
 \begin{align}
 A_2 \les \frac{M^{s_p}}{\Lambda^{1- \frac{2}{d}}} \norm{P_{> M/2p} u}_{L^{\infty}_t  L^2} . 
 \label{EQ4.45p}
 \end{align}
 
  Combining the estimates for $A_1$ and $A_2$ from Claim~\ref{claim_2}
 and \eqref{EQ4.45p}, and setting $\Lambda^{1- \frac{2}{d}}:= \eta^{-1}$ for $\eta \in (0,1)$ to be fixed below, we have
 \begin{align}
 \begin{split}
 \norm{A}_{\dot{H}^{s_p}} & \les A_1 + A_2
 \\ &  \les 
 \Lambda M^{s_p} \eta^{p-1} \norm{ P_{>M/2p} u}_{L^{\infty}_t L^2} 
 + \Lambda M^{s_p} N_0^{-s_p} 
 \\ & \quad 
 + \frac{M^{s_p}}{\Lambda^{1 - \frac{2}{d}}}  \norm{ P_{>M/2p} u}_{L^{\infty}_t L^2} 
\\ & \les 
\eta M^{s_p} \norm{ P_{>M/2p} u}_{L^{\infty}_t L^2} + \eta^{ - \frac{d}{d-2}} M^{s_p} N_0^{-s_p}
\period
 \end{split}
 \label{EQ4.46p}
 \end{align}
 Note that the right hand side of \eqref{EQ4.46p} also controls the tem $ \Vert {\tilde{A}} \Vert_{\dot{H}^{s_p}}$.

 Next, we consider $\langle A, \tilde{A} + \tilde{B} \rangle_{\dot{H}^{s_p}}$ and $\langle A+B, \tilde{A} \rangle_{\dot{H}^{s_p}}$. 
  We recall once again that $e^{it \sqrt{ - \Delta} } v(t) \rightharpoonup 0$ weakly in 
 $\dot{H}^{s_p} $ as $t \to \pm \infty$ by Lemma~\ref{Lemma_weak_limits}. Therefore, we have
 \begin{align}
 A+B & = P_M v(0) 
 \label{eq4.38} \\
 \tilde{A} + \tilde{B} & = P_M v(0)  .
 \label{eq4.39}
 \end{align}
% weakly in $\dot{H}^{s_p}$.   
 We may then estimate 
 \begin{align}
 \begin{split}
& 
 | \langle A, \tilde{A}+ \tilde{B} \rangle_{\dot{H}^{s_p}} |
 \les 
   \norm{A}_{\dot{H}^{s_p}} 
   \norm{P_M v }_{L_t^{\infty} \dot{H}^{s_p}}
 \\ & \quad \quad \les 
 \left( \eta  M^{s_p} \norm{ P_{>M/2p} u}_{L^{\infty}_t L^2} + \eta^{ - \frac{d}{d- 2}} M^{s_p} N_0^{-s_p} \right)
  M^{s_p}  \norm{P_M v }_{L_t^{\infty} L^2} .
 \end{split}
 \label{eq4.40}
\end{align}  
The same estimate holds for the term $ | \langle A+B , \tilde{A}  \rangle_{\dot{H}^{s_p}} |$ as well. 

Lastly, we estimate $\langle B, \tilde{B} \rangle_{\dot{H^{s_p}}}$ in \eqref{EQ4.2}. To simplify the notation, we denote by
\begin{align*}
f_M (t) = P_M ( |u|^{p-1} u (t) ).
\end{align*}
Note that $f_M$ is a real valued function as the definition of $P_M$ is based on a real radial bump function. 

Recalling \eqref{eq4.14}, we may express
  $\langle B, \tilde{B} \rangle_{\dot{H^{s_p}}}$ as 
\begin{align*}
\int_{\Lambda/M}^{\infty}
 \int_{-\infty}^{- \Lambda /M} 
\left\langle D^{s_p -1} {e^{i t \sqrt{-\Delta}}} 
  \chi \left(x/{ct} \right) f_M (t), D^{s_p -1}  
 {e^{i \tau \sqrt{-\Delta}}} \chi \left(x/{c \tau} \right) f_M (\tau) \right\rangle_{L^2} d\tau dt .
\end{align*}
%First, in the case 
 %$s_p -1$  is a positive integer, we proceed as follows: due to the Huygens Principle, when $c= 1/4$, we have
%\begin{align*}
% \mbox{ supp} \left( {e^{i ( \tau -t) \sqrt{-\Delta}}}  \chi \left(x/{c \tau} \right)  P_M (u | u |^{p-1} (\tau)) \right)  \subset
 %\{x:~ |x| \geq \frac{3}{4} |t- \tau| \}
 %\period
 %\end{align*}
 %Since $t > \Lambda /M >0$ and $\tau < - \Lambda /M <0$, the support of the function on the right side of the bracket 
 %is included in the set $|x| > \frac{3}{4} t$, whereas that of the function $\chi \left(x/ct \right) P_M (  |u|^{p-1} u (t))$
% is in the set $|x| < t/4$. Therefore, in this case $\langle B, \tilde{B} \rangle_{\dot{H^{s_p}}}$ = 0.
Next, we split $\langle B, \tilde{B} \rangle_{\dot{H^{s_p}}}$  into four pieces:
\begin{align}
\begin{split}
  & \int_{\Lambda/M}^{\infty}
    \int_{-\infty}^{- \Lambda /M}
\left\langle
  {e^{i t  \sqrt{-\Delta}}}  \chi \left(x/{ct} \right) D^{s_p -1 } f_M  (t), 
  {e^{i \tau  \sqrt{-\Delta}}}  \chi \left(x/{c \tau} \right)  D^{s_p -1} f_M  (\tau)  \right\rangle_{L^2} d\tau dt  
  \\ & \quad +
  \int_{\Lambda/M}^{\infty}
  \int_{-\infty}^{- \Lambda /M} 
\left\langle
e^{i t \sqrt{- \Delta}} T_{\chi} (t)   \left( f_M ( t)  \right), 
e^{i \tau \sqrt{- \Delta} } T_{\chi}  (\tau) \left(  f_M  (\tau) \right)
\right\rangle_{L^2} d\tau dt  
  \\ & \quad +
  \int_{\Lambda/M}^{\infty}
  \int_{-\infty}^{- \Lambda /M} 
\left\langle
e^{i t \sqrt{- \Delta}}  \chi \left(x/{ct} \right) D^{s_p -1 } f_M (t), 
e^{i \tau \sqrt{- \Delta} } T_{\chi}  (\tau) \left(  f_M (\tau) \right)
\right\rangle_{L^2} d\tau dt  
  \\ & \quad +
  \int_{\Lambda/M}^{\infty}
  \int_{-\infty}^{- \Lambda /M} 
\left\langle
e^{i t \sqrt{- \Delta}} T_{\chi} (t)   \left( f_M  ( t)  \right), 
e^{i \tau \sqrt{- \Delta} } \chi \left( x/{c \tau} \right)  D^{s_p -1} f_M  (\tau ) 
\right\rangle_{L^2} d\tau dt  
\\ &    =
 K_1 + K_2 + K_3 + K_4
\end{split}
 \label{eq4.41}
\end{align}
where  the commutator operator  $T_{\chi} $ evaluated at time $t$ and $\tau$ is defined as follows:
\begin{equation}
T_{\chi} (s)  \left( f_M  ( s) \right) = D^{s_p -1}  \left( \chi \left(x/{cs} \right) f_M (s)  \right)-  \chi \left(x/{cs} \right) D^{s_p -1} f_M  (s) 
.
\label{EQ4.47p}
\end{equation}
As we mentioned in the paragraph following  \eqref{eq4.14}, it is sufficient to estimate the real part of $\langle B, \tilde{B} \rangle_{\dot{H^{s_p}}}$.
To that end,  we will estimate the real part of $K_1$. We first use the commutativity of the multipliers to express $K_1$ as follows:
\begin{align*}
 K_1  = 
  \int_{\Lambda/M}^{\infty}
 \int_{-\infty}^{- \Lambda /M}
\left\langle
    \chi \left(x/{ct} \right) D^{s_p -1 } f_M  (t), 
  {e^{i ( \tau - t)  \sqrt{-\Delta}}}  \chi \left(x/{c \tau} \right)  D^{s_p -1} f_M  (\tau)  \right\rangle_{L^2} d\tau dt .
\end{align*}
Since  $D^{s_p -1}$ and $P_M$ preserve the class of real-valued functions, we get
\begin{align*}
\mbox{Re} (K_1) = 
\int_{\Lambda/M}^{\infty}
 \int_{-\infty}^{- \Lambda /M} 
\left\langle
    \chi \left(x/{ct} \right) D^{s_p -1 } f_M  (t), 
  {\cos (( t - \tau)  \sqrt{-\Delta}) }  \chi \left(x/{c \tau} \right)  D^{s_p -1} f_M  (\tau)  \right\rangle_{L^2} d\tau dt .
\end{align*}
Due to the Huygens Principle, when $c = 1/4$, we have
\begin{align*}
 \mbox{ supp} \left( {\cos {( ( t -\tau) \sqrt{-\Delta}} )}  \chi \left(x/{c \tau} \right) D^{s_p -1}  f_M (\tau ) \right)  \subset
 \left\lbrace~ |x| \geq  |t- \tau| - \frac{|\tau|}{4}\right\rbrace. 
 \end{align*}
 Since $t > \Lambda /M >0$ and $\tau < - \Lambda /M <0$, the support above
 is included in the set $|x| > \frac{3}{4} t$, whereas the support of the function $\chi \left(x/ct \right) f_M (t)$
 is in the set $|x| < t/4$. 
  As a result, we deduce that 
$\mbox{Re}(K_1) = 0$. 
 
 In order to treat the term $K_2$ in \eqref{eq4.41}, we appeal to Lemma~\ref{Com_Leibniz Rule} in the Appendix.  Letting $s_p -1 = m+s$, we write
 \begin{align}
 \norm{T_{\chi} (t)  \left( P_M  ( |u|^{p-1} u ( t) ) \right)}_{L^2} 
\les \sum_{i=1}^{m-1} \tilde{J}_i (t) + \tilde{J}_m (t)
\label{EQ4.6}
 \end{align}
 where 
 \begin{align*}
 \tilde{J}_i  (t) &  = \sup_{I_i}  \frac{1}{|ct|^i} \norm{\chi^{(i)} (r/ct)  D^{m+s-i-j} \partial_x^{\beta_j} P_M (|u|^{p-1} u (t))}_{L^2} \\
 \tilde{J}_m  (t) & = \frac{1}{|ct|^m} \norm{  \chi^{(m)} (r/ct) P_M (|u|^{p-1} u (t))}_{\dot{H}^s}
 \end{align*}
and the index set $I_i$ is defined as in Lemma~\ref{Com_Leibniz Rule}. Note that each $\tilde{J}_i$ may be treated similar to $J_i$ in \eqref{EQ4.30p}. 
By \eqref{EQ4.32p}, the radial Sobolev inequality, and Bernstein's inequality, we estimate
\begin{align}
\begin{split}
 \tilde{J}_i  (t) & \les \sup_{I_i}  \frac{1}{|ct|^{i+1}}  \norm{ |x|  D^{m+s-i-j} \partial_x^{\beta_j} P_M (|u|^{p-1} u (t))}_{L^2} \\
 & \les \sup_{I_i} \frac{M^{s_p -i-j}}{|ct|^{i+1}} \norm{  \partial_x^{\beta_j} P_M (|u|^{p-1} u (t))}_{L^{\frac{2d}{d+4}}} \\
 & \les \frac{M^{s_p -i}}{|ct|^{i+1}}  \norm{ P_M (|u|^{p-1} u (t))}_{L^{\frac{2d}{d+4}}} .
\end{split}
\label{EQ4.48p}
\end{align}
In addition, since $m \geq 1$, we observe that $\tilde{J}_m (t) = J_m (t)$ as  defined in \eqref{EQ4.30p}, which yields
\begin{equation}
\tilde{J}_m (t) \les  \frac{M^{1+s}}{|ct|^{m+1}}   \norm{ P_M (|u|^{p-1} u (t))}_{L^{\frac{2d}{d+4}}} .
\end{equation}
Back to $K_2$, we proceed as in \eqref{EQ4.43p}, and we get
\begin{align}
\begin{split}
K_2 & \les
 \int_{\Lambda / M}^{\infty} 
 \int_{- \infty}^{- \Lambda / M}  
 \left( \sum_{i=1}^{m-1} \tilde{J}_i (t) + \tilde{J}_m (t) \right) 
\left( \sum_{i=1}^{m-1} \tilde{J}_i (\tau) + \tilde{J}_m (\tau) \right) ~ d \tau d t \\
& \les \frac{M^{2 s_p}}{\Lambda^2}  \norm{ P_M (|u|^{p-1} u (t))}^2_{L^{\infty}_t L^{\frac{2d}{d+4}}} . 
\end{split}
\label{EQ4.49p}
\end{align}

For $K_3$ and $K_4$, we utilize the limit in \eqref{eq4.3p}. Beginning with $K_3$, we have
\begin{align}
\begin{split}
 K_3 & = 
\left\langle  
  \int_{\Lambda / M}^{\infty}
  {e^{i t \sqrt{-\Delta}}}
 D^{s_p -1} 
      f_M (t) ~ d t,
      \int_{-\infty}^{- \Lambda /M} 
 {e^{i \tau  \sqrt{-\Delta}}} T_{\chi} (\tau) f_M (\tau) ~d \tau \right\rangle_{L^2} 
 \\ & \quad -
 \left\langle  
    \int_{\Lambda / M}^{\infty}
   {e^{i t \sqrt{-\Delta}}} \left( 1 - \chi (x / c t) \right)    D^{s_p -1} 
      f_M (t) ~ d t,
    \int_{-\infty}^{- \Lambda /M}  
 {e^{i \tau  \sqrt{-\Delta}}} T_{\chi} (\tau) f_M (\tau) ~d \tau \right\rangle_{L^2} .
\end{split}
\label{EQ4.4}
\end{align}
By \eqref{eq4.3p}, the first term above may be estimated as follows:
\begin{align}
\begin{split}
&
\left\langle 
 \int_{\Lambda / M}^{\infty}
  D^{s_p -1} 
 {e^{-i( \frac{\Lambda}{M} - t) \sqrt{-\Delta}}}
      f_M (t) ~ d t,
 \int_{-\infty}^{- \Lambda /M}
 {e^{i (\tau - \frac{\Lambda}{M}) \sqrt{-\Delta}}} T_{\chi} (\tau) f_M (\tau) ~d \tau \right\rangle_{L^2} 
 \\ & 
 \quad = 
 \left\langle 
 D^{s_p} P_M (i v ( \Lambda /M)),
   \int_{-\infty}^{- \Lambda /M}
 {e^{i (\tau - \frac{\Lambda}{M}) \sqrt{-\Delta}}} T_{\chi} (\tau) f_M (\tau) ~d \tau \right\rangle_{L^2} 
 \\ & \quad \leq 
 M^{s_p} \norm{ P_M v}_{L^{\infty}_t L^2}  \int_{-\infty}^{- \Lambda /M} \norm{ T_{\chi } (\tau) f_M (\tau)}_{L^2} d \tau .
\end{split}
\label{EQ4.3}
\end{align}
We control the second term on the right side of \eqref{EQ4.4} by 
\begin{align}
\begin{split}
&
 \int_{\Lambda / M}^{\infty}
 \norm{\left( 1 - \chi (x / c t) \right)    D^{s_p -1} 
      f_M (t) }_{L^2} dt 
      \int_{-\infty}^{- \Lambda /M}       
\norm{ T_{\chi } (\tau) f_M (\tau)}_{L^2} d \tau 
\\ & \quad \les 
\frac{M^{s_p}}{\Lambda^{1 - \frac{2}{d}}}
  \norm{ P_M (|u|^{p-1} u (t))}_{L^{\infty}_t L^{\frac{2d}{d+4}}} 
   \int_{-\infty}^{- \Lambda /M}       
\norm{ T_{\chi } (\tau) f_M (\tau)}_{L^2} d \tau 
\end{split}
\label{EQ4.5}
\end{align}
where we used the radial Sobolev inequality and Bernstein's inequality on the first factor as done in \eqref{EQ4.36p}. Recalling \eqref{EQ4.6}--\eqref{EQ4.49p} and  adding up the bounds in \eqref{EQ4.5} and \eqref{EQ4.3}, we observe that
\begin{align}
|K_3| \les \frac{M^{ 2 s_p}}{\Lambda^{2 - \frac{2}{d}}}
  \norm{ P_M (|u|^{p-1} u (t))}^2_{L^{\infty}_t L^{\frac{2d}{d+4}}} 
  + \frac{M^{ 2 s_p}}{\Lambda}
  \norm{ P_M (|u|^{p-1} u (t))}_{L^{\infty}_t L^{\frac{2d}{d+4}}} \norm{P_M (v)}_{L^{\infty}_t L^2}.  
  \label{EQ4.8}
\end{align}
Analogously, we may bound the $K_4$ term by 
\begin{align}
|K_4| \les \frac{M^{ 2 s_p}}{\Lambda^{2 - \frac{2}{d}}}
  \norm{ P_M (|u|^{p-1} u (t))}^2_{L^{\infty}_t L^{\frac{2d}{d+4}}} 
  + \frac{M^{ 2 s_p}}{\Lambda}
  \norm{ P_M (|u|^{p-1} u (t))}_{L^{\infty}_t L^{\frac{2d}{d+4}}} \norm{P_M (v)}_{L^{\infty}_t L^2}.  
  \label{EQ4.7}
\end{align}
We then combine \eqref{EQ4.49p}, \eqref{EQ4.8}, \eqref{EQ4.7}, and \eqref{EQ4.44p} to estimate 
the real part of $\langle B, \tilde{B} \rangle_{\dot{H^{s_p}}} $. Also, we recall that 
\begin{equation*}
\Lambda > \Lambda^{1 - \frac{2}{d}} := \eta^{-1}
\end{equation*}
for some $\eta \in (0,1)$. Hence, we obtain 
\begin{align}
\begin{split}
|  \mbox{Re}( \langle B, \tilde{B} \rangle_{\dot{H^{s_p}}} )|  &  \les  
\frac{M^{2 s_p}}{\Lambda^{2 - \frac{2}{d}}}  \norm{ P_M (|u|^{p-1} u (t))}^2_{L^{\infty}_t L^{\frac{2d}{d+4}}} \\
& \quad \quad 
 + \frac{M^{ 2 s_p}}{\Lambda}
  \norm{ P_M (|u|^{p-1} u (t))}_{L^{\infty}_t L^{\frac{2d}{d+4}}} \norm{P_M (v)}_{L^{\infty}_t L^2}    \\
& \les   \eta^2 {M^{2 s_p}}  \norm{ P_{>M/2p} u}_{L^{\infty}_t L^2}^2 
   + \eta {M^{ 2 s_p}}
  \norm{ P_{>M/2p} u }_{L^{\infty}_t L^{\frac{2d}{d+4}}} \norm{P_M (v)}_{L^{\infty}_t L^2}.  
\end{split}
\label{EQ4.50p}
\end{align}

Finally, combining \eqref{EQ4.46p}, \eqref{eq4.40}, and \eqref{EQ4.50p}, we arrive at the conclusion that 
 \begin{align}
 \begin{split}
 \langle  P_M v(0), P_M v(0) \rangle_{\dot{H}^{s_p}} 
& \les 
\left( \eta M^{s_p} \norm{ P_{>M/2p} u}_{L^{\infty}_t L^2} + \eta^{ - \frac{d}{d-2}} M^{s_p} N_0^{-s_p} \right)^2
 \\ & \quad  +
\left( \eta  M^{s_p} \norm{ P_{>M/2p} u}_{L^{\infty}_t L^2} + \eta^{ - \frac{d}{d-2}} M^{s_p} N_0^{-s_p} \right) 
M^{s_p} \norm{ P_M v}_{L^{\infty}_t L^2}
 \period
 \end{split}
 \label{eq4.43}
 \end{align}
 As noted at the beginning of the proof, we may utilize the same logic and arguments to estimate the term $\norm{P_M v(t_0)}_{\dot{H}^{s_p}}$ for any $t_0 \in \mathbb{R}$,
 and therefore the upper bound in \eqref{eq4.43} provides us with a uniform control in time. Namely, we have
 \begin{align}
 \begin{split}
 \langle  P_M v(t), P_M v(t) \rangle_{\dot{H}^{s_p}} 
& \les 
\left( \eta M^{s_p} \norm{ P_{>M/2p} u}_{L^{\infty}_t L^2} + \eta^{ - \frac{d}{d-2}} M^{s_p} N_0^{-s_p} \right)^2
 \\ & \quad  +
\left( \eta  M^{s_p} \norm{ P_{>M/2p} u}_{L^{\infty}_t L^2} + \eta^{ - \frac{d}{d-2}} M^{s_p} N_0^{-s_p} \right) 
M^{s_p} \norm{ P_M v}_{L^{\infty}_t L^2}
 \period
 \end{split}
 \label{eq4.44}
 \end{align}
 
 We now go back to the proof of Claim~\ref{claim_1}. Setting $M= 2^j$ for $j \in \mathbb{Z}^{-}$ and recalling \eqref{eq4.10}, we get 
for~$j <0$
\begin{align*}
\begin{split}
a_j^2  & \les \left( \eta M^{s_p} \norm{ P_{>M/2p} u}_{L^{\infty}_t L^2} + \eta^{ - \frac{d}{d-2}} M^{s_p} N_0^{-s_p} \right)^2
 \\ & \quad  +
\left( \eta  M^{s_p} \norm{ P_{>M/2p} u}_{L^{\infty}_t L^2} + \eta^{ - \frac{d}{d-2}} M^{s_p} N_0^{-s_p} \right) 
M^{s_p} \norm{ P_M v}_{L^{\infty}_t L^2}
\\ & \les 
\left( \eta \sum_{i > j- \ell_p} 2^{s_p (j-i)} a_i + 2^{s_p j} \eta^{-\frac{d}{d-2}} N_0^{-s_p}  \right)^2 
\\ & \quad + 
a_j \left(    \sum_{i > j-\ell_p} 2^{s_p (j-i)} a_i  + 2^{s_p j} \eta^{-\frac{d}{d-2}} N_0^{-s_p} \right)
\end{split}
\end{align*}
where $\ell_p = [\log (2p)] +1$. It implies that
\begin{align*}
a_j \les \eta \sum_{i > j-\ell_p} 2^{s_p (j-i)} a_i  + 2^{s_p j} \eta^{-\frac{d}{d-2}} N_0^{-s_p}
\end{align*} 
 for $j <0$.
 For $j >0$, it suffices to use the estimate 
 \begin{align}
 a_j \cong \norm{P_j v}_{L^{\infty}_t \dot{H}^{s_p}} \les 1 \period
\end{align}  
 Recalling the definition of $\beta_k$ in \eqref{eq4.7}-\eqref{eq4.8}, we then deduce that  for $k <0$ 
  \begin{align}
  \begin{split}
  \beta_k &  \les \sum_{j>0} 2^{- |j-k|} + \eta  \left( \sum_{j <0} 2^{- |j-k|}  \left( \sum_{i > j- \ell_p} 2^{- s_p |j-i|} a_i \right) \right)
  \\ & \quad \quad \quad 
  + \eta^{-\frac{d}{d-2}} N_0^{-s_p}  \sum_{j <0} 2^{- |j-k|} 2^{s_p j} 
  \\ & \les 
  \eta \beta_k + \sum_{j>0} 2^{- |j-k|} + \eta^{-\frac{d}{d-2}} N_0^{-s_p}   \sum_{j<0} 2^{- |j-k|} 2^{s_p j} 
  \period
  \end{split}
  \label{eq4.45}
  \end{align}
  Selecting $\eta = 1/4C$, where $C$ is the implicit constant in \eqref{eq4.45}, we absorb the first term 
 on the last line above into the  left hand side, and we obtain 
 \begin{align*}
 \beta_k \les \sum_{j} 2^{-|j-k|} \min (1, 2^{s_p j})
 \end{align*}
 which yields
 \begin{align*}
 \beta_k \les 2^k \inon{ for}~ k <0 \period
 \end{align*}
 As we set $\beta_k =1$ for $k \geq 0$, we conclude that
 $\{ 2^{-3k/4} \beta_k\}_{k \in \mathbb{Z}} \in \ell^2$,
 which
 completes the proof of Proposition~\ref{decay_proposition}.
 \end{proof}
  
\section{Channels of energy for the linear radial wave equation}
\label{sec:linear estimates}
A key ingredient in the rigidity argument is the exterior energy estimates for radial solutions of the wave equation. Next, we state the general form of these estimates 
in all odd dimensions for the radial linear wave equation, taken directly from \cite[Theorem~2]{KLLS15}. The statement of the result below includes minor changes to avoid notational overlaps with previous sections. 
 %is proved in \cite{KLLS15}.
%The next result is taken verbatim from \cite{KLLS15} for convenience. 
%A proof of this result for general odd 
%dimensions may be found in \cite{KLLS15}. 
\begin{Theorem}
\label{thm_linear_estimates}
In any odd dimension $d>0$, every 
radial solution of   
\begin{align}
   \begin{split}
    & \partial^2_{t} V
    - \Delta V 
     = 0   \comma \mbox{~ } x \in \mathbb{R}^d , t \in \mathbb{R}
    \\
   &  \vec{V} (0)
    = (V_0, V_1) \in \dot{H}^1 \times L^2 (\mathbb{R}^d)
    \end{split}
   \label{eq5.1}
  \end{align}
  satisfies the following estimate: 
For every $R>0$, 
 \begin{align}
\begin{split}
\max_{\pm} \lim_{t \to \pm \infty} 
& \int_{r \geq |t|+R} |\nabla_{x,t} V(r,t)|^2 r^{d-1} dr   \\
& \geq \frac{1}{2} \Vert \pi^{\perp}_{R} (V_0, V_1) \Vert^2_{\dot{H}^1 \times L^2 (r \geq R,~r^{d-1} dr)}
\end{split}
\label{eq5.2}
\end{align}
Here
\begin{equation}
P(R) := \mbox{span} \left\lbrace \left( r^{2k_1 - d}, 0 \right), \left( 0, r^{2k_2 -d} \right):~ k_1= 1, 2, \cdots \left[ \frac{d+2}{4} \right];~ k_2= 1, 2, \cdots \left[ \frac{d}{4} \right]  \right\rbrace
\label{eq5.3}
\end{equation}
and $\pi^{\perp}_{R}$ denotes the orthogonal projection onto the complement of the plane $P(R)$ in ${\dot{H}^1 \times L^2 (r \geq R,~r^{d-1} dr)}$. 

The inequality becomes an equality for data of the form $(0, V_1)$ and $(V_0, 0)$. 
Moreover,
the left-hand side 
of \eqref{eq5.2} vanishes  exactly for all data in $P(R)$.  
\end{Theorem}

%\begin{Theorem}
%\label{thm_linear_estimates}
%Let $V$ be a radial solution of 
%\begin{align}
 %  \begin{split}
  %  & \partial^2_{t} V
  %  - \Delta V 
   %  = 0   \comma \mbox{~ } x \in \mathbb{R}^7 , t \in \mathbb{R}
    %\\
   %&  \vec{V} (0)
   % = (V_0, V_1) \in \dot{H}^1 \times L^2 (\mathbb{R}^7)
   %\period 
   % \end{split}
   %\label{eq5.1}
  %\end{align}
%For every $R>0$, 
 %\begin{align}
%\begin{split}
%\max_{\pm} \lim_{t \to \pm \infty} 
%& \int_{r \geq |t|+R} |\nabla_{x,t} V(r,t)|^2 r^6 dr   \\
%& \geq \frac{1}{2} \Vert \pi^{\perp}_{R} (V_0, V_1) \Vert^2_{\dot{H}^1 \times L^2 (r \geq R,~r^6 dr)}
%\end{split}
%\label{eq5.2}
%\end{align}
%where $\pi_R = \mbox{Id} - \pi^{\perp}_R$ is the orthogonal projection onto the plane 
%\begin{equation}
%P(R) := \mbox{span} \left\lbrace \left( 1/r^5, 0 \right), \left(1/r^3, 0 \right), \left( 0, 1/r^5 \right)  \right\rbrace
%\label{eq5.3}
%\end{equation}
%in the space $\dot{H}^1 \times L^2 (r \geq R, r^6 dr)$. 
%\end{Theorem}
%The left hand side 
%of \eqref{eq5.2} vanishes for all data in $P(R)$.  
%Moreover, in \eqref{eq5.2} equality holds  for data of the form $(0, V_1)$ and $(V_0, 0)$. 

\subsection{Algebraic identities for the projection}
\label{sec: orthogonal proj}
In this part,  we discuss the orthogonal projection onto the plane $P(R)$ 
in  $\dot{H}^1 \times L^2 (r\geq \nobreak{R},r^{d-1} dr)$. 
Besides Theorem~\ref{thm_linear_estimates}, 
we recall a few more results  from \cite[Section~4]{KLLS15} below for convenience. 
%Similar to Theorem~\ref{thm_linear_estimates}, 
%the results in this section  may be found in \cite[Section~4]{KLLS}.
We also introduce the notation
$\mathcal{H} = \dot{H}^1 \times L^2 (\mathbb{R}^d )$ that will be commonly used for the rest of the discussion. 
%Here, $B(0,1)$ denotes the unit ball in $\mathbb{R}^7$ and can be replaced with any ball around zero as we are only interested in the 
%behaviour of solutions for $r \geq R_0$, for some fixed radius $R_0 >0$. 

First we concentrate on the coefficients of the orthogonal projection onto the space $P(R)$ of $\vec{u} (t,r)$ with respect to a suitable basis. Below, we recall a few facts 
regarding these projection coefficients such as their explicit formulas and various algebraic identities that highlight the relationship between the exterior energy of the projected solutions $\Vert \pi_R \vec{u} (t) \Vert_{{ \mathcal{H}(r \geq R)}}$ and the projection coefficients. The proofs of these results may be found in \cite[Section~4]{KLLS15}. 

Fixing  $R>0$ in \eqref{eq5.2}--\eqref{eq5.3},  
we denote the  coefficients of the orthogonal projection onto $P(R)$ of $\vec{u} (t,r)$ by $\lambda_i (t, R)$ and $\mu_j (t, R)$, and we have
\begin{align}
\pi_R { \vec{u} (t,r)} & = \left( \sum_{i=1}^{\tilde{k}}  \lambda_i (t,R) r^{2i -d} ,~ \sum_{j=1}^{k} \mu_j (t,R) r^{2j -d}  \right)
\label{eq5.4.1}
\\
\pi^{\perp}_R { \vec{u} (t, r)} & = \left( u (t,r) - \sum_{i=1}^{\tilde{k}}  \lambda_i (t,R) r^{2i -d},~    u_t (t, r) -  \sum_{j=1}^{k} \mu_j (t,R) r^{2j -d}    \right)
.
\label{eq5.4.2}
\end{align}
Here
\begin{equation}
k:= \left[ \frac{d}{4} \right] \comma \tilde{k} := \left[  \frac{d+2}{4} \right] .
\label{EQ5.1p}
\end{equation}

The explicit formulas for  each $\lambda_i$ and $\mu_j$ may be derived via \textit{Cauchy matrix} and some additional tools  from linear algebra as shown in \cite[Section~4]{KLLS15}. We recall these formulas directly below
\begin{align}
\lambda_j (t,R) = \sum_{i=1}^{\tilde{k}} \frac{- R^{d+2-2i -2j }}{(d- 2j)( d+2 -2i -2j)} d_i d_j \int_R^{\infty} u_{r} (t,r) r^{2i-2}~dr \comma 1 \leq j \leq \tilde{k}
\label{EQ5.2p}
\end{align}
and 
\begin{align}
\mu_j (t,R) = \sum_{i=1}^{k} \frac{R^{d-2i -2j}}{d-2i -2j} c_i c_j \int_R^{\infty} u_t (t,r)~r^{2i -1}~dr \comma 1 \leq j \leq k
\label{EQ5.3p}
\end{align}
with $c_j$ and $d_j$ are defined as follows:
\begin{equation}
c_j = \frac{\prod_{1 \leq \ell \leq k}  (d-2j - 2 \ell)}{\prod_{1 \leq \ell \leq k, \ell \neq j}   (2 \ell - 2 j)} \comma 1 \leq j \leq k
\label{EQ5.4p}
\end{equation}
and 
\begin{equation}
d_j = \frac{\prod_{1 \leq \ell \leq \tilde{k}}  (d+2 -2j - 2 \ell)}{\prod_{1 \leq \ell \leq \tilde{k}, \ell \neq j}   (2 \ell - 2 j)} \comma 1 \leq j \leq \tilde{k} .
\label{EQ5.5p}
\end{equation}
Furthermore, similar computations also yield 
\begin{align}
\int_R^{\infty} u_r (t,r) r^{2i -2}~dr & = - \sum_{j=1}^{\tilde{k}} \frac{R^{2i+ 2j -d-2}  (d-2j)}{d+2 - 2i - 2j} \lambda_j (t, R) \comma 1 \leq i \leq \tilde{k}
\label{EQ5.6p} \\
\int_R^{\infty} u_t (t,r) r^{2i -1}~dr & = \sum_{j=1}^{{k}} \frac{R^{2i+ 2j -d} }{d - 2i - 2j} \mu_j (t, R) \comma 1 \leq i \leq {k}.
\label{EQ5.7p}
\end{align}
Next, we state the following identity on $\lambda_j$ which may be obtained by applying  integration by parts to \eqref{EQ5.2p} and using some algebraic identities
satisfied by $d_j$, which leads to 
\begin{equation}
\lambda_j (t,R) = \frac{d_j}{d -2j} \left(  u (R) R^{d-2j} + \sum_{i=1}^{\tilde{k}-1}  \frac{2i ~d_{i+1}~  R^{d-2i -2j}}{(d-2i -2j)} \int_R^{\infty} u(t,r) r^{2i-1}~dr  \right). \label{EQ5.8p}
\end{equation}
Lastly, we recall Lemma~4.5 in \cite{KLLS}, which states the formulas for the norms of the orthogonal projections $\pi_R \vec{u}$ and $\pi^{\perp}_R \vec{u}$. 
\begin{Lemma}
\label{Lem5}
Given $\vec{u} (t) \in \mathcal{H}$, let $\lambda_j (t,R)$ and $\mu_j (t,R)$ 
be the projection coefficients defined as in \eqref{EQ5.2p} and \eqref{EQ5.3p}, respectively. Then, we have
\begin{align}
\Vert \pi_R \vec{u} (t) \Vert^2_{\mathcal{H} (r \geq R)} 
\cong   \sum_{i=1}^{\tilde{k}} \left(   \lambda_i (t,R) R^{2i - \frac{d+2}{2}} \right)^2
+ \sum_{i=1}^k \left(   \mu_i (t,R) R^{2i - \frac{d}{2}} \right)^2
\label{eq5.10}
\end{align}
and 
\begin{align}
\Vert {\pi^{\perp}_R} \vec{u} (t) \Vert^2_{\mathcal{H} (r \geq R)} 
\cong \int_{R}^{\infty} \left(    \sum_{i=1}^{\tilde{k}} \left(  \partial_r \lambda_i (t,r) r^{2i - \frac{d+1}{2}}\right)^2
+   \sum_{i=1	}^{k} \left( \partial_r \mu_i (t,r) r^{2i - \frac{d-1}{2}} \right)^2 \right)~dr 
\period
\label{eq5.11}
\end{align}
Here, $X \cong Y$ means $c_1 Y \leq X \leq c_2 Y$ for some positive constants $c_1$, $c_2$, which depend only in $d$. In particular,
 the constants are independent of $t$ and $R$. 
\end{Lemma}

\section{Singular stationary solutions}
\label{sec:stationary solutions}
In this section, 
we cook up a one-parameter family of singular stationary solutions to the equation \eqref{wave equations} whose asymptotic behaviour 
resemble that of a nonzero solution to \eqref{wave equations} with the compactness property. 
By construction, these singular stationary solutions do not lie in the critical Sobolev space $\dot{H}^{s_p} \times \dot{H}^{s_p -1} (\mathbb{R}^d)$.
We will utilize this fact to close the contradiction argument in the next section. 
\begin{Proposition}
\label{prop_singular_solns}
For any $l \in \mathbb{R}\backslash \{0\}$ there exists a radial $C^2$ solution $Z_l$ of 
\begin{equation}
\Delta Z_l + |Z_l|^{p-1} Z_l =0 \inon{ in }~  \mathbb{R}^d\backslash \{0\}
\label{SS_eq1}
\end{equation}
with the asymptotic behaviour 
\begin{equation}
r^{d-2} Z_l (r) = l + O \left( r^{d - p (d-2)}\right)  \inon{ as }~  r \to \infty
\period
\label{SS_eq2}
\end{equation}
Furthermore, $Z_l \notin L^{q_p} (\mathbb{R}^d)$, where $q_p := \frac{d (p-1)}{2}$ is the critical Sobolev exponent corresponding to $\dot{H}^{s_p}$. This implies that $Z_l \notin \dot{H}^{s_p} (\mathbb{R}^d)$. 
\end{Proposition}

\begin{Remark}
In Proposition~\ref{prop_singular_solns}, we assume that $d \geq 7$ and $p \geq 3$ is an odd integer. The proof below is also valid for $d > 4$ and $p \geq 3$ or 
$d=4$ and $p >3$.  In the energy-critical case, the ground state solution 
\begin{equation*}
W(r) = \frac{1}{\left( 1+ \frac{r^2}{d (d-2)} \right)^{\frac{d}{2} -1}}
\end{equation*}
is a stationary radial solution that belongs to the critical Sobolev space $\dot{H}^1 (\mathbb{R}^d)$, $p = \frac{d+2}{d-2}$, so that $p<3$ for $d\geq 5$. 
\end{Remark}
 
\begin{proof}
Our goal is to show the existence of a radial function  $\varphi \in C^2 (\mathbb{R}^d \backslash \{0\})$ that solves the equation \eqref{SS_eq1}, i.e., 
\begin{align}
- \partial_{rr} \varphi - \frac{d-1}{r} \partial_r \varphi = | \varphi |^{p-1} \varphi \comma r>0. 
\label{SS_eq3}
\end{align}
Setting $\omega (r) = r \varphi (r)$, we note that \eqref{SS_eq3} is equivalent to 
\begin{align*}
- \partial_{rr} \omega - \frac{d-3}{r} \partial_r \omega  + \frac{(d-3) r^{p-3} \omega (r) - |\omega|^{p-1} \omega }{r^{p-1}} = 0. 
\end{align*}
In order to guarantee that $\varphi$ satisfies \eqref{SS_eq2} with $l \in \mathbb{R} \backslash \{0\}$, we will need to impose the condition 
\begin{align}
\lim_{r \to \infty} \omega (r) =0.
\label{SS_eq4}
\end{align}
Next, we introduce the new variables $s= \log (r)$ and $\phi (s) = \omega (r)$ and obtain a non-autonomous differential equation for $\phi$. We get 
\begin{equation}
\ddot{ \phi} + (d-4) \dot{\phi} - (d-3) \phi + |\phi |^{p-1} \phi~ e^{- (p-3) s} =0 .
\label{SS_eq5}
\end{equation}
We may rewrite the equation above as a $2 \times 2$ system by setting 
$x (s) = \phi (s)$ and $y (s) = \dot{\phi} (s)$. We then obtain
\begin{align}
\begin{pmatrix}
\dot{x}  \\
\dot{y} 
\end{pmatrix}
= 
\begin{pmatrix}
y  \\
-(d-4) y + (d-3) x - |x|^{p-1} x ~ e^{- (p-3)s}
\end{pmatrix}
=: F (x, y).
\label{SS_eq6} 
\end{align}
Denote by $\Phi_s ((x_0, y_0))$  the flow 
associated to this system with the initial data $(x_0, y_0)$. Our hypothesis on $p$ implies that $F$ is a $C^2$ vector field,
which, in turn, guarantees the existence of  a unique solution for every initial data $(x_0, y_0) \in \mathbb{R}^2$ that lives in a maximal time interval around the origin. 
Next, we study the phase portraits associated to the system \eqref{SS_eq6} to find a class of initial data that yields a solution with the desired properties. 
%Furthermore, the flow $\Phi_s ((x_0, y_0))$ is as smooth as the nonlinear term in \eqref{SS_eq6}. 
%In particular, as  $p$ is an odd integer,  we guarantee that our solution is $C^2$.  

We note that $(0,0)$ is the only equilibrium point \eqref{SS_eq6}, i.e., $(0,0)$ is the unique point in $\mathbb{R}^2$ that solves  $F(\tilde{x},\tilde{y}) = (0,0)$. 
 Checking the linearized system associated to \eqref{SS_eq6} at $(0,0)$, we find that
\begin{align*}
D F ((0,0)) = \begin{pmatrix}
0 & 1 \\
(d-3) & -(d-4)
\end{pmatrix}
\end{align*}
with eigenvalues $\lambda_1 = 1$ and $\lambda_2 = -(d-3)$. 
Denote by $E_{1}$ and $E_{2}$ the corresponding eigenspaces. More precisely, we have
\begin{align*}
E_{1} = \left\lbrace c \begin{pmatrix}
1 \\ 1
\end{pmatrix}: c \in \mathbb{R}  \right\rbrace
\end{align*}
and 
\begin{align*}
E_{2} = \left\lbrace c \begin{pmatrix}
1 \\ -(d-3)
\end{pmatrix}: c \in \mathbb{R}  \right\rbrace . 
\end{align*}
We then write the formula for solutions to the linearized system 
\begin{align}
\begin{pmatrix}
x_L \\
y_L
\end{pmatrix}
= c_1 \begin{pmatrix}
1 \\ 1
\end{pmatrix}
e^s 
+ 
c_2 \begin{pmatrix}
1 \\ -(d-3)
\end{pmatrix}
e^{- (d-3) s}
\comma 
s \in \mathbb{R}
\label{SS_eq7}
\end{align} 
Note that $E_{2}$ denotes the stable subspace of the space of solutions to the linear system 
given in \eqref{SS_eq7}. 

Due to the hyperbolic nature of the matrix $DF ((0,0))$, the stable manifold theorem (cf. \cite{HSD} and \cite{KH97}) yields a neighborhood $B$ of $(0,0)$ 
on which we have a one-dimensional stable manifold
\begin{equation*}
W^s : = \left\lbrace (x,y) \in B:~  \Phi_s ((x,y)) \to (0,0)  \quad s\to \infty  \right\rbrace
\end{equation*}
that satisfies the following properties:  $W^s$ is 
tangent to the subspace $E_{2}$ at the origin and 
 $B \cap W^s$ is positively invariant, i.e., 
\begin{align*}
\Phi_s (B \cap W^s) \subset B \cap W^s \comma s \in [0, \infty). 
\end{align*}
Moreover, 
 for all $(x_0, y_0) \in B \cap W^s$ we have
\begin{align}
\left| \Phi_s ((x_0, y_0)) - (x_0, -(d-3) x_0) e^{-(d-3) s} \right| = O \left( e^{(-(d-2) p +3)s } \right) \qquad s \to \infty.
\label{SS_eq8}
\end{align}
Recalling that 
\begin{align*}
 \Phi_s ( (x_0, y_0) ) = \left( \phi (s) , \dot{\phi} (s) \right)
 \end{align*}
we may characterize the stable manifold $W^s$ as follows:
\begin{align*}
W^s = \left\lbrace x_0 \in B_1:~ \phi (s) = x_0 e^{-(d-3) s} + O \left( e^{(-(d-2) p +3)s } \right) \comma s \in [0, \infty)  \right\rbrace
\end{align*}
 where the higher order terms are determined in terms of $x_0$. Here, $B_1$ denotes the projection of $B$ onto its first coordinate. 
We also remark that  with the choice of $x_0 =0$, the higher order terms vanish, and the corresponding solution coincides with the equilibrium solution $(0,0)$. 
Next, we fix a nonzero $x_0 \in B_1$ and consider the behaviour of the nonlinear flow. 
Changing back to 
 $\omega (r) = \phi (s)$ and  $r = e^s$
 we obtain 
 \begin{align}
 \omega (r) = \phi (\log (r)) \qquad \omega (1) = x_0.  
 \end{align}
 with the asymptotic estimate 
 \begin{align}
\left|  \omega (r) - x_0 r^{-(d-3)} \right| = O \left( r^{-(d-2)p +3} \right) \qquad r \to \infty.
\label{SS_eq9}
 \end{align}
 Moreover, we claim that $\varphi (r) = \omega (r) / r$ fails to belong to $L^{{d (p-1)}/{2}} (\mathbb{R}^d)$ due to a singularity at the origin. 
To that end, we check the energy identity  corresponding to \eqref{SS_eq5}; multiplying \eqref{SS_eq5} with $\dot{\phi}$ and integrating from $s$ to $\infty$, we obtain
\begin{align*}
\begin{split}
&
(d-4) \int_s^{\infty} \dot{\phi}^2 (\rho) ~d\rho + \int_s^{\infty} \phi^{p} (\rho) e^{- (p-3) \rho} \dot{\phi} (\rho) ~d\rho
\\ & \quad \quad 
= \frac{\dot{\phi}^2 (s)}{2} - (d-3) \phi^2 (s) .
\end{split}
\end{align*}
Applying integration by parts for the second term on the left hand side, we rewrite the identity above as 
\begin{align}
\begin{split}
&
(d-4) \int_s^{\infty} \dot{\phi}^2 (\rho) ~d\rho + a \int_s^{\infty} \phi^{p+1} (\rho) e^{- (p-3) \rho} ~d\rho
\\ & \quad \quad 
= \frac{\dot{\phi}^2 (s)}{2} - (d-3) \phi^2 (s)  + \frac{(\phi (s) e^{- a s })^{p+1}}{p+1}
\end{split}
\label{SS_eq10}
\end{align}
 where $a = \frac{p-3}{p+1} \geq 0$. Since the left hand side is positive, it is impossible to have 
 \begin{equation*}
 (\phi (s), \dot{\phi} (s)) = (0,0)
 \end{equation*}
 for any $s \in (-\infty, \infty)$. 
 On one hand, in the case
 \begin{equation}
 \lim_{s \to - \infty} (\phi (s), \dot{\phi} (s)) = (0,0) , 
 \label{SS_eq11}
 \end{equation}
\eqref{SS_eq10} implies that
\begin{equation}
\lim_{s \to - \infty} \phi (s) e^{- a s} \neq 0.
\label{SS_eq12}
\end{equation}
Switching back to the original variables, $\log (r) = s$, the limit in  \eqref{SS_eq12} yields
\begin{align*}
\lim_{r \to 0} \frac{\omega (r)}{r^a} = \lim_{r \to 0} \frac{\phi (\log (r))}{r^a} \neq 0
\end{align*}
and so, we deduce that 
$ \varphi (r) = \omega (r)/r \notin L^{{d (p-1)}/{2}} (\mathbb{R}^d)$ holds under \eqref{SS_eq11}.  
On the other hand, it is clear that if 
 \begin{equation}
  \lim_{s \to - \infty} (\phi (s), \dot{\phi} (s)) \neq (0,0) , 
 \end{equation}
 then we get
 \begin{equation}
 \lim_{s \to - \infty} \phi (s) \neq 0
 \end{equation}
 which, in turn, implies that $\omega (r)/r \notin L^{{d (p-1)}/{2}} (\mathbb{R}^d)$. 
 
 In conclusion, for every non-zero $x_0$ that belongs to a small neighborhood around $0$, we obtain a solution $\varphi$ that safisfies the assumptions of Proposition~\ref{prop_singular_solns} with $\ell =x_0$. In general,  we may scale our solutions by
  $\varphi_l = {\lambda^{- \frac{2}{p-1}}} \varphi (\frac{r}{\lambda})$.  This way, we obtain a solution of the equation \eqref{SS_eq1} that satisfies \eqref{SS_eq2} with  
  $l = x_0 \lambda^{(d-2)- \frac{2}{p-1}}$.  
%We switch back to the initial setting by $\varphi (r) = \omega (r) / r$. 
%By \eqref{SS_eq10}, we deduce that 
 %$ \varphi \notin L^{{d (p-1)}/{2}} (\mathbb{R}^d)$. Also by \eqref{SS_eq9}, our solution $\varphi (r)$  to \eqref{SS_eq3} satisfies \eqref{SS_eq2}
 %with $x_0 \neq 0$. As $x_0$ might vary in a bounded neighborhood around $0$, we may scale our solutions by
 % $\varphi_l = {\lambda^{- \frac{2}{p-1}}} \varphi (\frac{r}{\lambda})$. 
 % This way, we obtain a solution of the equation \eqref{SS_eq1} that satisfies \eqref{SS_eq2} with  
  %$l = x_0 \lambda^{(d-2)- \frac{2}{p-1}}$.  
 
\end{proof}

\section{Rigidity argument}
\label{sec: rigidity}
In this section, we prove Proposition~\ref{Proposition1}. The proof proceeds in three main steps and follows the line of arguments
presented in \cite{KLLS} for solutions to exterior wave maps in all equivariance classes. 

First, we state an important outcome of the decay results obtained in Section~\ref{sec:Decay}. We show that boundedness in 
$\dot{H}^{{3}/{4}} \times \dot{H}^{-{1}/{4}} (\mathbb{R}^d)$ combined with the pre-compactness in $\dot{H}^{s_p} \times \dot{H}^{s_{p-1}} (\mathbb{R}^d)$ yields  
pre-compactness in the energy space. 
\begin{Corollary}
\label{cor8.1}
Let $\vec{u} (t)$ be a solution to \eqref{wave equations} as in Proposition~\ref{Proposition1}.
Then, we have $\vec{u} (t) \in \dot{H}^1 \times L^2 (\mathbb{R}^d)$ for all $t \in \mathbb{R}$. Moreover, the trajectory 
\begin{equation}
K_1 = \{ \vec{u}(t): \inon{$t \in \mathbb{R}$}\}
\label{eq7.1}
\end{equation}
is pre-compact in $\dot{H}^1 \times L^2 (\mathbb{R}^d)$. As a result, we have for all $R >0$
\begin{equation}
\limsup_{t \to +\infty} \Vert \vec{u} (t) \Vert_{\mathcal{H} (r \geq R +|t|)} = \limsup_{t \to - \infty} \Vert \vec{u} (t) \Vert_{\mathcal{H} (r \geq R +|t|)}   =  0
\label{eq7.2}
\end{equation}
\end{Corollary}

\begin{proof}
The proof of Corollary~\ref{cor8.1} is similar to the proof of Lemma~{6.1} in \cite{DL15}. 
We first prove that the trajectory $K_1$ is pre-compact in $\dot{H}^1 \times L^2 (\mathbb{R}^d)$. 
We take a sequence $\{ t_n \} \subset \mathbb{R}$ and show that $\{\vec{u} (t_n)\}$ has a convergent sequence. 
The argument below shows that it suffices to consider $t_n \to \pm \infty$. Without loss of generality, 
we let $t_n \to \infty$. 

Firstly, we consider the case $\{ \lambda (t_n)\}$ remains bounded, which implies that $\{\lambda (t_n)\}$ is a pre-compact sequence. 
Note that in this case 
the sequence $\{\vec{u} (t_n)\}$
is pre-compact in $\dot{H}^{s_p} \times \dot{H}^{s_{p-1}}$ if and only if the sequence 
\begin{equation}
\left\lbrace \left(\frac{1}{\lambda (t_n)^{\frac{2}{p-1}}} u \left( \frac{x}{\lambda (t_n)}, t_n \right), 
\frac{1}{\lambda (t_n)^{\frac{2}{p-1} +1}} \partial_t u \left(\frac{x}{\lambda (t_n)}, t_n \right)   \right) \right\rbrace
\label{eq6.3}
\end{equation}
is pre-compact  in $\dot{H}^{s_p} \times \dot{H}^{s_{p-1}}$ , where the latter fact is guaranteed by hypothesis. 

Using interpolation we control the norm in energy space  by
\begin{align*}
\norm{\vec{u} (t_n) - \vec{u} (t_m)}_{\dot{H}^1 \times \dot{L^2}} 
 \les 
\norm{\vec{u} (t_n) - \vec{u} (t_m)}_{\dot{H}^{\frac{3}{4}} \times \dot{H}^{- \frac{1}{4}}}^{\alpha_p}
\norm{\vec{u} (t_n) - \vec{u} (t_m)}_{\dot{H}^{s_p} \times \dot{H}^{s_{p-1}}}^{1 -\alpha_p}
\end{align*}
where $\alpha_p  \in (0,1)$.  Then, by Proposition~\ref{decay_proposition} we get
\begin{align*}
\norm{\vec{u} (t_n) - \vec{u} (t_m)}_{\dot{H}^1 \times \dot{L^2}} 
 \les  
 \norm{\vec{u} (t_n) - \vec{u} (t_m)}_{\dot{H}^{s_p} \times \dot{H}^{s_{p-1}}}^{1 -\alpha_p}
 \period
\end{align*}
Since the sequence on the right hand side is precompact as discussed above, so is the sequence 
on the left hand side.

Next, we consider the case $\lambda (t_n) \to \infty$. We will show that in this case 
\begin{equation}
\vec{u} (t_n) \to 0 \quad \mbox{in } \dot{H}^1 \times L^2
\period
\label{eq6.1} 
\end{equation}
Let $\eta>0$ be given. By Lemma~\ref{Lemma_tails} there is $c(\eta)>0$ so that
\begin{equation}
\int_{|\xi| \leq c (\eta) \lambda (t)} |\xi|^{2 s_p} |\hat{u} (\xi, t)|^2 d\xi \leq \eta
\label{eq6.4}
\end{equation} 
for all $t \in \mathbb{R}$.
Then we get
\begin{align}
\begin{split}
\norm{ u (t_n)}^2_{\dot{H}^1} 
& = \int_{|\xi| \leq c (\eta) \lambda (t_n)} |\xi|^{2} |\hat{u} (\xi, t_n)|^2 d\xi
+ \int_{|\xi| \geq c (\eta) \lambda (t_n)} |\xi|^{2} |\hat{u} (\xi, t_n)|^2 d\xi  
\\
& \leq 
\int_{|\xi| \leq c (\eta) \lambda (t_n)} |\xi|^{2} |\hat{u} (\xi, t_n)|^2 d\xi
+\left( c(\eta) \lambda (t_n) \right)^{2 -2 s_p } \norm{u (t_n)}^2_{\dot{H}^{s_p}}. 
\end{split}
\label{eq6.2}
\end{align}
Recall that the norm $\Vert u (t_n) \Vert_{\dot{H}^{s_p}}$ is invariant when scaled by $\lambda (t_n)$ as in \eqref{eq6.3}, and therefore is bounded 
in $n$. 
Combined with our assumption that $\lambda (t_n) \to \infty$, we find that the  second term in \eqref{eq6.2}  approaches zero as $n \to \infty$. 
The first term is controlled by interpolation as done above. We get
\begin{align*}
&
\int_{|\xi| \leq c (\eta) \lambda (t_n)} |\xi|^{2} |\hat{u} (\xi, t_n)|^2 d\xi
\\ &  \indeq  \les
\left( \int_{|\xi| \leq c (\eta) \lambda (t_n)} |\xi|^{\frac{3}{2}} |\hat{u} (\xi, t_n)|^2 d\xi \right)^{\alpha_p}
\left( \int_{|\xi| \leq c (\eta) \lambda (t_n)} |\xi|^{2 s_p} |\hat{u} (\xi, t_n)|^2 d\xi   \right)^{1- \alpha_p}
\\ & \indeq \les 
\eta^{1 - \alpha_p} \norm{u (t_n)}_{\dot{H}^{3/4}}^{\alpha_p} \les \eta^{1 - \alpha_p}
\end{align*}
where we used \eqref{eq6.4} and Proposition~\ref{decay_proposition} in the last line. 
As a result, $u (t_n)$ tends to zero in $\dot{H}^1$. Using the same line of arguments, we may also get
$\partial_t u(t_n) \to 0$ in $L^2$, which completes the proof of the first claim that the set $K_1$ is precompact
in $\dot{H}^1 \times L^2$.  

We note that the pre-compactness of $K_1$ implies that 
\begin{equation*}
\norm{\vec{u} (t) }_{\dot{H}^1 \times L^2 (r \geq R)} \to 0 \quad \mbox{as } R \to \infty
\end{equation*}
uniformly in $t \in \mathbb{R}$. Therefore, it
leads to the fact that the energy of $\vec{u} (t)$ 
  on the exterior cone $\{r > R + |t| \}$ vanishes as $t \to \pm \infty$. 
\end{proof}

\subsection{Step~1} 
Let $\vec{u} (t)$ satisfy the assumptions of Proposition~\ref{Proposition1}.
The goal of this part is to estimate $\pi^{\perp}_R \vec{u} (t)$ in $\mathcal{H} (r \geq R)$. 
We combine the linear estimates in Theorem~\ref{thm_linear_estimates} with Corollary~\ref{cor8.1}
to obtain the following result. 

\begin{Lemma} 
\label{Lem4}
There exists $R_0 >0$ such that for all $R \geq R_0$ and for all $t \in \mathbb{R}$ we have
\begin{align}
\Vert \pi^{\perp}_R \vec{u}(t) \Vert^2_{\mathcal{H} (r \geq R)}  
\les {R^{- 2 \beta }} \Vert \pi_R \vec{u}(t) \Vert^{2p}_{\mathcal{H} (r \geq R)}  
\label{eq7.3}
\end{align}
where $\beta = \frac{d-2}{2} \left( p- \frac{d+2}{d-2} \right)$ and  the projections $\pi$, $\pi^{\perp}$ are as in Section~\ref{sec:linear estimates}. 
\end{Lemma}

First, we prove a preliminary result concerning a Cauchy problem for finite energy solutions 
away from the origin. 

\begin{Notation}
Let $I \subset \mathbb{R}$ be an interval with $0 \in I$. For  $q \in [1, \infty]$, denote by 
$L^q_{I} := L^q \left( \mathbb{R}^d \times I\right)$. 
\end{Notation}

First, we take a radial cut-off function $\chi \in C^{\infty} (\mathbb{R}^{d})$ so that 
\begin{align}
\chi (r) = 
\begin{cases}
1 & \inon{if  }~ r \geq 1/2, \\
0 & \inon{if  }~ r \leq 1/4 \period
\end{cases}
\label{eq7.4}
\end{align}
For $r_0 > 0$, denote by $\chi_{r_0} (r) = \chi \left( r/r_0 \right)$ and consider the Cauchy problem
\begin{align}
   \begin{split}
    & \partial^2_{t} h
    - \Delta h
     =  \left| V + \chi_{r_0}  h \right|^{p-1} \left( V + \chi_{r_0}  h  \right)  - |V|^{p-1} V
     \inon{ in }~ \mathbb{R}^d \times I ,
    \\
   &   \evaluat[\big]{(h, \partial_{t} h)}_{t=0}
    = (h_0, h_1) \in \mathcal{H} \period
    \end{split}
   \label{eq7.5}
  \end{align}

We introduce the following exponents in order to simplify the notation in Lemma~\ref{Lem17}.  Denote by
\begin{align}
q_0 = \frac{d+1}{2} (p-1) \comma q_1 = \frac{2 (d+1)}{d-2} \comma q_2 = \frac{2 (d+1)}{d-1} .
\label{EQ7.1p}
\end{align}  
  
\begin{Lemma} 
\label{Lem17}
There exists $\delta_0 > 0$  satisfying the following property:
let $V \in L^{q_0}_{I}$ be radial in the $x$~variable satisfying 
\begin{align}
\Vert D^{1/2} V \Vert_{L^{q_2}_I} \leq \delta_0 r_0^{\beta / (p-1)}~\inon{ and } \inon{~}  \Vert V \Vert_{L^{ q_0}_I} \leq \delta_0 
\label{eq7.6}
\end{align}
where $\beta = \frac{d-2}{2} \left( p- \frac{d+2}{d-2} \right)$.
Furthermore,  $(h_0, h_1) \in \mathcal{H}$ be radial functions with 
\begin{equation}
\Vert (h_0, h_1) \Vert_{\mathcal{H}} \leq \delta_0 r_0^{\beta / (p-1)} \period
\label{eq7.7}
\end{equation}
Then, the Cauchy problem \eqref{eq7.5} is well-posed on the interval $I$, and we have
\begin{align}
\sup_{t \in I} \Vert h (t) - S(t) (h_0, h_1) \Vert_{\mathcal{H}} \leq 
\frac{1}{2^d \cdot 100} \Vert (h_0, h_1) \Vert_{\mathcal{H}} \period
\label{eq7.8}
\end{align}
Moreover, if $V =0$, we may take $I = \mathbb{R}$ and we obtain
\begin{align}
\sup_{t \in \mathbb{R}} \Vert h (t) - S(t) (h_0, h_1) \Vert_{\mathcal{H}} \les
\frac{1}{r_0^{\beta }} \Vert (h_0, h_1) \Vert_{\mathcal{H}}^{p} \period
\label{eq7.9}
\end{align}
\end{Lemma}

\begin{proof}
Let $F_{V} (h) =  \left| V + \chi_{r_0}  h \right|^{p-1} \left( V + \chi_{r_0}  h  \right)  - |V|^{p-1} V$. We apply a fixed point argument to show that 
the formula 
\begin{align*}
h(t) = S(t) (h_0, h_1) + \int_0^t \frac{\sin ((t-s) \sqrt{- \Delta}) }{\sqrt{- \Delta}} F_V (h(s))~ds
\end{align*}
holds for $t \in I$. 
We define the norm
\begin{align*}
\| h \|_{\mathcal{S}} := \norm{h}_{L^{q_1}_I} + {\| D^{1/2} h\|}_{L^{q_2}_I} + \sup_{t \in I} \norm{(h, h_t)}_{\mathcal{H}}
\end{align*}
where the exponents $q_1$, $q_2$ are as in \eqref{EQ7.1p} and for $\alpha >0$ we denote by
\begin{align*}
B_{\alpha} := \{h \in L^{q_1}_{I}: h \mbox{ is radial},~ \| h \|_{\mathcal{S}}  \leq \alpha \}. 
\end{align*}
Now, for $v \in B_{\alpha}$ we set
\begin{align*}
\Phi (v) (t) := S(t) (h_0, h_1) +   \int_0^t \frac{\sin ((t-s) \sqrt{- \Delta}) }{\sqrt{- \Delta}} F_V (h(s))~ds .
\end{align*}
We will show that if \eqref{eq7.6}--\eqref{eq7.7} hold, we can set $\alpha >0$ small enough that 
$\Phi$ is a contraction on $B_{\alpha}$. 

By Strichartz estimates, 
\begin{align}
\| \Phi (v) \|_{\mathcal{S}}  \leq C \left( \norm{(h_0, h_1)}_{\mathcal{H}} + \|{ D^{1/2} F_V (v)}\|_{L^{q'_2}_I}  \right)
\label{eq7.10}
\end{align}
where $q'_2 = 2 (d+1) / (d+3)$. 
We estimate the second term above using the chain rule for fractional derivatives (cf. \cite[Lemma~2.2]{KM11}). First, we let 
$G(h) = |h|^{p-1} h$ so that we may write 
\begin{equation*}
F_V (v) = G (V + \chi_{r_0} v) - G(V). 
\end{equation*}
We get
\begin{align}
\begin{split}
& \| D^{1/2} F_V (v) \|_{L^{q'_2}_I} 
=
\| D^{1/2} (G (V + \chi_{r_0} v) - G(V)) \|_{L^{q'_2}_I} 
\\ & 
 \indeq \indeq \indeq \leq
C \left( \|  G' (V + \chi_{r_0} v) \|_{L^{(d+1)/2}_I}    +   \| G' (V) \|_{L^{(d+1)/2}_I} \right) \| D^{1/2} (\chi_{r_0} v) \|_{L^{q_2}_I}
\\ & 
 \indeq \indeq \indeq +
 C \left( \|  G'' (V + \chi_{r_0} v) \|_{L^{{q_0}/{(p-2)} }_I}    +   \| G'' (V) \|_{L^{{q_0}/{(p-2)}}_I} \right)  \| \chi_{r_0} v \|_{L^{q_0}_I}
 \\ &  \indeq \indeq \indeq \indeq \indeq  \indeqtimes
 \left( \| D^{1/2}  (V + \chi_{r_0} v)\|_{L^{q_2}_I}  + \| D^{1/2} (V) \|_{L^{q_2}_I}  \right).
\end{split}
\label{eq7.11}
\end{align}
Note that we may express
\begin{align}
\begin{split}
\| \chi_{r_0} v \|_{L^{q_0}_I}^{p-1}
&  =
\left(  \int_I \int_{r_0}^{\infty}  |v(r) |^{\frac{d+1}{2} (p-1)} r^{d-1} ~dr~dt  \right)^{2/(d+1)}
\\ &
= \left(  \int_I \int_{r_0}^{\infty} \frac{ | r^{(d-2)/2} v(r) |^{q_0 - q_1} }{r^{\frac{d-2}{2}  (q_0 -q_1)}}  |v(r) |^{q_1} r^{d-1} ~dr~dt  \right)^{2/ (d+1)} 
\end{split}
\label{eq7.12}
\end{align}
where
\begin{align*}
\frac{d-2}{2} (q_0 - q_1) = (d+1) \left( \frac{d-2}{4} (p-1) -1  \right).
\end{align*}
 Recalling the Sobolev inequality for  radial functions $f \in \dot{H}^1 (\mathbb{R}^d)$
 \begin{equation}
 \| r^{(d-2)/2} f \|_{L^{\infty} (\mathbb{R}^d)}  \les \| f \|_{\dot{H}^1 (\mathbb{R}^d)},
\end{equation} 
we estimate \eqref{eq7.12} from above by
\begin{align}
\begin{split}
\| \chi_{r_0} v \|_{L^{q_0}_I}^{p-1}
&  \les
{r_0^{ - \frac{(d-2)}{2} (p-1) + 2 } } \norm{v}^{(p- (d+2)/ (d-2))}_{L^{\infty} (I, \dot{H}^1 )} \norm{v}_{L^{q_1}_I}^{4/(d-2)} 
\\ & \les
 {r_0^{ - \frac{(d-2)}{2} (p-1) + 2 } } \| v \|_{\mathcal{S}}^{p-1}.
\end{split}
\label{eq7.13}
\end{align}
Therefore, we may control the right hand side of \eqref{eq7.11} by 
\begin{align}
\begin{split}
&
C \left( \| V \|^{p-1}_{L^{ q_0}_I} +  r_0^{- \frac{(d-2)}{2} (p-1) + 2 }  {\| v \|^{p-1}_{\mathcal{S}}}  \right) \|  v \|_{\mathcal{S}} \\
& +  C \left(    \| V \|^{p-2}_{L^{ q_0}_I} +  r_0^{- \frac{(d-2)}{2} (p-2) + \frac{ 2 (p-2)}{p-1} }  {\| v \|^{p-2}_{\mathcal{S}}}    \right)   { r_0^{ - (d-2)/2 + 2/(p-1)}}  \|  v \|_{\mathcal{S}} 
\\ & 
 \indeq \indeq \indeq \indeq \indeq  \indeqtimes
 \left(   \| D^{1/2}  V \|_{L^{q_2}_I} + \| v \|_{\mathcal{S}} \right)
 \\ & \leq C 
 \| v \|_{\mathcal{S}}  \left(  \| V \|^{p-1}_{L^{ q_0}_I}  + r_0^{  - \frac{(d-2)}{2} (p-1) +2}   \| D^{1/2}  V \|_{L^{q_2}_I}^{p-1}   + r_0^{  - \frac{(d-2)}{2} (p-1) +2}   {\| v \|^{p-1}_{\mathcal{S}}}   \right) 
\end{split}
\label{eq7.14}
\end{align}
where we applied Young's inequality on the second term in order to obtain the upper bound on the last line.

Combining the bound in \eqref{eq7.14} with \eqref{eq7.10} we get 
\begin{align*}
\| \Phi (v) \|_{\mathcal{S}} \leq C_0 \norm{(h_0, h_1)}_{\mathcal{H}}
          + C_0 \alpha (    \| V \|^{p-1}_{L^{ q_0}_I}  + r_0^{  -  \beta}   \| D^{1/2}  V \|_{L^{q_2}_I}^{p-1}   + r_0^{  - \beta }  \alpha^{p-1}   ) 
\end{align*}
%\label{eq7.15}
for some $C_0 > 0$. We set 
\begin{equation}
\alpha = 2 C_0 \norm{(h_0, h_1)}_{\mathcal{H}} \leq 2 C_0 \delta_0 r_0^{\beta/ (p-1)}.  
\label{eq7.16}
\end{equation}
By \eqref{eq7.6}--\eqref{eq7.7}, we then obtain
\begin{align}
\| \Phi (v) \|_{\mathcal{S}} \leq \frac{\alpha}{2} +  \alpha (2 C_0 \delta_0^{p-1} + 2^{p-1} C_0^{p} \delta_0^{p-1}  ) .
\end{align}
Selecting $\delta_0  > 0 $ sufficiently small we guarantee that $\Phi (v) \in B_{\alpha}$ for every $v \in B_{\alpha}$. 

   The contraction property may be proved using the same arguments. For each $v, \omega \in B_{\alpha}$ the difference 
\begin{align*}
   \| D^{1/2} \left( F_V( V+ \chi_{r_0} v )  - F_V (V + \chi_{r_0} \omega) \right) \|_{L^{q'_2}_I}
   \end{align*}
is estimated by using once again the chain rule for fractional derivatives. Namely, we have
\begin{align}
\begin{split}
& \| D^{1/2}  \left( F_V( V+ \chi_{r_0} v )  - F_V (V + \chi_{r_0} \omega) \right)  \|_{L^{q'_2}_I} 
=
\| D^{1/2} (G (V + \chi_{r_0} v) - G(V + \chi_{r_0} \omega)) \|_{L^{q'_2}_I} 
\\ & 
 \indeq \leq 
C \left( \|  G' (V + \chi_{r_0} v) \|_{L^{(d+1)/2}_I}    +   \| G' (V+ \chi_{r_0} \omega ) \|_{L^{(d+1)/2}_I} \right) \| D^{1/2} (\chi_{r_0} v - \chi_{r_0} \omega ) \|_{L^{q_2}_I}
\\ & 
 \indeq \indeq  +
 C \left( \|  G'' (V + \chi_{r_0} v) \|_{L^{ q_0 / {(p-2)} }_I}    +   \| G'' (V+ \chi_{r_0} \omega) \|_{L^{ q_0 /{(p-2)}}_I} \right)  \| \chi_{r_0} (v - \omega) \|_{L^{q_0}_I}
 \\ &  \indeq \indeq \indeq \indeq \indeq  \indeqtimes
 \left( \| D^{1/2}  (V + \chi_{r_0} v)\|_{L^{q_2}_I}  + \| D^{1/2} (V + \chi_{r_0} \omega ) \|_{L^{q_2}_I}  \right)
 \\ & 
 \indeq  \leq 
 C \| (v - \omega) \|_{\mathcal{S}} \left( \| V \|^{p-1}_{L^{q_0}_I}  +  r_0^{  - \frac{(d-2)}{2} (p-1) +2} (  \| D^{1/2}  V \|_{L^{q_2}_I}^{p-1}   
 +    \| v \|^{p-1}_{\mathcal{S}} + \| \omega \|^{p-1}_{\mathcal{S}} ) \right).
\end{split}
\end{align}

Therefore, we obtain $\Phi (h(t)) = h(t)$. Moreover, by Strichartz estimates and  \eqref{eq7.14} we get
\begin{align}
\| h - S(t) (h_0, h_1) \|_{\mathcal{S}} \leq C_0~ \alpha \left(  \| V \|^{p-1}_{L^{ q_0}_I}  +  r_0^{  - \frac{(d-2)}{2} (p-1) +2} (  \| D^{1/2}  V \|_{L^{q_2}_I}^{p-1}  + \alpha^{p-1} ) \right)
\label{eq7.17}
\end{align}    
which implies \eqref{eq7.8}  with our choice of $\alpha$ in \eqref{eq7.16} provided that $\delta >0$ is sufficiently small. Similarly, in the case $V=0$, 
the inequality \eqref{eq7.17} yields \eqref{eq7.9}.
   \end{proof}

Going back to the proof of Lemma~\ref{Lem4}, we follow the ideas demonstrated in \cite[Prop.~5.3]{KLLS} 
\begin{proof}[Proof of Lemma~\ref{Lem4}] 
First we prove the inequality \eqref{eq7.3} for $t=0$. We take $R>0$ and denote the truncated initial data by
\begin{equation}
\vec{u}_R (0) : = (u_{0,R}, u_{1,R})
\label{eq7.18} 
\end{equation}
where 
\begin{align*}
u_{0,R} := 
\begin{cases}
u_0 (r) & \mbox{for}~ r\geq R \\
u_0 (R) & \mbox{for}~r \leq R
\end{cases}
\end{align*}
and
\begin{align*}
u_{1,R} := 
\begin{cases}
u_1 (r) & \mbox{for}~ r\geq R \\
0 & \mbox{for}~r \leq R.
\end{cases}
\end{align*}
Observe that 
\begin{equation}
\| \vec{u}_R (0) \|_{\mathcal{H}} \leq \|\vec{u} (0) \|_{\mathcal{H} (r \geq R)}
\end{equation}
which implies that we may select $R_0 >0$ sufficiently large so that for all $R \geq R_0$ the truncated initial data is small in $\dot{H}^1 \times L^2$. 
In particular, fixing $\delta < \min (\delta_0, 1)$, where $\delta_0$ denotes the positive constant given in Lemma~\ref{Lem17}, we may guarantee that
\begin{align*}
\| \vec{u}_R (0) \|_{\mathcal{H}} \leq \delta
\end{align*}
for all $R > R_0$. 

Let $\vec{u}_R (t)$ denote the solution to the equation 
\begin{align*}
  & \partial^2_{t} h
    - \Delta h
     =  \chi_{R}   \left|   h \right|^{p-1}   h   
     \inon{ in }~ \mathbb{R}^d \times I ,
    \\
   &   \evaluat[\big]{(h, \partial_{t} h)}_{t=0}
    = (u_{0,R}, u_{1,R}) \in \mathcal{H} 
\end{align*}
given by \eqref{eq7.5} in the case $V=0$. Note that in this case the solution $\vec{u}_R (t)$ exists for all 
$t \in \mathbb{R}$. Moreover, by finite speed of propagation, 
\begin{align}
\vec{u}_R (t,r) = \vec{u} (t,r)   
\label{eq7.21}
\end{align}
for all $t \in \mathbb{R}$ and $r \geq R+ |t| $. 
Next, we define $\vec{u}_{R, L} (t) = S(t) \vec{u}_R (0)$ and note that 
\begin{align}
\begin{split}
&\| \vec{u} (t) \|_{\mathcal{H} (r \geq R + |t| )} = \| \vec{u}_R (t) \|_{\mathcal{H} (r \geq R + |t| )}  \\
& \indeq \indeq  \geq  \| \vec{u}_{R, L} (t) \|_{\mathcal{H} (r \geq R + |t| )} - \| \vec{u}_R (t) - \vec{u}_{R, L} (t) \|_{\mathcal{H} (r \geq R + |t| )}.
 \end{split}
 \label{eq7.19}
\end{align}
By Lemma~\ref{Lem17}, 
\begin{align*}
\sup_{t \in \mathbb{R}} \| \vec{u}_R (t) - \vec{u}_{R,L} (t) \|_{\mathcal{H}} \leq \frac{C_0}{R^{\beta }} \| \vec{u}_R (0) \|^{p}_{\mathcal{H}} .
\end{align*}
Combining this estimate with \eqref{eq7.19} we have
\begin{align}
\| \vec{u} (t) \|_{\mathcal{H} (r \geq R + |t|)}  \geq \| \vec{u}_{R, L} (t) \|_{\mathcal{H} (r \geq R + |t| )}  
- \frac{C_0}{R^{  \beta }} \| \vec{u} (0) \|^{p}_{\mathcal{H} (r \geq R)}.
\label{eq7.20}
\end{align}
Recall that the linear estimates in Theorem~\ref{thm_linear_estimates} yields a lower bound for the term 
$\| \vec{u}_{R,L}  (t) \|_{\mathcal{H} (r \geq R+ |t|)}$, namely we have
\begin{align*}
\| \pi_R^{\perp} \vec{u}_R (0) \|^2_{\mathcal{H}} \leq 
\max_{\pm} \lim_{t \to \pm \infty}  
\| \vec{u}_{R,L} \|^2_{\mathcal{H} (r \geq R + |t|)} \period
\end{align*}
We then let $|t| \to \infty$ according to the choice of sign dictated by Theorem~\ref{thm_linear_estimates},
which leads to the vanishing of the left hand side in \eqref{eq7.20}. Therefore  we have
\begin{align*}
\| \pi_R^{\perp} \vec{u}_R (0) \|^2_{\mathcal{H}} \leq  
\frac{C^2_0}{R^{ 2 \beta}} \| \vec{u} (0) \|^{2p}_{\mathcal{H} (r \geq R)}.
\end{align*}
Once again using \eqref{eq7.21} we note that $\| \pi_R^{\perp} \vec{u}_R (0) \|^2_{\mathcal{H}}  = \| \pi_R^{\perp} \vec{u} (0) \|^2_{\mathcal{H} (r \geq R)} $,
which gives us 
\begin{align*}
& \| \pi_R^{\perp} \vec{u}  (0) \|^2_{\mathcal{H} (r \geq R)} \leq  
\frac{C^2_0}{R^{ 2 \beta }} \| \vec{u} (0) \|^{2p}_{\mathcal{H} (r \geq R)} \\
& \indeq \indeq \indeq 
\leq \frac{C^2_0}{R^{ 2 \beta }}  (\| \pi_R \vec{u}  (0) \|_{\mathcal{H} (r \geq R)}  + \| \pi_R^{\perp} \vec{u}  (0) \|_{\mathcal{H} (r \geq R)} )^{2p}.
\end{align*}
Then, we choose $R_0$ large enough to absorb $C_0^2 R^{ - 2 \beta} \| \pi_R \vec{u}  (0) \|_{\mathcal{H} (r \geq R)}^{2p} $ on the left side, which completes the proof of Lemma~\ref{Lem4} for $t=0$. 

We utilize Corollary~\ref{cor8.1}  to prove the  inequality \eqref{eq7.3} for all $t \in \mathbb{R}$. 
By the pre-compactness of $K_1$ we may select $R_0 = R_0 (\delta_0)$ such that 
\begin{equation*}
\| \vec{u} (t) \|_{\mathcal{H}(r \geq R)} \leq \min (\delta_0, 1)
\end{equation*}
uniformly in $t \in \mathbb{R}$. 
For fixed $t_0 \in \mathbb{R}$, we take 
\begin{align*}
u_{t_0,R} (r) := 
\begin{cases}
u (t_0, r) & \mbox{for}~ r\geq R \\
u (t_0, R) & \mbox{for}~r \leq R
\end{cases}
\end{align*}
and
\begin{align*}
\tilde{u}_{ t_0,R} (r) := 
\begin{cases}
u_t (t_0, r) & \mbox{for}~ r\geq R \\
0 & \mbox{for}~r \leq R.
\end{cases}
\end{align*}
as the truncated inital data and repeat the same steps to obtain \eqref{eq7.3} for $t= t_0$. 
\end{proof}

\subsection{Step~2}
Next, we aim to investigate the asymptotic behaviour of $(u_0 (r), u_1 (r))$ as $r \to \infty$. 
Our goal is to establish the asymptotic rates given in the following proposition. 
\begin{Proposition}
\label{Proposition2}
Let $\vec{u} (t)$ be as in Proposition~\ref{Proposition1} with $\vec{u} (0) = (u_0, u_1)$. 
Then, there exists $ \ell \in \mathbb{R}$ so that
\begin{align}
& r^{d-2} u_0 (r)  = \ell + O \left( r^{-(d-2) p +d} \right)   \inon{ as }~ r \to \infty    \label{eq8.1} \\
& \int_r^{\infty} u_1(s) s^{2i-1}~ds  = O \left( r^{- (d-2)p +2i +1} \right) \inon{ as }~ r \to \infty \label{eq8.2} 
\end{align}
for each $1\leq i \leq k$. 
\end{Proposition}

First, we recall the bounds on the norms $\Vert \pi_R \vec{u} (t) \Vert_{\mathcal{H}(r \geq R)}$ 
and $\Vert  \pi^{\perp}_R \vec{u} (t) \Vert_{\mathcal{H}(r \geq R)}$ given in Lemma~\ref{Lem5} and 
rewrite Lemma~\ref{Lem4} in terms of the projection coefficients $\lambda_i (t, r)$  and $\mu_i (t,r)$. 

\begin{Lemma}
\label{Lem6}
There exists $R_0 > 0$ so that for all $R > R_0$ we have 
\begin{align}
\begin{split}
& 
\int_{R}^{\infty}
\sum_{i=1}^{\tilde{k}}  \left( \partial_r \lambda_i (t,r) r^{2i - \frac{(d+1)}{2}} \right)^2 
+   \sum_{i=1}^k \left( \partial_r \mu_i (t,r) r^{2i - \frac{(d-1)}{2}}  \right)^2 
~dr
 \\ &
\indeq \indeq \indeq  \les 
 \sum_{i=1}^{\tilde{k}} \lambda_i^{2p} (t,R) R^{  (4i -2d) p+ d+2} + \sum_{i=1}^{k} \mu_i^{2p} (t,R) R^{ (4i - 2 d+2 ) p  +d+2   }  \period
\label{eq8.3}
\end{split}
\end{align}
where the implicit constant on the right hand side is uniform in $t \in \mathbb{R}$. 
\end{Lemma}

\begin{Remark}
{\rm
Lemma~\ref{Lem6} yields uniform in time estimates on the projection coefficients, which then leads to difference estimates. 
Let $\delta_0$ and $R_0$ denote the constants introduced in the proof of Lemma~\ref{Lem4}. We take $\delta_1 \in (0, \delta_0)$ to be determined
below. By the pre-compactness of the set $K_1$ in \eqref{eq7.1}, we may find $R_1  > R_0$ such that 
for all $R \geq R_1$, 
\begin{align}
\Vert  \vec{u}(t) \Vert^{p-1}_{\mathcal{H} (r \geq R)}  \leq \delta_1,  \inon{  }~ t\in   \mathbb{R} 
\label{eq8.4}
\end{align}
and 
\begin{align}
1/R_1 \leq \min( \delta_1, 1) \period
\label{eq8.5}
\end{align}
Consequently, we obtain the following estimates that hold uniformly in time: for every $r \geq R_1$ and for all $t \in \mathbb{R}$
\begin{align}
\frac{|\lambda_i (t,r)|^{p-1}}{r^{\left( -2i +\frac{d+2}{2} \right) (p-1) }} \leq \delta_1\comma
\frac{|\mu_i (t,r)|^{p-1}}{r^{\left(  -2i +\frac{d}{2} \right) (p-1)}} \leq \delta_1 \period
\label{eq8.6}
\end{align}}
\end{Remark}

\begin{Lemma}
\label{Lem7}
Let $R_1$ be as in \eqref{eq8.5}. For all $r, r'$ such that $R_1 \leq r \leq r' \leq 2r$, the following difference estimates 
hold uniformly in time.  We have for each $1 \leq j \leq \tilde{k}$
\begin{align}
\left| \lambda_j (t,r) - \lambda_j (t,r') \right| \les r^{d+2 -pd -2j} \left(  \sum_{i=1}^{\tilde{k}} \left| \lambda_i (t,r) \right|^p r^{2ip}
+    \sum_{i=1}^k \left| \mu_i (t,r) \right|^p r^{2ip +p} \right)
\label{eq8.7}
\end{align}
and for each $1 \leq j \leq k$,
\begin{align}
\left| \mu_j (t,r) - \mu_j (t,r') \right|  \les r^{d+1 -pd -2j} \left(  \sum_{i=1}^{\tilde{k}} \left|\lambda_i (t,r) \right|^p  r^{2 i p} 
+  \sum_{i=1}^k \left| \mu_i (t,r) \right|^p r^{2ip+p} \right) 
\label{eq8.9}
\end{align}
for all $t \in \mathbb{R}$ uniformly. 
\end{Lemma}

\begin{proof}
The inequalities \eqref{eq8.7}--\eqref{eq8.9} follow directly from Lemma~\ref{Lem6}. 
First, we consider \eqref{eq8.7}. We express difference on the left hand side as an integral from $r$ to $r'$ and  apply the inequality \eqref{eq8.3} in Lemma~\ref{Lem6}.
More precisely, we get
\begin{align*}
\begin{split}
  \left|  \lambda_j (t,r) - \lambda_j (t,r') \right|^2 & = \left| \int_r^{r'} \partial_s \lambda_j (t,s)~ds \right|^2  
% \\ & \indeq \indeq \indeq
 \les  \left( \int_r^{r'} \left( \partial_s \lambda_j (t,s) s^{2j - \frac{d+1}{2}} \right)^2~ds  \right)
  \left( \int_r^{r'} s^{d+1 - 4j}~ds  \right) \\
& \les   
 r^{d+2 - 4j}  
 \sum_{i=1}^{\tilde{k}} \lambda_i^{2p} (t,r) r^{  (4i -2d) p+ d+2} + 
  r^{d+2 - 4j}  
  \sum_{i=1}^{k} \mu_i^{2p} (t,r) r^{ (4i - 2 d+2 ) p  +d+2  }
 \period 
\end{split}              
\end{align*}
In the same fashion, we obtain
\begin{align*}
\begin{split}
  \left|  \mu_j (t,r) - \mu_j (t,r') \right|^2 & = \left| \int_r^{r'} \partial_s \mu_j (t,s)~ds \right|^2  
% \\ & \indeq \indeq \indeq
 \les     \left( \int_r^{r'} \left( \partial_s \mu_j (t,s) s^{2j - \frac{d-1}{2}} \right)^2~ds  \right)   \left( \int_r^{r'} s^{-4j + d-1}~ds  \right) \\
& \les    r^{ d-4j} \sum_{i=1}^{\tilde{k}} \lambda_i^{2p} (t,r) r^{  (4i -2d) p+ d+2} + r^{ d-4j} \sum_{i=1}^{k} \mu_i^{2p} (t,r) r^{ (4i - 2 d+2 ) p  +d+2  }
\period
\end{split}              
\end{align*}
\end{proof}

Recalling the setting in \eqref{eq8.4}--\eqref{eq8.6}, we state a direct consequence of the difference estimates above. 

\begin{Corollary}
\label{Cor1}
Let $R_1$ and $\delta_1$ be defined as in \eqref{eq8.5}. Then, for all $r$ and $r'$ 
with $R_1 < r < r'< 2r$ and for all $t \in \mathbb{R}$,  we have
\begin{align}
\begin{split}
 \left|  \lambda_j (t,r) - \lambda_j (t,r') \right| 
& \les  \delta_1  \left(  \sum_{i=1}^{\tilde{k}} |\lambda_i (t,r)| r^{2i - 2j } r^{\frac{d+2}{2} - \frac{pd}{2} +p}    + \sum_{i=1}^k |\mu_i (t,r)| r^{2i - 2j} 
r^{\frac{d}{2} - \frac{pd}{2} +p +2}  \right)  \\
\left| \mu_j (t,r) - \mu_j (t,r') \right| 
 & \les \frac{ \delta_1}{r}   \left(  \sum_{i=1}^{\tilde{k}} |\lambda_i (t,r)| r^{2i - 2j } r^{\frac{d+2}{2} - \frac{pd}{2} +p}    + \sum_{i=1}^k |\mu_i (t,r)| r^{2i - 2j} 
r^{\frac{d}{2} - \frac{pd}{2} +p +2}  \right) \period
\end{split}
\end{align} 
\end{Corollary}

Next, we state a result from \cite{KLLS}  which contains a list of formulas that relate $\vec{u} (t,r)$ to the projection coefficients $\lambda_j (t,r)$ and $\mu_j (t,r)$ directly in the region
\begin{equation}
\Omega_R = \{r \geq R+ |t| \}.
\label{EQ8.1p}
\end{equation}
These formulas will be useful in the proof of Proposition~\ref{Proposition2}. 
\begin{Lemma}
\label{Lem8.10}
For each fixed $R>0$ and for all $(t,r) \in \Omega_R$ we have
\begin{align}
& u(t,r) = \sum_{j=1}^{\tilde{k}} \lambda_j (t,r) r^{2j -d}   \label{EQ8.2p} \\
& \int_r^{\infty} u_t (t,s) s^{2i -1}~ds = \sum_{j=1}^k \mu_{j} (t,r) \frac{r^{2i+2j -d}}{d- 2i -2j} \comma 1 \leq i \leq k  \label{EQ8.3p} \\
& \mu_j (t,r) = \sum_{i=1}^k \frac{r^{d-2i -2j }}{d-2i -2j} c_i c_j \int_r^{\infty} u_t (t,s) s^{2i-1}~ds \comma 1 \leq j \leq k  \label{EQ8.4p}\\
& \lambda_j (t,r) = \frac{d_j}{d- 2j} 
\left( u(t,r) r^{d-2j} + \sum_{i=1}^{\tilde{k}-1}  \frac{2i ~ d_{i+1} ~r^{d-2i -2j}}{d- 2i -2j} \int_r^{\infty} u (t,s) s^{2i -1} ~ds	  \right) \label{EQ8.5p}
\end{align}
where the last line holds for all $1 \leq j \leq \tilde{k}$. Moreover, for any $1 \leq j, j' \leq \tilde{k}$, we have
\begin{align}
\begin{split}
|\lambda_j (t_1, r)   - \lambda_j (t_2, r) |  & \les r^{2j' - 2j} | \lambda_{j'} (t_1, r) - \lambda_{j'} (t_2, r)  |   \\
& \indeq + \sum_{m=1}^k \left|   r^{2m - 2j}   \int_{t_2}^{t_1} \mu_m (t,r)~dt   \right|  .
\end{split}
\label{EQ8.6p}
\end{align}
\end{Lemma}
The proofs of \eqref{EQ8.2p}--\eqref{EQ8.6p} depend on the algebraic identities provided in Section~\ref{sec: orthogonal proj}, and may be found in \cite[Lemma~5.10]{KLLS}.  
As a result of \eqref{EQ8.4p}, 
we obtain a version of Lemma~$5.12$ from \cite{KLLS}.   
 \begin{Lemma}
\label{Lem15.p}
For each $1 \leq j \leq k$, $t_1 \neq t_2$, and $R$ large enough,  we have
\begin{align}
\frac{1}{R} \int_R^{2R} \mu_j (t_1, r) - \mu_j (t,r)~dr
= \sum_{i=1}^k \frac{c_i c_j}{d-2 - 2j} \int_{t_1}^{t_2} ( I (i,j) + II (i,j)  )~dt
\label{EQ8.7p}
\end{align}
with $I(i,j)$ and $II (i,j)$ are given as follows
\begin{align}
\begin{split}
I (i,j) & = \frac{1}{R} \int_{R}^{2R} r^{d-2i - 2j} \int_r^{\infty} \left( u_{ss} (t,s) + \frac{d-1}{s} u_s (t,s) \right) s^{2i -1} ds ~dr \\
& = - \frac{1}{R} \evaluat[\big]{  (u (t,r) r^{d- 2j -1})   }_{r= R}^{r= 2R}
+ (2i - 2j -1) \frac{1}{R} \int_R^{2R} u (t,r) r^{d- 2j -2}~dr \\
& \indeq - \frac{(d- 2i)(2i-2)}{R} \int_R^{2R} r^{d- 2i - 2j} \int_r^{\infty} u (t,s) s^{2i -3} ds~dr
\end{split}
\label{EQ8.8p}
\end{align}
and 
\begin{align}
II (i,j) = \frac{1}{R} \int_{R}^{2R} r^{d-2i - 2j} \int_r^{\infty} \left( |u|^{p-1} u (t,s) \right) s^{2i-1} ds~dr . 
\label{EQ8.9p}
\end{align}
\end{Lemma}
The proof of Lemma~\ref{Lem15.p} follows from the same 
arguments employed in that of \cite[Lemma~5.12]{KLLS}. The only change 
arises in $II (i,j)$ and is due to the difference in nonlinearity on the respective equations.   
Here, $\vec{u} (t,r)$ is a solution of \eqref{wave equations} whereas \cite[Lemma~5.12]{KLLS} is based on solutions 
of wave maps. 

We now concentrate on the proof of Proposition~\ref{Proposition2}. 
We aim towards understanding the asymptotic behaviour of $(u_0, u_1)$ through the projection coefficients $\lambda_j$, $\mu_j$. To that end, we will utilize the algebraic identities and difference estimates collected above to estimate the asymptotic behaviour of these projection coefficients. 
For the rest of this section, we will apply the program demonstrated in \cite[Section~5.2]{KLLS}. 

\subsubsection{Proof of Proposition~\ref{Proposition2} when $d \equiv 3  \pmod 4$}
By \eqref{EQ5.1p}, we have  
\begin{align}
\begin{split}
& d= 4 \tilde{k} -1= 4k +3 \\
& \tilde{k} = k+1. 
\end{split}
\end{align}
The argument is based on $k$ iterative steps. For each 
$2 \leq i \leq \tilde{k}$, we show that both  $\lambda_i (t,r), \mu_{i-1} (t,r)$ converge to $0$ as $r \to \infty$. Finally, for the leading coefficient $\lambda_1 (t,r)$ we set $t=0$, and prove that 
\begin{equation*}
|\lambda_1 (0, r) - \ell | = O(r^{-(d-2) p +d})   \inon{ as }~ r \to \infty
\end{equation*}
where $\ell$ is the limit in Proposition~\ref{Proposition2}.

For notational simplicity, we introduce for each $1 \leq j \leq \tilde{k}$
\begin{align}
\begin{split}
\alpha_j & = d+2 - pd -2j + 2 \tilde{k}p  \\
& = 4 \tilde{k} +1 - 2j - p (2 \tilde{k} -1 ).
\end{split}
\label{EQ8.14p}
\end{align}
Additionally, we denote by
\begin{equation}
\eta_p := (2 \tilde{k} +1) - p (2 \tilde{k} -1) \label{EQ8.eta.2}
\end{equation}

\begin{Remark}\label{R_eta}
{\em    Recalling $p \geq 3$, we deduce that $\eta_p$ is a negative quantity that is bounded from above by 
\begin{align*}
\begin{split}
\eta_p & = (2 \tilde{k} +1) - 3 (2 \tilde{k} -1)  - (p-3) (2 \tilde{k} -1) \\
& \leq -4 .
\end{split}
\end{align*}
The right hand side is achieved in the case $\tilde{k}=2$ and $p=3$, i.e., when $d=7$ and $p=3$. For higher dimensions, we refer back to the definition of $\eta_p$. 
Additionally, this implies that 
\begin{align}
\begin{split}
\alpha_j &= 2 \tilde{k} + \eta_p -2j \\
& = 4 - 2 \tilde{k} - (p-3) (2 \tilde{k} -1) - 2 j \\
& < 0
\end{split}
\label{EQ8.43} 
\end{align}
for every $1 \leq j \leq \tilde{k}$. 
 }
\end{Remark}

\begin{Lemma}
\label{Lem8}
The projection coefficients defined in \eqref{EQ5.2p}--\eqref{EQ5.3p} satisfy the following estimates uniformly in time: 
\begin{align}
\left| \lambda_{j} (t,r) \right| & \les 1  \comma 1 \leq j \leq \tilde{k}    \label{EQ8.12p} \\
\left| \mu_j (t,r) \right| & \les 1  \comma 1 \leq j \leq  k .
\label{EQ8.13p}
\end{align}
\end{Lemma} 

\begin{proof}
Let $\epsilon >0$ be a small number.    We fix $r_0 > R_1$ and set 
$r = 2^n r_0$, $ r' = 2^{n+1} r_0 $
in Corollary~\ref{Cor1}. We then obtain for all $t \in \mathbb{R}$, 
\begin{align*}
  |\lambda_j (t, 2^{n+1} r_0)| &  \leq  |\lambda_j (t, 2^n r_0)| 
  					+ C \delta_1 \left(   \sum_{i=1}^{\tilde{k}} |\lambda_i (t, 2^n r_0)|  (2^n r_0 )^{2i -2j} (2^n r_0)^{\frac{d+2}{2} - \frac{pd}{2} +p}   \right) \\
& \indeq + C \delta_1 \left(    \sum_{i=1}^{{k}} |\mu_i (t, 2^n r_0)|  (2^n r_0 )^{2i -2j} (2^n r_0)^{\frac{d}{2} - \frac{pd}{2} +p+2}   \right)
\end{align*}
and 
\begin{align*}
  |\mu_j (t, 2^{n+1} r_0)| &  \leq       |\mu_j (t, 2^n r_0)| 
  					 + C \delta_1 \left(   \sum_{i=1}^{\tilde{k}} |\lambda_i (t, 2^n r_0)|  (2^n r_0 )^{2i -2j} (2^n r_0)^{\frac{d}{2} - \frac{pd}{2} +p}   \right) \\
& \indeq + C \delta_1 \left(    \sum_{i=1}^{{k}} |\mu_i (t, 2^n r_0)|  (2^n r_0 )^{2i -2j} (2^n r_0)^{\frac{d}{2} - \frac{pd}{2} +p+1}   \right)
\period
\end{align*}
Denote by
\begin{equation*}
H_n := 
\sum_{j=1}^{\tilde{k}}  |\lambda_j (t, 2^n r_0)|  (2^n r_0 )^{2j -2 \tilde{k}}  +   \sum_{j=1}^{{k}} |\mu_j (t, 2^n r_0)|  (2^n r_0 )^{2j +1 - 2 \tilde{k}} .
\end{equation*}
We deduce that
\begin{align}
\begin{split}
H_{n+1} & \leq  \left(1+  C(k+ \tilde{k}) \delta_1  \right) H_n \\
& \leq \left(1+  C(k+ \tilde{k}) \delta_1  \right)^{n+1} H_0 
\end{split}
\label{eq8.14}
\end{align}
since 
\begin{align*}
\frac{d}{2} - \frac{pd}{2} +p +1   = \frac{d+2}{2} - p \left( \frac{d-2}{2}\right)  \leq 0 \period
\end{align*}
We then pick $\delta_1 >0$ sufficiently small that 
\begin{align}
 1 + C (k + \tilde{k})  \delta_1 < 2^{\epsilon}
 \label{EQ8.16p}
\end{align}
and the inequality \eqref{eq8.14} above yields that
\begin{align*}
H_n \leq C \left( 2^n r_0  \right)^{\epsilon} \period
\end{align*}
Note that  the constant $C$ depends only on $H_0$, which is uniformly bounded for all $t \in \mathbb{R}$ by 
\eqref{eq8.6} once $r_0 > R_1$ is fixed.  
It then follows that 
\begin{align}
\begin{split}
|\lambda_i (t, 2^n r_0) | & \leq (2^n r_0)^{2 \tilde{k} - 2i + \epsilon} \comma 1 \leq i \leq \tilde{k} \\
|\mu_i (t, 2^n r_0)|  & \leq (2^n r_0)^{2 \tilde{k} - 2i -1 + \epsilon} \comma 1 \leq i \leq k .
\end{split}
\label{eq8.15p}
\end{align}
Next we use these estimates in difference inequalities \eqref{eq8.7}--\eqref{eq8.9}.
For $1 \leq j \leq \tilde{k}$, \eqref{eq8.7} yields
\begin{align*}
\begin{split}
| \lambda_j (t, 2^{n+1} r_0)  - \lambda_j (t, 2^n r_0) | & \les   
 (2^n r_0)^{d+2 - pd - 2j} \sum_{i=1}^{\tilde{k}} |\lambda_i (t, 2^n r_0)|^p (2^n r_0)^{2ip} \\
& \indeq + (2^n r_0)^{d+2 - pd - 2j} \sum_{i=1}^{{k}} |\mu_i (t, 2^n r_0)|^p (2^n r_0)^{2ip +p} 
\period
\end{split}
\end{align*}
By \eqref{eq8.15p}, we then have
\begin{align}
|\lambda_j (t, 2^{n+1} r_0)| \leq |\lambda_j (t, 2^n r_0)| + C (2^n r_0 )^{\alpha_j + p \epsilon}.
\label{EQ8.15p}
\end{align}
Iterating \eqref{EQ8.15p}, we obtain
\begin{align}
|\lambda_j (t, 2^{n+1} r_0)| \leq  |\lambda_j (t, r_0)| + C \sum_{m=1}^{n} (2^m r_0)^{\alpha_j + p \epsilon} .  
\label{EQ8.17p}
\end{align}
As discussed in Remark~\ref{R_eta}, we have $\alpha_j <0$ for each $1 \leq j \leq \tilde{k}$. Also, recalling the uniform in time bounds on the projection coefficients in\eqref{eq8.6}, we deduce that, for sufficiently small $\epsilon >0$,  \eqref{EQ8.17p} implies 
\begin{align}
|\lambda_j (t, 2^n r_0)| & \les 1  \comma 1 \leq j \leq \tilde{k}.
\label{EQ8.18p}
\end{align}

Similarly, plugging \eqref{eq8.15p} into \eqref{eq8.9} leads to 
\begin{align}
|\mu_j (t, 2^{n+1} r_0)| \leq  |\mu_j (t, 2^n r_0)| 
+ C  (2^n r_0 )^{\alpha_j -1 +p \epsilon}
\label{EQ8.19p}
\end{align}
for each $1 \leq j \leq k$. 
We proceed as above and deduce that
 \begin{align}
|\mu_j (t, 2^n r_0)| & \les 1  \comma 1 \leq j \leq k .
\label{EQ8.20p}
 \end{align}
Finally,  combining these growth estimates with the difference estimates in Lemma~\ref{Lem7} one more time for any $2^n r_0 < r < 2^{n+1} r_0$,  
we obtain the result for arbitrary $r> R_1$. 
\end{proof}

\begin{Proposition}
\label{Proposition4}
The following estimates for the projection coefficients $\lambda_j (t,r)$ and $\mu_j (t,r)$ 
hold true for $S= 0, \ldots, k$,  uniformly in time:
  We have
  \begin{align}
  \begin{split}
  \begin{array}{llr}
   |\lambda_j (t,r) | & \les r^{\alpha_j  - 2p S} \comma  & \tilde{k} - S < j \leq \tilde{k} \\
  |\lambda_j (t,r)| & \les 1      \comma & 1 \leq j  \leq \tilde{k} - S .
  \end{array}
  \end{split}
  \label{EQ8.21p}
  \end{align}
Similarly, we have
   \begin{align}
  \begin{split}
  \begin{array}{llr}
  |\mu_j (t,r) | & \les r^{\alpha_j -1  - 2p S} \comma    & {k} - S < j \leq {k} \\
  |\mu_j (t,r)| &\les 1  \comma & 1 \leq j \leq  {k} - S .
\end{array}  
  \end{split}
  \label{EQ8.22p}
  \end{align}
\end{Proposition}

\begin{proof}[Proof of Proposition~\ref{Proposition4}]
We established the case $S=0$ by Lemma~\ref{Lem8}. The rest of the argument will be justified by 
induction, which is built on a series of lemmas. More precisely, throughout Lemmas~\ref{Lem9}--\ref{Lem13} 
we assume that the estimates in  \eqref{EQ8.21p}--\eqref{EQ8.22p} 
 are true for $S$ with $0 \leq S \leq k-1$, and prove that they must also hold for $S+1$. 
\begin{Lemma}
\label{Lem9}
There exist bounded functions $\theta_{\tilde{k} - S} (t)$ and $\rho_{k- S} (t)$ such that 
\begin{align}
|\lambda_{\tilde{k} - S} (t,r) - \theta_{\tilde{k} - S} (t) | & \les O \left( r^{\alpha_{(\tilde{k} - S)}  - 2 p S  } \right)  \label{EQ8.23p} \\
|\mu_{k -S} (t,r)   - \rho_{k-S} (t)|  & \les O \left( r^{\alpha_{(k-S)} -1 - 2p S }\right) \label{EQ8.24p}
\end{align}
uniformly in time. 
\end{Lemma}
 
 \begin{proof}
 We plug our induction hypothesis into the difference estimates \eqref{eq8.7}--\eqref{eq8.9}. Letting
 $R_1 < r < r' < 2r$, with $R_1$  as in \eqref{eq8.5}, we have
 \begin{equation}
 |\lambda_j (t,r) - \lambda_j (t, r')| \les r^{d+2 - pd - 2j}   (I_{\lambda} (r)  + I_{\mu} (r))
\label{EQ8.25p}
 \end{equation}
 where $I_{\lambda}$ and $I_{\mu}$ represent the sums corresponding to $\lambda_i (t,r)$ and $\mu_i (t,r)$ in \eqref{eq8.7}, respectively.  Splitting both $I_{\lambda}$ and $I_{\mu}$ into two parts, we get
 \begin{align*}
 I_{\lambda} (r)  \les 
 r^{2p (\tilde{k} - S)} 
  +  \sum_{i> \tilde{k} - S}  r^{p \alpha_i  - 2 p^2 S}  r^{2ip}  
 \end{align*}
 and 
  \begin{align*}
 I_{\mu} (r)  \les 
  r^{p (2 {k} -2 S +1)} 
   + \sum_{i> {k} - S}  r^{p \alpha_i - p  - 2 p^2 S}  r^{p (2 i +1 )}.
 \end{align*}
 By \eqref{EQ8.14p}, we get
 \begin{equation*}
 p \alpha_i + 2 i p = (4 \tilde{k} +1) p - p^2 (2 \tilde{k} +1) 
 \end{equation*}
 which implies that it suffices to control the $I_{\lambda} (r)$ term in \eqref{EQ8.25p}. 
 To that end, 
 we  estimate the growth rate in  \eqref{EQ8.25p} by simply comparing the exponents in $ I_{ \lambda} (2^n r_0)$. 
 By Remark~\ref{R_eta}, we have $\eta_p <0$. More precisely, plugging in the definition of $\eta_p$ we obtain
 \begin{align}
 \begin{split}
 (4 \tilde{k} +1)p - p^2 (2 \tilde{k} +1) -2p^2 S  & =
 2 \tilde{k} p + p \eta_p - 2 p^2 (1 +S ) \\
 & < 0 \\
 & <  2p (\tilde{k} - S). 
 \end{split}
 \label{EQ8.26.2p}
 \end{align}
 Therefore, fixing $r_0 > R_1$ we may estimate
 \begin{align}
 |\lambda_j (t, 2^{n+1} r_0) - \lambda_j (t, 2^n r_0)| \les (2^n r_0)^{d+2 - pd - 2j}  (2^n r_0)^{2p (\tilde{k} - S) }.
 \label{EQ8.27p}
 \end{align}
 Setting $j = \tilde{k} - S$ in \eqref{EQ8.27p}, using \eqref{EQ8.eta.2} and $d = 4 \tilde{k} -1$,  we find
 \begin{align*}
 |\lambda_{\tilde{k} - S} (t, 2^{n+1} r_0) - \lambda_{\tilde{k}- S} (t, 2^n r_0)| \les 
 (2^n r_0)^{\eta_p  + 2 S - 2p S }
 \end{align*}
 which comes with a negative power leading to
 \begin{equation*}
 \sum_n  |\lambda_{\tilde{k} - S} (t, 2^{n+1} r_0) - \lambda_{\tilde{k}- S} (t, 2^n r_0)| < \infty .
 \end{equation*}
 As a result, there exists a function $\theta_{\tilde{k} - S} (t)$ so that
 \begin{align*}
 \lim_{n \to \infty} \lambda_{\tilde{k} - S } (t, 2^n r_0) = \theta_{\tilde{k} - S} (t).
 \end{align*}
 Moreover, we have 
 \begin{align}
 \begin{split}
 |\theta_{\tilde{k}- S} (t) - \lambda_{\tilde{k} - S} (t, r_0)| & = \lim_{n \to \infty}   |\lambda_{\tilde{k} - S} (t, 2^{n+1} r_0) - \lambda_{\tilde{k}- S} (t, r_0)| \\
 & \les   \lim_{n \to \infty} \sum_{m=1}^n  |\lambda_{\tilde{k} - S} (t, 2^{n+1} r_0) - \lambda_{\tilde{k}- S} (t, 2^n r_0)| \\
 & \les r_0^{\alpha_{\tilde{k} - S}} r_0^{ - 2 p S} \sum_{m=1}^n  (2^m)^{\alpha_{\tilde{k} - S}}  (2^m)^{ -2 p S} .
 \end{split}
 \label{EQ8.28p}
 \end{align}
Since $\lambda_{\tilde{k} -S} (t, r_0)$ is uniformly bounded in $t$, we conclude that $\theta_{\tilde{k} - S} (t)$ is also uniformly bounded. Using the difference estimate once again we get
\begin{align*}
\lim_{r \to \infty} \lambda_{\tilde{k} - S} (t,r) = \theta_{\tilde{k} - S} (t).
\end{align*}

 In the same fashion, we may show that
 \begin{align}
 |\mu_{j} (t, 2^{n+1} r_0) - \mu_j (t, 2^n r_0)|  \les (2^n r_0)^{d+1 - pd -2j} (2^n r_0)^{2p (\tilde{k} - S) } 
 \label{EQ8.29p}
 \end{align}
 Setting $j = k-S$ in \eqref{EQ8.29p}, we obtain
 \begin{align*}
 |\mu_{k- S} (t, 2^{n+1} r_0) - \mu_{k -S} (t, 2^n r_0)|  \les (2^n r_0)^{\alpha_{k-S}} (2^n r_0)^{-1 - 2pS}. 
 \end{align*}
 Since the power  $\alpha_{k-S} -1 - 2pS$ is also negative, 
 arguing as above, we deduce that there exists a function $\rho_{k-S} (t)$ such that
 \begin{align*}
 \lim_{r \to \infty} \mu_{k- S} = \rho_{k - S} (t)
 \end{align*}
 and both $\mu_{k-S} (t,r)$ and $\rho_{k-S} (t)$ are in fact uniformly bounded in time. 
 
 Using the uniform boundedness of $\lambda_{\tilde{k} -S} (t,r)$ and $\mu_{k-S} (t,r)$, we refine the difference estimates 
 \eqref{EQ8.25p} and \eqref{EQ8.29p}. Recalling the asymptotic rates in  $I_{ \lambda} (t, 2^n r_0)$, we obtain
 \begin{align}
 \begin{split}
&  |\lambda_{\tilde{k} - S} (t, 2^{n+1} r_0) - \lambda_{\tilde{k}- S} (t, 2^n r_0)| \\
& \indeq  \les 
 (2^n r_0)^{d+2 - pd - 2 (\tilde{k} - S)}  
 \left(      (2^n r_0)^{(4 \tilde{k} +1) p - p^2 (2 \tilde{k} +1) - 2 p^2 S } + (2^n r_0)^{2p (\tilde{k} - S)}    \right) \\
 & \indeq \les (2^n r_0)^{\alpha_{\tilde{k} - S}} (2^n r_0)^{- 2p S }
\end{split} 
\label{EQ8.30p}
 \end{align}
 where we used the  comparison \eqref{EQ8.26.2p} at the last step. 
 
 Repeating the same arguments to improve the difference estimate for $\mu_{k - S} (t,r)$, we have
 \begin{equation}
 |\mu_{k-S} (t, 2^{n+1} r_0) - \mu_{k- S} (t, 2^n r_0)| \les (2^n r_0)^{\alpha_{k-S}} (2^n r_0)^{-1 - 2p S }. 
 \label{EQ8.31p}
 \end{equation}
 From  \eqref{EQ8.30p} and \eqref{EQ8.31p}, we then obtain the desired asymptotic estimates. Namely, for any $r > R_1$
 \begin{align*}
  |\lambda_{\tilde{k} - S} (t, r) - \theta_{\tilde{k} - S} (t) | & \les \sum_{n=0}^{\infty} (2^n r)^{\alpha_{\tilde{k} - S}} (2^n r)^{- 2p S }
  \les r^{\alpha_{\tilde{k}- S} ~-2pS }  \\
   |\mu_{k-S} (t, r) - \rho_{k- S} (t)| & \les \sum_{n=0}^{\infty} (2^n r)^{\alpha_{k-S}} (2^n r )^{-1 - 2p S } \les
   r^{\alpha_{k-S} ~-1 -2pS }. 
 \end{align*}
 \end{proof}
 
 Combining the asymptotic estimate  for $\lambda_{\tilde{k} - S} (t,r)$  with the expansion  formula for $u$ as given in  Lemma~\ref{Lem8.10}, we obtain the following result.
\begin{Lemma}
\label{Lem12}
The following holds uniformly in time: 
\begin{align}
r^{- 2 (\tilde{k} - S)  +d} u (t,r) =  \theta_{\tilde{k} - S} (t) + O\left( r^{  -2} \right).
\end{align}
\end{Lemma}

\begin{proof}
 We expand $u(t,r)$ by \eqref{EQ8.2p}, and split it into three parts
\begin{align}
\begin{split}
u (t,r) & = \sum_{j=1}^{\tilde{k}} \lambda_j (t,r) r^{2j -d}  \\
& = \sum_{j > \tilde{k} - S}   I_{1,j} (r) +  \lambda_{\tilde{k} - S} (t,r) r^{2(\tilde{k} - S) -d} 
+ \sum_{1 \leq j < \tilde{k} - S} I_{3,j} (r). 
\end{split}
\label{EQ8.32p}
\end{align}
Using the induction hypothesis \eqref{EQ8.21p} on $I_{1, j} (r)$ and $I_{3, j}$ we get
\begin{align*}
I_{1,j} &  \les r^{\alpha_j  - 2p S} r^{2j -d} \les r^{2 - p (2 \tilde{k} -1)  - 2 p S} \\
I_{3, j} & \les  r^{2j -d}
\end{align*}
which implies that the first  sum in \eqref{EQ8.32p} will be dominated by the latter sum with $I_{3,j}$. 
Next, we apply \eqref{EQ8.23p} on the middle term in \eqref{EQ8.32p}, which yields
\begin{align*}
 \lambda_{\tilde{k} - S} (t,r)  r^{2(\tilde{k} - S) -d} 
=     \theta_{\tilde{k} - S} (t)  r^{2 (\tilde{k} - S) -d}  +   O \left( r^{\alpha_{\tilde{k} -S} - 2p S } ~ r^{2 (\tilde{k} - S) -d} \right)
\end{align*}
Multiplying both side of \eqref{EQ8.32p} with $r^{- 2 (\tilde{k} - S) +d}$, recalling \eqref{EQ8.14p} and \eqref{EQ8.eta.2}, we then obtain
\begin{align}
\begin{split}
 r^{- 2 (\tilde{k} - S) +d}  u (t,r) &  =  \theta_{\tilde{k} - S} (t) + O \left(  r^{\alpha_{\tilde{k} -S} - 2p S }     \right)  \\
 & \indeq  + \sum_{1 \leq j < \tilde{k} - S } \left(  r^{2j -d}     \right)   r^{- 2 (\tilde{k} - S) +d}  \\
 & = \theta_{\tilde{k} - S} (t) +  O \left(   r^{\eta_p + 2S - 2p S} \right) \\
 & \indeq + O \left(  r^{ -2} \right) .
\end{split}
\label{EQ8.33p}
\end{align}
By Remark~\ref{R_eta}, the asymptotic rate on the right hand side above is dominated by $O \left(  r^{-2} \right) $. 
\end{proof}
 \begin{Lemma}
\label{Lem10}
The limit $\theta_{\tilde{k} - S} (t)$ is independent of time. 
\end{Lemma}

 \begin{proof}
The result follows from  \eqref{EQ8.6p} in Lemma~\ref{Lem8.10}.
We take $t_1$ and $t_2 \neq t_1$ and check the difference of $\theta_{\tilde{k} - S} (t_1)$ and $\theta_{\tilde{k} - S} (t_2)$. 
By \eqref{EQ8.23p} and \eqref{EQ8.6p},
\begin{align}
\begin{split}
|\theta_{\tilde{k} - S} (t_1) - \theta_{\tilde{k} - S} (t_2)| & = |\lambda_{\tilde{k} - S} (t_1,r) - \lambda_{\tilde{k} - S} (t_2, r)| 
+ O \left(  r^{\alpha_{\tilde{k} -S} - 2p S }     \right) \\
& \les  r^{-2}  |\lambda_{{k} - S} (t_1,r) - \lambda_{{k} - S} (t_2, r)| 
+ \sum_{m=1}^{k} \int_{t_2}^{t_1}    r^{2m - 2 (\tilde{k} - S)} |\mu_m (t,r)|~dt \\
&  \indeq + O \left(  r^{\alpha_{\tilde{k} -S} - 2p S }     \right) 
 \end{split}
 \label{EQ8.34p}
\end{align} 
We once again use the induction hypothesis on the asymptotic growth rates of $\lambda_{{k} - S} (t,r)$
and $\mu_{m} (t,r)$, and compare these growth rates. More precisely, by \eqref{EQ8.21p}, we have
\begin{equation*}
 r^{-2}  |\lambda_{{k} - S} (t_1,r) - \lambda_{{k} - S} (t_2, r)|  \les r^{-2 }. 
\end{equation*}
Also, using  Lemma~\ref{Lem9} and \eqref{EQ8.21p}
\begin{align*}
\sum_{m=1}^{k}    r^{2m - 2 (\tilde{k} - S)} |\mu_m (t,r)|~dt 
& \les \sum_{m=1}^{k- S-1}   r^{2m - 2 (\tilde{k} - S)}  \\
& \indeq + r^{-2 } +
\sum_{m= k-S +1}^{k} r^{2m - 2(\tilde{k} - S)} r^{\alpha_m -1 - 2pS } \\
& \les  r^{-2 } .
\end{align*}
The last step above follows from a direct comparison of the asymptotic rates. By \eqref{EQ8.14p} and Remark~\ref{R_eta}, we get
\begin{align*}
\alpha_m + 2m -1 - 2p S - 2 (\tilde{k} - S) & = \eta_p - 2p  -1 - 2pS + 2S \\
& < -4 .
\end{align*}
Also, the term $O \left(  r^{\alpha_{\tilde{k} -S} - 2p S }     \right) $ is dominated by $O (r^{-2 i})$ as shown in \eqref{EQ8.33p}.

Therefore, collecting the estimates on each term on the right side of \eqref{EQ8.34p}, we obtain
\begin{align*}
|\theta_{\tilde{k} - S} (t_1) - \theta_{\tilde{k} - S} (t_2)|  \les  O \left( r^{-2 }  (1 + |t_1 - t_2|) \right) .
\end{align*}
Letting $r \to \infty$, we arrive at the conclusion that $\theta_{\tilde{k} - S} (t_1)  = \theta_{\tilde{k} - S} (t_2)$.
\end{proof}
From here on, we will denote the limit  by $\theta_{\tilde{k} - S}$.

\begin{Lemma}
\label{Lem11}
The limit $\rho_{k - S} (t)$ is independent of time.  Furthermore,  $\theta_{\tilde{k} -S} =0$. 
%Moreover, the limit   
\end{Lemma}

The proof of Lemma~\ref{Lem11} is parallel to that of Lemma~$5.28$ in \cite{KLLS}. We will provide the details that arise from the change in nonlinearity and the difference of asymptotic rates of $u (t,r)$. 
\begin{proof}
We begin with the first claim. Selecting two distinct times $t_1$ and $t_2$, by \eqref{EQ8.24p}, we have
\begin{align*}
\begin{split}
\rho_{k- S } (t_1) - \rho_{k- S} (t_2) & = \frac{1}{R} \int_R^{2R} \left(   \rho_{k- S } (t_1) - \rho_{k- S} (t_2)  \right) ~dr \\
& = \frac{1}{R} \int_R^{2R} \left(   \mu_{k- S } (t_1, r) - \mu_{k- S} (t_2, r)  \right) ~dr 
+ O \left(    R^{\alpha_{k- S}  -1 - 2p S } \right) \\
& = I_{\mu} (R) + O \left(    R^{\alpha_{k- S}  -1 - 2p S } \right).
\end{split}
\end{align*}
The first term above may be expressed by Lemma~\ref{Lem15.p}
\begin{align}
I_{\mu} (R)
= \sum_{i =1}^{k}  \frac{c_i c_j}{d - 2i - 2 (k - S)} \int_{t_1}^{t_2} (I (i, k-S) + II (i, k - S))~dt
\label{EQ8.35p}
\end{align}
where the formulas for $I (i, k-S)$ and $II(i, k-S)$ are as in \eqref{EQ8.8p}--\eqref{EQ8.9p}. 

First, we consider  $I (i, k-S)$ which includes three terms:
\begin{align*}
\begin{split}
I (i, k- S) & = - \frac{1}{R} \evaluat[\big]{  (u (t,r) r^{d- 2(k - S) -1})   }_{r= R}^{r= 2R}
+ (2i - 2(k - S) -1) \frac{1}{R} \int_R^{2R} u (t,r) r^{d- 2(k - S) -2}~dr \\
& \indeq - \frac{(d- 2i)(2i-2)}{R} \int_R^{2R} r^{d- 2i - 2(k - S)} \int_r^{\infty} u (t,s) s^{2i -3} ds~dr .
\end{split}
\end{align*}
Recalling Lemma~\ref{Lem12} 
\begin{align*}
r^{d - 2 (k - S) - 2} u (t,r) = \theta_{\tilde{k} - S} + O \left( r^{-2 } \right) 
\end{align*}
we have
\begin{equation*}
r^{d - 2i - 2(k - S) } \int_r^{\infty} u (t,s) s^{2i -3} ~ds
=  \frac{\theta_{\tilde{k} - S}}{d - 2i - 2 (k - S)} + r^{ - 2 } . 
\end{equation*}
Collecting these estimates, we arrive at the equality
\begin{align*}
I (i, k- S) & =  - \frac{1}{R} \evaluat[\big]{ ( \theta_{\tilde{k} - S} r + O ( r^{-1 })  ) }_{r= R}^{r= 2R} \\ 
& \indeq +  (2i - 2(k - S) -1)(   \theta_{\tilde{k} - S}  + O ( R^{-2 }) ) \\
& - \frac{(d - 2i) (2i - 2)}{d - 2i - 2 (k - S)} \theta_{\tilde{k} - S} + O (R^{- 2 }). 
\end{align*}
We then add up the coefficients of $\theta_{\tilde{k} - S}$ in the formula above, and we find
\begin{equation}
I (i, k- S) =   \frac{-2 (k- S) (d-2 -2(k-S))}{d - 2i - 2 (k - S)}   \theta_{\tilde{k} - S} + O (R^{- 2 }). 
\end{equation}
As a result, the contribution from the term involving $I (i, k-S)$ in \eqref{EQ8.35p} may be controlled by
\begin{align}
\sum_{i =1}^{k}  \frac{c_i c_j}{d - 2i - 2 (k - S)}   \frac{-2 (k- S) (d-2 -2(k-S))}{d - 2i - 2 (k - S)} 
\theta_{\tilde{k} - S} (t_1 - t_2) + O (R^{-2 }) .
\label{EQ8.36p}
\end{align}
We note that  the coefficient in front of $\theta_{\tilde{k} - S} (t_1 - t_2)$ is nonzero  by \cite[Remark~.5.29]{KLLS}.

Next, we go back to \eqref{EQ8.35p} and  estimate the contribution of the term involving $II(i, k- S)$ in $I_{\mu} (R)$. 
We employ the asymptotic estimates for $u (t,r)$ in Lemma~\ref{Lem12} once again, and we find
\begin{align*}
\int_r^{\infty} |u|^{p-1} u (t,s) s^{2i -1} ~ds  = O \left(  r^{p (-d + 2 (k-S) +2)  + 2i}  \right)
\end{align*}
which yields
\begin{align*}
II(i, k - S) & = \frac{1}{R} \int_R^{2R} r^{d- 2i - 2(k-S)} \int_r^{\infty}   |u|^{p-1} u (t,s) s^{2i -1} ~ds \\
& \les O \left(  R^{- (p-1)d + (p-1) (2 (k-S) )  + 2p} \right).
\end{align*}
Note that the right hand side above may be controlled by the asymptotic term in \eqref{EQ8.36p}. 

Finally, collecting the estimates on $I_{\mu} (R)$ and selecting the dominant asymptotic error,  we bound the right side of \eqref{EQ8.35p} by
\begin{align}
| \rho_{k- S } (t_1) - \rho_{k- S} (t_2) | =  C |t_1 - t_2| |\theta_{\tilde{k} - S}| + O \left(  R^{-2 } (1 + |t_1 - t_2|) \right)
\label{EQ8.37p}
\end{align}
where $C$ denotes the nonzero prefactor in front of  $\theta_{\tilde{k} - S}$.
Letting $R \to \infty$ we find that
\begin{align*}
| \rho_{k- S } (t_1) - \rho_{k- S} (t_2) | =  C |t_1 - t_2| |\theta_{\tilde{k} - S}|  . 
\end{align*}
Recalling the boundedness of $\rho_{k-S} (t)$, we deduce that 
\begin{align*}
|\theta_{\tilde{k} - S}| = \frac{1}{C } \frac{| \rho_{k- S } (t_1) - \rho_{k- S} (t_2) |}{|t_1 - t_2|} \to 0  \inon{ as } ~ |t_1 - t_2| \to \infty.
\end{align*}
Having shown $\theta_{\tilde{k} - S} =0$, we refine \eqref{EQ8.37p}, and for fixed $t_1 \neq t_2$ we use 
\begin{equation}
| \rho_{k- S } (t_1) - \rho_{k- S} (t_2) |  = O \left(  R^{-2 }  (1 + |t_1 - t_2|) \right) \to 0 \inon{ as } ~ R \to \infty 
\end{equation}
which implies that $\rho_{k-S} (t) $ must be independent of time. 
\end{proof}

\begin{Lemma}
\label{Lem13}
The limit $\rho_{k - S} = 0$. 
\end{Lemma}

\begin{proof}
By the way of contradiction, we assume $\rho_{k - S} \neq 0$. 
Using the asymptotic estimate  \eqref{EQ8.24p}, 
we may guarantee that $ \mu_{k-S} (t, R)$ and $\rho_{k - S}$ share the same sign as well as 
\begin{align}
\left| \mu_{k-S} (t,R) \right|  > \frac{1}{2} |\rho_{k- S} | 
\label{EQ8.38p}
\end{align}
for  $R >0$ sufficiently large. 
Next, we integrate  \eqref{EQ8.38p} from $0$ to $T$, which yields
\begin{align}
\left|  \int_0^{T}    \mu_{k-S} (t,R)~dt \right|  > \frac{T}{2} |\rho_{k -S} | . 
\label{eq8.24}
\end{align}
At the same time, applying formula \eqref{EQ8.4p}  we may bound the left hand side in \eqref{eq8.24} from above
\begin{align}
\left|  \int_0^{T}    \mu_{k-S} (t,R)~dt \right|  \les 
\sum_{i=1}^k R^{d-2i - 2 (k- S)}  \left|  \int_{R}^{\infty} \int_{0}^T    u_t (t,s) s^{2i -1} ~dt ds    \right| . 
\label{EQ8.39p}
\end{align}
Furthermore,  we utilize the asymptotic estimate in Lemma~\ref{Lem12} with $\theta_{\tilde{k} - S} =0$, which leads to 
\begin{align*}
\left|  \int_{R}^{\infty} \int_{0}^T    u_t (t,s) s^{2i -1} ~dt ds   \right| 
& \les \int_R^{\infty} | u(T, s) - u (0,s) | s^{2i-1} ~ds \\
& \les O \left(  R^{-d + 2 (k- S)  + 2i  } \right) .
\end{align*}
 As a result, feeding this  asymptotic estimate  into \eqref{EQ8.39p}, and fixing a sufficiently large $R$,  we obtain a uniform in time control on the right hand side of 
 \eqref{eq8.24}. We obtain
\begin{align*}
\frac{T}{2} |\rho_{k - S} | \les 1. 
\end{align*}
As  we may take $T \to \infty$, we find that $\rho_{k - S}$ must be  zero. 
\end{proof}

We may now complete the proof of Proposition~\ref{Proposition4}. First, we highlight a couple of arguments 
in the remark below that were used in the proof of Lemma~\ref{Lem9} and will be useful again as we close the inductive step. 
%We may now complete the proof of Proposition~\ref{Proposition4}. First, we recall Remark~$5.31$ from \cite{KLLS} 
%highlighting an argument on the use of difference estimates such as \eqref{eq8.7}--\eqref{eq8.9}. 
%that were used  both  in Section  Lemma~\ref{Lem9} 
%will be useful again as we close the inductive step. 
\begin{Remark}[\text{See also \cite[Remark~5.31]{KLLS}}] 
\label{R5}
{\em    Let $\delta_1$ and $R_1$ be as in \eqref{eq8.4}--\eqref{eq8.5}, and fix $r_0 > R_1$.
Suppose the difference estimates \eqref{eq8.7}--\eqref{eq8.9} yields
\begin{equation}
|\lambda_{j} (t, 2^{n+1} r_0) - \lambda_{j} (t, 2^n r_0) | \les (2^n r_0)^a
\label{EQ8.40p}
\end{equation}
for some $a < 0$. 
Arguing as in the proof of Lemma~\ref{Lem9}, we may guarantee the existence of a limit $\theta_{j} (t)$ with the asymptotic estimates
\begin{equation*}
|\lambda_j (t,r) - \theta_j (t)| \les r^{a}
\end{equation*}
uniformly in time. Moreover, $a <0$ also leads to the uniform boundedness of $|\lambda_j (t,r)|$ and $\theta_j (t)$. 
}
\end{Remark}

Going back to the proof Proposition~\ref{Proposition4}, under the induction hypothesis that  \eqref{EQ8.21p}--\eqref{EQ8.22p}
hold true for $S$ with $0 \leq S \leq k-1$, we improved the asymptotic estimates on $\lambda_{\tilde{k} - S} (t,r)$ and $\mu_{k-S} (t,r)$, namely 
we have
\begin{align}
\begin{split}
& |\lambda_{\tilde{k} - S} (t,r)| = O \left(r^{\alpha_{\tilde{k} - S} - 2 p S } \right)  \inon{ as }~ r \to \infty   \\
& |\mu_{k -S} (t,r)| = O \left( r^{\alpha_{k- S} -1 - 2 p S } \right) \inon{ as }~ r \to \infty 
\end{split}
\label{eq8.25}
\end{align}
by Lemmas~\ref{Lem9}--\ref{Lem13}. 

We now feed the induction hypothesis along with \eqref{eq8.25} into the difference estimate \eqref{eq8.7}. 
For $R_1 <r < r' < 2r$, we have
\begin{align*}
\begin{split}
|\lambda_j (t, r') - \lambda_j (t,r)| &\les r^{d+2 - pd - 2j} 
\left(\sum_{i=1}^{\tilde{k}} |\lambda_{i} (t,r)|^{p} r^{2ip}  
+ \sum_{i=1}^k   |\mu_{i} (t,r)|^{p} r^{2ip}  \right)
\\ 
 &\les r^{d+2 - pd - 2j} 
\left(   \sum_{i \geq \tilde{k} -S } r^{\alpha_{i} p + 2i p  - 2 p^2 S }  + \sum_{i < \tilde{k} - S} r^{ 2 ip }  \right).
\end{split}
\end{align*}
We once again observed that the contribution from $|\mu_{i} (t,r)|^{p} r^{2ip}$  for each $1\leq i\leq k$ is dominated by the first sum on $|\lambda_{i} (t,r)|^{p} r^{2ip}$. 
Recalling the definitions 
of $\alpha_i$ and $\eta_p$  in \eqref{EQ8.14p} and \eqref{EQ8.eta.2}, and using the computations in \eqref{EQ8.26.2p} we deduce that
\begin{align}
\begin{split}
|\lambda_j (t, r') - \lambda_j (t,r)| & \les r^{d+2 - pd - 2j}  \left(   r^{2 \tilde{k} p + p \eta_p - 2 p^2 (S+1 )} 
+ r^{ 2p (k-S)} \right) \\
& \les r^{d+2 - pd - 2j}   r^{ 2 p (k-S)} \\
& \les r^{\alpha_j - 2 p S -2p} .
\end{split}
\label{EQ8.42p}
\end{align}
By Remark~\ref{R_eta}, the power on the asymptotic rate in \eqref{EQ8.42p} is negative.       
In particular, setting $j = \tilde{k} - S -1$ and  following the argument in Remark~\ref{R5} we obtain
\begin{equation*}
|\lambda_{\tilde{k} - S -1} (t, r)| \les O (1) . 
\end{equation*}
 On the other hand, 
for $j > \tilde{k} - S -1$, we already know that $\lambda_j (t,r) \to 0$ as $r \to \infty$ uniformly in time. As a result, \eqref{EQ8.42p} implies that
\begin{equation*}
|\lambda_j (t,r)| \les O \left( r^{\alpha_j - 2 p (S+1) }   \right) \comma j \geq  \tilde{k} - S .
\end{equation*}
In the same fashion, we may show that
\begin{align*}
|\mu_{j} (t, r') - \mu_{j} (t,r)| \les O \left( r^{\alpha_{j} -1 - 2p (S+1)} \right) 
\end{align*}
for $1 \leq j \leq k$. We proceed as above and conclude that
  \begin{align}
  \begin{split}
  \begin{array}{lll}
  |\mu_j (t,r) | & \les r^{\alpha_j -1  - 2p (S+1)}  ,    & j \geq {k} - S  \\
  |\mu_j(t,r)| &\les 1,     & j =  {k} - S -1 .
\end{array}  
  \end{split}
  \label{EQ8.43p}
  \end{align}
  Finally, the asymptotic decay rates for $\lambda_j (t,r)$ and $\mu_{i} (t,r)$ with $j \geq \tilde{k} - (S+1)$  and $i \geq k - (S+1)$  may be refined by $p \epsilon$.
  That follows from running the difference estimates one more time with the fact that
\begin{equation*}
|\lambda_{\tilde{k} - S -1} (t,r) | \les O (1) \comma   |\mu_{k - S-1} (t,r)|  \les O(1). 
\end{equation*}
This completes the proof of Proposition~\ref{Proposition4}.
\end{proof}

  Consequently, setting $S = k$ in Proposition~\ref{Proposition4}, we have
\begin{align}
\begin{split}
|\lambda_j (t,r) | & \les r^{\alpha_j  - 2pk} \comma j = 2, \ldots, \tilde{k} \\
|\mu_j (t,r) | & \les r^{\alpha_j -1  - 2pk} \comma j = 1, \ldots, {k} \\
|\lambda_{1} (t,r)| & \les 1 . 
\end{split}
\label{EQ8.44p}
\end{align}
Next, we  check the asymptotic decay rate of the leading coefficient $\lambda_1 (t,r)$.  
 As we  see from the statement of Proposition~\ref{Proposition2},
 it will be sufficient to obtain this result at time $t=0$. 
 For that reason, we simplify the notation and denote by 
 \begin{align}
 \begin{split}
  \lambda_j (r) & := \lambda_j (0,r) \comma 1 \leq j \leq  \tilde{k}        \\
  \mu_j (r) & := \mu_j (0,r) \comma 1 \leq j \leq k.
  \end{split}
  \label{eq8.26p}
 \end{align}

\begin{Lemma}
\label{Lem14}
There exists $\ell \in \mathbb{R}$ so that 
\begin{align}
 \left|  \lambda_1 (r)  - \ell  \right| = O\left( r^{-(d-2) p +d} \right) \inon{ as }~ r \to \infty .
 \label{eq8.ell}
\end{align}
%where the exponent $\alpha_1= \min (5p -7, 3p^2 -2p -7)$. 
\end{Lemma}

\begin{proof}
We insert the asymptotic decay rates in \eqref{EQ8.44p} into 
 the difference equation \eqref{eq8.7} for $\lambda_1$, and  we obtain 
\begin{align}
\begin{split}
|\lambda_1 (2^{n+1} r_0) - \lambda_1 (2^n r_0) | & \les r^{d+2 - pd - 2} \left(  r^{ 2p} + r^{(4 \tilde{k} +1)p - p^2 (4 k +1)} \right) \\
& \les   r^{d+2 - pd - 2 } r^{ 2p} 
\end{split}
\label{eq8.27p}
\end{align}
where $r_0 > R_1$ is fixed. Noting that the power on the right side above is negative, we once again apply
the arguments in Remark~\ref{R5}, and  we deduce that there exists 
$\ell := \lim_{n\to \infty} \lambda_1 (2^n r_0)$ satisfying
\begin{equation}
|\lambda_1 (r)  - \ell| \les r^{d- p (d-2) } .
\end{equation}
\end{proof}

In addition, we recall the difference estimates at $t=0$:
\begin{align}
\begin{split}
|\lambda_j (2^{n+1} r_0) - \lambda_j (2^n r_0)| & \les (2^n r_0)^{d+2 - p (d-2) - 2 j } \\
|\mu_j (2^{n+1} r_0) - \mu_j (2^n r_0)| & \les (2^n r_0)^{d+1 - p (d-2) - 2j} .
\end{split}
\label{EQ8.45p}
\end{align}
By \eqref{EQ8.44p}, we also know that $\lambda_j (r) \to 0$ for $2 \leq j \leq \tilde{k}$ and $\mu_j (r) \to 0$ for 
$1 \leq j \leq k$. Hence, we have
\begin{align}
\begin{split}
|\lambda_j ( r) | & \les r^{d+2 - p (d-2) - 2 j }  \comma j = 2, \ldots, \tilde{k} \\
|\mu_j ( r) | & \les r^{d+1 - p (d-2) - 2j} \comma j = 1, \ldots, k.
\end{split}
\label{EQ8.46p}
\end{align}

\begin{proof}[Proof of Proposition~\ref{Proposition2} when $d \equiv 3  \pmod 4$]  
We complete the proof of Proposition~\ref{Proposition2}
by combining Lemma~\ref{Lem14} and \eqref{EQ8.46p} with the identities  \eqref{EQ8.2p} and \eqref{EQ8.3p}. 
%$\lambda_1 (r)$,  $\lambda_2 (r)$ and $\mu (r)$ 
\end{proof}

\subsubsection{Proof of Proposition~\ref{Proposition2} when $d \equiv 1  \pmod 4$}
Now, we assume that $d \equiv 1 \pmod 4$. Recalling \eqref{EQ5.1p}, we then obtain 
\begin{equation}
k = \tilde{k} \comma d = 4k +1.
\label{EQ10.0.p}
\end{equation}
This implies that we have the same number of $\lambda_j (t,r)$
 and $\mu_j (t,r)$. Nevertheless, the proof  may still be adapted to the same line of arguments employed in the previous case. Here, we will first show that $\mu_k (t,r) \to 0$
 as $r \to \infty$ uniformly in time. This will let us pair up $\lambda_i (t,r)$ and $\mu_{i-1} (t,r)$, and prove that $(\lambda_{i} (t,r) , \mu_{i-1} (t,r)) \to 0$ as $r \to \infty$ by using the ideas demostrated in the case $d \equiv 3 \pmod 4$. We will then observe the asymptotic rate of the leading coefficient $\lambda_1 (t,r)$ at $t=0$. 
  
Throughout this subsection, the definition of $\alpha_j$ and $\eta_p$ in \eqref{EQ8.14p} and \eqref{EQ8.eta.2} will still be useful. Nevertheless, we need to keep in mind the change in numerology regarding $d, k , \tilde{k}$. In order to avoid confusion, we introduce for each $1 \leq j \leq k$,
\begin{align}
\begin{split}
\beta_j & = d+2 - pd -2j + (2 {k} +1 )p  \\
& = 4 {k} +3 - 2k p - 2j . 
\end{split}
\label{EQ10.1p}
\end{align}

Throughout this subsection, we in fact have $d \geq 9$ and $p \geq 3$. As a result,  we get
\begin{equation}
\beta_j = - 2k + 3 - (p-3) 2k - 2j <0
\end{equation}
for every $1 \leq j \leq k$. 

We begin with a review of Lemma~\ref{Lem8}. Here, the same result follows from a few adjustments in the proof. 

\begin{Lemma}
\label{Lem8p}
The projection coefficients defined in \eqref{EQ5.2p}--\eqref{EQ5.3p} satisfy the following estimates uniformly in time: 
\begin{align}
\left| \lambda_{j} (t,r) \right| & \les 1  \comma 1 \leq j \leq {k}    \label{EQ10.2p} \\
\left| \mu_j (t,r) \right| & \les 1  \comma 1 \leq j \leq  k .
\label{EQ10.3p}
\end{align}
\end{Lemma} 

\begin{proof}
Let $\epsilon >0$ be a small number. 
Similar  to the proof of Lemma~\ref{Lem8}, we first appeal to Corollary~\ref{Cor1}. We have
\begin{align*}
  |\lambda_j (t, 2^{n+1} r_0)| &  \leq  |\lambda_j (t, 2^n r_0)| 
  					+ C \delta_1 \left(   \sum_{i=1}^{{k}} |\lambda_i (t, 2^n r_0)|  (2^n r_0 )^{2i -2j} (2^n r_0)^{\frac{d+2}{2} - \frac{pd}{2} +p}   \right) \\
& \indeq + C \delta_1 \left(    \sum_{i=1}^{{k}} |\mu_i (t, 2^n r_0)|  (2^n r_0 )^{2i -2j} (2^n r_0)^{\frac{d}{2} - \frac{pd}{2} +p+2}   \right)
\end{align*}
and 
\begin{align*}
  |\mu_j (t, 2^{n+1} r_0)| &  \leq       |\mu_j (t, 2^n r_0)| 
  					+ C \delta_1 \left(   \sum_{i=1}^{{k}} |\lambda_i (t, 2^n r_0)|  (2^n r_0 )^{2i -2j} (2^n r_0)^{\frac{d}{2} - \frac{pd}{2} +p}   \right) \\
& \indeq + C \delta_1 \left(    \sum_{i=1}^{{k}} |\mu_i (t, 2^n r_0)|  (2^n r_0 )^{2i -2j} (2^n r_0)^{\frac{d}{2} - \frac{pd}{2} +p+1}   \right)
\period
\end{align*}
We adjust the definition of $H_n$ as follows: in this case, we denote by
\begin{equation*}
H_n := 
\sum_{j=1}^{{k}}  |\lambda_j (t, 2^n r_0)|  (2^n r_0 )^{2j -2 k -1}  +   \sum_{j=1}^{{k}} |\mu_j (t, 2^n r_0)|  (2^n r_0 )^{2j  - 2 k} .
\end{equation*}
This adjustment still allows us to get
\begin{align*}
\begin{split}
H_{n+1} & \leq  \left(1+  C(k+ \tilde{k}) \delta_1  \right) H_n \\
& \leq \left(1+  2 C k \delta_1  \right)^{n+1} H_0 
\end{split}
\end{align*}
by iterating in $n$. Selecting $\delta_1 >0$ to be sufficiently small so that $1 + 2C k \delta_1 < \epsilon$, we obtain
\begin{equation}
H_n \leq C (2^n r_0 )^{\epsilon}
\end{equation}
which yields
\begin{align}
\begin{split}
|\lambda_i (t, 2^n r_0) | & \leq (2^n r_0)^{2 {k} - 2i +1 + \epsilon} \comma 1 \leq i \leq k \\
|\mu_i (t, 2^n r_0)|  & \leq (2^n r_0)^{2{k} - 2i  + \epsilon} \comma 1 \leq i \leq k .
\end{split}
\label{EQ10.4p}
\end{align}
We close the proof in the same fashion. Inserting the estimates \eqref{EQ10.4p} into the difference inequalities 
\eqref{eq8.7}--\eqref{eq8.9} we obtain
\begin{align*}
|\lambda_j (t, 2^{n+1} r_0)| \leq  |\lambda_j (t, 2^n r_0)| + C (2^n r_0 )^{\beta_j + p \epsilon}
\end{align*}
and 
\begin{align*}
|\mu_j (t, 2^{n+1} r_0)| \leq  |\mu_j (t, 2^n r_0)| + C (2^n r_0 )^{\beta_j -1 + p \epsilon}.
\end{align*}
Since $\beta_j <0$, by iterating the above estimates we deduce that
\begin{align}
 |\lambda_j (t, 2^n r_0)| \les 1 \comma  |\mu_j (t, 2^n r_0)| \les 1
 \label{EQ10.5p}
\end{align}
for every $1 \leq j \leq k$. Running the difference inequalities once again for $2^n r_0 < r < 2^{n+1} r_0$, we obtain
\eqref{EQ10.2p}--\eqref{EQ10.3p}.
\end{proof}

Before we discuss the inductive step, we show that $\mu_k (t,r) \to 0$ as $r \to \infty$. 
\begin{Lemma}
\label{Lem9p}
There exists a bounded function $\rho_k (t)$ such that
\begin{equation}
|\mu_k (t,r)  - \rho_k (t)| \les O \left(  r^{\beta_k -1 } \right)
\end{equation}
uniformly in $t \in \mathbb{R}$.
\end{Lemma}

\begin{proof}
We fix $ r_0 > R_1$ and feed the  estimates in Lemma~\ref{Lem9p} into the difference inequality for $\mu_k (t,r)$. Recalling \eqref{EQ10.1p} we get
\begin{align*}
|\mu_k (t, 2^{n+1} r_0) - \mu_k (t, 2^n r_0)| & \les (2^n r_0)^{d+1 - pd - 2k} \left( \sum_{i=1}^k  (2^n r_0)^{2ip}  +
\sum_{i=1}^k (2^n r_0)^{2ip + p}  \right) \\
& \les (2^n r_0)^{d+1 - pd - 2k} (2^n r_0)^{2kp +p } \\
& \les (2^n r_0)^{\beta_k -1}. 
\end{align*}
Note that the right hand side above is dominated by the asymptotic rate of $\mu_k (t, 2^n r_0)$.  Since $\beta_k <0$, we follow the arguments in Remark~\ref{R5} and we deduce that there exists a bounded function $\rho_k (t)$ such that
\begin{equation*}
\lim_{n \to \infty} \mu_k (t, 2^n r_0) = \rho_k (t)
\end{equation*}
uniformly in time. 
From \eqref{eq8.9}, we conclude that $\lim_{r \to \infty} \mu_k (t,r) = \rho_{k} (t)$. 
\end{proof}

\begin{Lemma}
\label{Lem10p}
The limit $\rho_k (t)$ is independent of time. 
\end{Lemma}

\begin{proof}
Combining  \eqref{EQ8.2p} with Lemma~\ref{Lem8p} we may estimate
\begin{align}
|u (t,r)| \leq \sum_{i=1}^k | \lambda_i (t,r) | r^{2i - d}  \les r^{2k - d}.
\label{EQ10.6p}
\end{align}
We then expand by Lemma~\ref{Lem9p}
\begin{align}
\begin{split}
\rho_k (t_1) - \rho_k (t_2) & = \frac{1}{R} \int_R^{2R} ( \mu_k (t_1, r) - \mu_k (t_2, r) ) dr \\
& \indeq + O \left(  R^{\beta_k -1 } \right) \\
& = I_{\mu_{k}} (R) +  O \left(  R^{\beta_k -1 } \right) .
\end{split}
\end{align}
The rest of the proof follows a simpler version of the proof of Lemma~\ref{Lem11}. 
Similarly, we appeal to  Lemma~\ref{Lem15.p} to express
\begin{align*}
I_{\mu_k } (R)
= \sum_{i =1}^{k}  \frac{c_i c_j}{d - 2i - 2k} \int_{t_1}^{t_2} (I (i, k) + II (i, k ) )~dt
\end{align*}
where the formulas for $I (i, k)$ and $II(i, k)$ are as in \eqref{EQ8.8p}--\eqref{EQ8.9p}. 
By \eqref{EQ10.6p}, we compute the asymptotic rates of $I (i, k)$ and $II (i, k)$ 
as follows: 
\begin{align*}
I (i, k) &  \les O \left( R^{-2 } \right) \\
II (i, k) & \les O \left( R^{2k +1 - p (2k +1)} \right).
\end{align*}
As we analyze the decay rates, we find that 
\begin{align*}
\beta_k -1  &  < -2  \\
2k +1 - p (2k+1) & < -2 . 
\end{align*}
Therefore,
\begin{equation}
| \rho_{k } (t_1) - \rho_{k} (t_2) |  \les O \left(  R^{-2}  (1 + |t_1 - t_2|) \right) \to 0 \inon{ as } ~ R \to \infty 
\end{equation}
which implies that $\rho_{k} (t) $ must be independent of time. 
\end{proof}

\begin{Lemma}
\label{Lem13p}
The limit $\rho_k =0$. 
\end{Lemma}

\begin{proof}
The proof follows from the same arguments presented in that of Lemma~\ref{Lem13}. 
\end{proof}
Next, we state the analogous version of Proposition~\ref{Proposition4}.
\begin{Proposition}
\label{Proposition4p}
The following estimates for the projection coefficients $\lambda_j (t,r)$ and $\mu_j (t,r)$ 
hold true for $Q= 0, \ldots, k-1$,  uniformly in time:
  We have
  \begin{align}
  \begin{split}
  \begin{array}{llr}
   |\lambda_j (t,r) | & \les r^{\beta_j - 2p Q} \comma  & {k} - Q < j \leq {k} \\
  |\lambda_j (t,r)| & \les 1    \comma & 1 \leq j  \leq k-Q .
  \end{array}
  \end{split}
  \label{EQ10.7p}
  \end{align}
Similarly, we have
   \begin{align}
  \begin{split}
  \begin{array}{llr}
  |\mu_j (t,r) | & \les r^{\beta_j -1  - 2p Q} \comma    & {k} - Q \leq  j \leq {k} \\
  |\mu_j (t,r)| &\les 1  \comma & 1 \leq j <  {k} - Q .
\end{array}  
  \end{split}
  \label{EQ10.8p}
  \end{align}
\end{Proposition}

\begin{proof}[Proof of Proposition~\ref{Proposition4p}]
The case $Q=0$ is already achieved  by Lemmas~\ref{Lem8p}--\ref{Lem9p}, and the fact that $\rho_k = 0$.
For the induction step, we assume that \eqref{EQ10.7p} and \eqref{EQ10.8p} are true for $Q$ with $0 \leq Q \leq k-2$, and we prove that
\eqref{EQ10.7p} and \eqref{EQ10.8p} hold for $Q+1$. 
The remainder of the proof may be analogously  achieved by following the same line of arguments as in the proof of Proposition~\ref{Proposition4}.
We first use the induction hypothesis to improve the asymptotic rates of $\lambda_{k-Q} (t,r)$ and $\mu_{k-Q-1} (t,r)$.
\\
\textbf{Claim 1.} There exist bounded functions $\theta_{k-P} (t)$ and $\rho_{k-P-1} (t)$ such that
\begin{align*}
|\lambda_{{k} - Q} (t,r) - \theta_{k - Q} (t) | & \les O \left( r^{\beta_{(k - Q)}  - 2 p Q  } \right)   \\
|\mu_{k -Q-1} (t,r)   - \rho_{k-Q-1} (t)|  & \les O \left( r^{\beta_{(k-Q-1)} -1 - 2p Q  }\right) 
\end{align*}
uniformly in time. 

The proof of Claim~1 follows from the same arguments employed in the proof of Lemma~\ref{Lem9}. 
\\
\textbf{Claim 2.} The following asymptotics for $u (t,r)$ hold uniformly in time:
\begin{equation}
r^{- 2 (k-Q) +d } u(t,r) = \theta_{k-Q} (t) + O \left( r^{ - 2}  \right).
\label{EQ10.9p} 
\end{equation}
 Similar to the proof of Lemma~\ref{Lem12}, we  use \eqref{EQ8.2p} and expand $u(t,r)$ into a sum of  $\lambda_j (t,r) r^{2j-d}$.  The asymptotic rate 
 in \eqref{EQ10.9p} is then achieved by utilizing the estimates in the induction hypothesis as well as Claim~$1$.
\\
\textbf{Claim 3.} $\theta_{k-Q} (t) = 0$ and $\rho_{k-Q-1} (t) =0$ for all $t \in \mathbb{R}$.

The proof of Claim~$3$ is obtained by proceeding as in the previous case. We apply the analogous arguments in Lemma~\ref{Lem10}, Lemma~\ref{Lem11}, and Lemma~\ref{Lem13}. 

The last step of the proof is  the verification of \eqref{EQ10.7p}--\eqref{EQ10.8p} for $Q+1$, which may also be established as done in the case $d \equiv 3 \pmod 4$. 
\end{proof}

Next, we summarize asymptotics we achieve by Proposition~\ref{Proposition4p}. More specifically, we set $Q = k-1$, and we find
\begin{align}
\begin{split}
|\lambda_j (t,r) | & \les r^{\beta_{j}  - 2p(k-1)} \comma j = 2, \ldots, {k} \\
|\mu_j (t,r) | & \les r^{\beta_j -1  - 2p(k-1)} \comma j = 1, \ldots, {k} \\
|\lambda_{1} (t,r)| & \les 1 . 
\end{split}
\label{EQ10.10p}
\end{align}
Arguing as in Lemma~\ref{Lem14}, for  $1 \leq j \leq k$ we denote by 
\begin{equation*}
\lambda_j (r) := \lambda_j (0,r) \comma \mu_j (r) := \mu_j (0,r).
\end{equation*}
\begin{Lemma}
\label{Lem14p}
There exists $\ell \in \mathbb{R}$ so that 
\begin{align*}
 \left|  \lambda_1 (r)  - \ell  \right| = O\left( r^{-(d-2) p +d} \right) \inon{ as }~ r \to \infty .
\end{align*}
%where the exponent $\alpha_1= \min (5p -7, 3p^2 -2p -7)$. 
\end{Lemma}
\begin{proof}
Recalling the definition of $\beta_j$ in \eqref{EQ10.1p}, we insert the asymptotic rates in \eqref{EQ10.10p} into the difference inequality for $\lambda_1 (r)$, and we find
\begin{align*}
|\lambda_1 (2^{n+1} r_0) - \lambda_1 (2^n r_0) | 
 \les   r^{d+2 - pd - 2 } r^{ 2p} 
\end{align*}
where $r_0 > R_1$ is fixed. The rest of the proof is identical to that of Lemma~\ref{Lem14}.
\end{proof}

Next, we recall  the decay rates of $\lambda_j (r)$ and $\mu_j (r)$ below. We have
\begin{align}
\begin{split}
|\lambda_j ( r) | & \les r^{d+2 - p (d-2) - 2 j }  \comma j = 2, \ldots, {k} \\
|\mu_j ( r) | & \les r^{d+1 - p (d-2) - 2j} \comma j = 1, \ldots, k.
\end{split}
\label{EQ10.11p}
\end{align}

\begin{proof}[Proof of Proposition~\ref{Proposition2} when $d \equiv 1  \pmod 4$]  We recall that $\tilde{k} = k$ and $d= 4k +1$. 
We then close the proof 
by combining Lemma~\ref{Lem14p} and \eqref{EQ10.11p} with the identities  \eqref{EQ8.2p} and \eqref{EQ8.3p}. 
%$\lambda_1 (r)$,  $\lambda_2 (r)$ and $\mu (r)$ 
\end{proof}

\subsection{Step~3}
Here, we show that $\vec{u} (t) \equiv (0,0)$ and close the proof of Proposition~\ref{Proposition1}.
Recalling the asymptotic rates in Proposition~\ref{Proposition2}, 
%\begin{align*}
%& r^5 u_0 (r)  = \ell + O(r^{-5p +7})   \inon{ as }~ r \to \infty    \\
%& \int_r^{\infty} u_1(s) s~ds  = O(r^{-5p+3}) \inon{ as }~ r \to \infty \period 
%\end{align*}
we consider the cases $\ell =0$ and $\ell \neq 0$ separately. 
\begin{Lemma}
\label{Lem1}
Let $\vec{u} (t)$ and $\ell$ be as in Proposition~\ref{Proposition2}. 
Suppose $\ell =0$.  Then $\vec{u} (0) = (u_0, u_1)$ is compactly supported. 
\end{Lemma}

\begin{proof}
Assume that  $\ell =0$. Then, by \eqref{EQ8.46p} and \eqref{eq8.ell}  (or by \eqref{EQ10.11p})  we get
\begin{align*}
\begin{split}
|\lambda_j (r)| & \les O \left( r^{d+2 - (d-2)p - 2j } \right) \comma    j = 1, \ldots, \tilde{k}
 \\
 |\mu_j (r)| & \les    O \left(  r^{d+1 - (d-2)p - 2j } \right)  \comma j = 1 \ldots, k
 \end{split}
\end{align*}
for $r \geq R_1$.  Taking $r_0 \geq R_1$, we have
\begin{align}
 \sum_{j=1}^{\tilde{k}}    |\lambda_j (2^n r_0 )|  + \sum_{j=1}^k |\mu_j (2^n r_0 ) |   \les   O \left( 2^n r_0 \right)^{d - (d-2)p}
 \label{eq9.1}
\end{align}
for every $n$.  On the other hand, the difference estimates in Corollary~\ref{Cor1} yield 
\begin{align*}
\begin{split}
|\lambda_j (2^{n+1} r_0)| & \geq 
\left( 1 - C_1 \delta_1  \right) |\lambda_j (2^n r_0)| \\
& \indeq - C_1 \delta_1 (2^n r_0)^{\alpha (p,d)}  \left( \sum_{{i = 1 , i \neq j}}^{\tilde{k}}     |\lambda_i (2^n r_0)| \right)   \\
& \indeq - C_1 \delta_1 (2^n r_0)^{\alpha (p,d) -1} \left( \sum_{i=1, i \neq j}^k |\mu_i (2^n r_0)|  \right)
\end{split}
\end{align*}
where
\begin{align*}
\alpha (p,d) = \frac{d+2}{2} - \frac{pd}{2} +p + 2 k < 0 .
\end{align*} 
Additionally, we get
\begin{align}
\begin{split}
|\mu_j (2^{n+1} r_0)| & \geq 
\left( 1 - C_1 \delta_1  \right) |\mu_j (2^n r_0)| \\
& \indeq - C_1 \delta_1 (2^n r_0)^{\alpha (p,d)-1}  \left( \sum_{{i = 1 , i \neq j}}^{\tilde{k}}     |\lambda_i (2^n r_0)| \right)   \\
& \indeq - C_1 \delta_1 (2^n r_0)^{\alpha (p,d) -2} \left( \sum_{i=1, i \neq j}^k |\mu_i (2^n r_0)|  \right)
\end{split}
\end{align}
Then, setting $\delta_1>0$ sufficiently small and $r_0$ sufficiently large so  that 
$C_1 \delta_1 (1+ (k+\tilde{k}) r_0^{\alpha (p,d)}) < 1/4$, we iterate the lower bounds above to get 
\begin{align}
\sum_{j =1}^{\tilde{k}}  |\lambda_j (2^n r_0)| + \sum_{j=1}^{k} |\mu_j (2^n r_0 ) |   
\geq 
\left( 3/4 \right)^{n}
\left( \sum_{ j=1}^{\tilde{k}}  | \lambda_j ( r_0)| 
+    \sum_{j=1}^k | \mu_j ( r_0)| \right).
\label{eq9.2}
\end{align}
Combining \eqref{eq9.1}, \eqref{eq9.2}, and the fact that $p \geq 3$,  we arrive at
\begin{align*}
\sum_{ j=1}^{\tilde{k}}  | \lambda_j ( r_0)| 
+    \sum_{j=1}^k | \mu_j ( r_0)| 
 \les \frac{4^n}{3^n \cdot 2^{n } } r_0^{d - (d-2)p}
\end{align*}
for every $n \in \mathbb{N}$, which leads to $
|\lambda_j (r_0)| =  | \mu_i (r_0)| = 0$ 
for every $j=1, \ldots, \tilde{k}$ and $i = 1, \ldots, k$.
%does not work for p = 9/5 +epsilon.
It then follows from \eqref{eq5.10} and Lemma~\ref{Lem4} that
\begin{align*}
\Vert \pi^{\perp}_{r_0} \vec{u}(0) \Vert^2_{\mathcal{H} (r \geq r_0)}  
\les  {r_0}^{- 2 \beta}   \Vert \pi_{r_0} \vec{u}(0) \Vert^{2 p}_{\mathcal{H} (r \geq r_0)}  =0.
\end{align*}
Therefore, 
\begin{align*}
\Vert \vec{u}(0) \Vert_{\mathcal{H} (r \geq r_0)} =0.
\end{align*}
In other words $(\partial_r u_0, u_1)$ is compactly supported. Since we have
\begin{align*}
\lim_{r \to \infty} u_0 (r) = 0
\end{align*}
we may conclude that $\vec{u} (0)$ is compactly supported.
\end{proof}
\begin{Lemma}
\label{Lem2}
Let $\vec{u} (t)$ and $\ell$ be as in Proposition~\ref{Proposition2}. 
Suppose $\ell =0$.  Then $\vec{u} (0) = (0, 0)$. 
\end{Lemma}

\begin{proof}
Assuming $\ell = 0$, we deduce from Lemma~\ref{Lem1} that the inital data $(u_0, u_1)$ must be compactly supported. 
Furthermore, if $(u_0, u_1) \neq (0,0)$, then there exists a positive radius $\rho_0$ such that
\begin{align*}
\rho_0 := \inf \{\rho: \norm{\vec{u} (0) }_{\mathcal{H} (r \geq \rho)} = 0 \}. 
\end{align*}
Let $\delta_1$  be as in \eqref{eq8.4}--\eqref{eq8.6}. Additionally, we take a small number $\epsilon >0$ , to be determined below, and find $\rho_1= \rho_1 (\epsilon)$ with 
$\frac{1}{2}  \rho_0 < \rho_1 < \rho_0$ such that 
\begin{align*}
0 < \norm{\vec{u} (0) }^{p-1}_{\mathcal{H} (r \geq \rho_1)}  < \epsilon \leq \delta_1 
\period
\end{align*}
By Lemma~\ref{Lem5}, we have 
\begin{align}
\begin{split}
& \norm{\vec{u} (0) }^2_{\mathcal{H} (r \geq R)} \cong  \sum_{i=1}^{\tilde{k}} \left(   \lambda_i ( R) R^{2i - \frac{d+2}{2}} \right)^2
+ \sum_{i=1}^k \left(   \mu_i (R) R^{2i - \frac{d}{2}} \right)^2 \\
&  \indeq \indeq +  \int_{R}^{\infty} \left(    \sum_{i=1}^{\tilde{k}} \left(  \partial_r \lambda_i (r) r^{2i - \frac{d+1}{2}}\right)^2
+   \sum_{i=1	}^{k} \left( \partial_r \mu_i (r) r^{2i - \frac{d-1}{2}} \right)^2 \right)~dr 
\period
\end{split}
\label{eq9.3}
\end{align}
Note that by setting $ R = \rho_0$ above  we get 
\begin{align}
\lambda_j (\rho_0)  = \mu_i (\rho_0) =0 \comma 1 \leq j \leq \tilde{k}, ~1 \leq i \leq k.
\label{eq9.4}
\end{align}
Also, by Lemma~\ref{Lem6} we may bound the integral on the right hand side above  
as follows:  
\begin{align}
\begin{split}
& 
\int_{R}^{\infty}
\sum_{i=1}^{\tilde{k}}  \left( \partial_r \lambda_i (r) r^{2i - \frac{(d+1)}{2}} \right)^2 
+   \sum_{i=1}^k \left( \partial_r \mu_i (r) r^{2i - \frac{(d-1)}{2}}  \right)^2 
~dr
 \\ &
\indeq \indeq \indeq  \les 
 \sum_{i=1}^{\tilde{k}} \lambda_i^{2p} (R) R^{  (4i -2d) p+ d+2} + \sum_{i=1}^{k} \mu_i^{2p} (R) R^{ (4i - 2 d+2 ) p  +d+2   }  \period
 \end{split}
\label{eq9.6}
\end{align}
We then argue as in the proofs of Lemma~\ref{Lem7} and Corollary~\ref{Cor1}, and estimate the differences 
\begin{align}
\begin{split}
 \left|  \lambda_j (\rho_1) - \lambda_j (\rho_0) \right| 
& \les  \epsilon  \left(  \sum_{i=1}^{\tilde{k}} |\lambda_i (\rho_1)| {\rho_1}^{2i - 2j } {\rho_1}^{\frac{d+2}{2} - \frac{pd}{2} +p}    + \sum_{i=1}^k |\mu_i (\rho_1)| 
{\rho_1}^{2i - 2j} 
{\rho_1}^{\frac{d}{2} - \frac{pd}{2} +p +2}  \right)  \\
\left| \mu_j (\rho_1) - \mu_j (\rho_0 ) \right| 
 & \les \frac{ \epsilon}{\rho_1}   \left(  \sum_{i=1}^{\tilde{k}} |\lambda_i (\rho_1)| {\rho_1}^{2i - 2j } {\rho_1}^{\frac{d+2}{2} - \frac{pd}{2} +p}   
  + \sum_{i=1}^k |\mu_i (\rho_1)| {\rho_1}^{2i - 2j} 
{\rho_1}^{\frac{d}{2} - \frac{pd}{2} +p +2}  \right) \period
\end{split}
\label{eq9.5}
\end{align}
Next, we set 
\begin{align*}
H  = \sum_{j=1}^{\tilde{k}} {\rho_1}^{2j} |\lambda_j (\rho_1) | + \sum_{j=1}^{k} {\rho_1}^{2j +1} |\mu_j (\rho_1)| . 
\end{align*}
Recalling \eqref{eq9.4} and the fact that $\frac{1}{2}  \rho_0 < \rho_1 < \rho_0$ we may rewrite equation \eqref{eq9.5} as
\begin{align*}
H \leq C \epsilon H
\end{align*}
where the constant $C$ depends only on $\rho_0$  and the uniform implicit  constant in \eqref{eq9.5} due to $\les$. As $\rho_0$ is fixed, 
we may select $\epsilon  \in (0, C^{-1})$ and deduce that $H=0$. By setting $R= \rho_1$ in \eqref{eq9.3} and \eqref{eq9.6} we find that
\begin{align*}
\norm{ \vec{u} (0)}_{\mathcal{H} (r \geq \rho_1)} =0.
\end{align*}
However, this leads to a contradiction as $\rho_1 < \rho_0$. 
\end{proof}

\begin{Lemma}
\label{Lem3}
Let $\vec{u} (t)$ and $\ell$ be as in Proposition~\ref{Proposition2}. 
Then, $\ell =0$. 
\end{Lemma}

In order to prove Lemma~\ref{Lem3} we show that the case $\ell \neq 0$ leads to a contradiction. If the limit $\ell $ is nonzero, 
then we may consider the difference $\vec{u} (t) - (Z_{\ell}, 0)$, where $Z_{\ell}$ is the corresponding stationary solution  constructed 
in Proposition~\ref{prop_singular_solns}. Below, we will argue that the results we obtained in Step~$2$ leads to $u(t,r) = Z_{\ell} (r)$, which gives us  contradiction
since $Z_{\ell} \notin \dot{H}^{s_p} (\mathbb{R}^d)$. 

Recalling \eqref{SS_eq1}--\eqref{SS_eq2}, we define $\vec{\omega_{\ell}} (0) = (\omega_{\ell,0}, \omega_{\ell,1})$ by
\begin{align}
\begin{split}
\omega_{\ell,0} & := u_0 (r) - Z_{\ell} (r) \\
\omega_{\ell,1} & := u_1 (r)
\end{split}
\label{eq10.1}
\end{align}
and we consider for all $t \in \mathbb{R}$
\begin{align}
\begin{split}
 \vec{\omega}_{\ell} (t)  & =  (\omega_{\ell} (t,r), \partial_t \omega_{\ell} (t,r) ) \\
   & : = ( u (t,r) - Z_{\ell} (r), \partial_t u (t,r) ).
\end{split}
\label{eq10.2}
\end{align}
Note that we may directly utilize the asymptotic decay rates obtained for $\vec{u} (0)$ and $Z_{\ell}$ and estimate
\begin{align}
\begin{split}
& r^{d-2} \omega_{\ell,0} (r) = O (r^{-(d-2)p +d })  \inon{as~}   r \to \infty \\
& \int_r^{\infty} \omega_{\ell,1} (\rho) \rho^{2i -1} ~d\rho = O (r^{-(d-2)p +2i+1})  \inon{as~} r \to \infty.
\end{split}
\label{eq10.3}
\end{align} 
Next, we check the equation for $\vec{\omega}_{\ell} (t,r)$. Since $\vec{u}$ and $Z_{\ell}$ are solutions to \eqref{wave equations}
and \eqref{SS_eq1} respectively, we get
\begin{align}
\begin{split}
\partial_{tt} \omega_{\ell} - \partial_{rr} \omega_{\ell} - \frac{d-1}{r} \partial_{r} \omega_{\ell}   
% & = \omega_{\ell} \left( u^4 + u^3 Z_{\ell} + u^2 Z_{\ell}^2 + u Z_{\ell}^3 + Z_{\ell}^4    \right) \\
 & = | \omega_{\ell} + Z_{\ell} |^{p-1} ( \omega_{\ell} + Z_{\ell}) - |Z_{\ell}|^{p-1} Z_{\ell}
\end{split}
\label{eq10.4}
\end{align}
As $Z_{\ell}$ is stationary, $\vec{\omega}_{\ell}$ verifies the latter conclusion of Corollary~{\ref{cor8.1}}, i.e.,  we simply get
\begin{equation}
\limsup_{t \to +\infty} \Vert \vec{\omega}_{\ell} (t) \Vert_{\mathcal{H} (r \geq R +|t|)} = \limsup_{t \to - \infty} \Vert \vec{\omega}_{\ell} (t) \Vert_{\mathcal{H} (r \geq R +|t|)}   =  0.
\label{eq10.5}
\end{equation}

\begin{Lemma}
\label{Lem16}
Suppose $\ell \neq 0$, and let  $\vec{\omega}_{\ell} (t)$ be defined as in \eqref{eq10.2}. Then, we must have $\vec{\omega}_{\ell} (0) \equiv (0,0)$.
\end{Lemma}
The proof of Lemma~\ref{Lem16} follows  from the same line of arguments presented in the first two steps. Firstly, as done in Step~1,  we will obtain an analogous version of Lemma~\ref{Lem4} and express that in terms of projection coefficients of $\vec{\omega}_{\ell} (t)$, which will then lead to corresponding difference estimates.  As we already established the asymptotic decay of $\vec{\omega}_{\ell} (0)$
in \eqref{eq10.3}, we will close the proof by showing that $\vec{\omega}_{\ell} (0) $ must be compactly supported. 
Below, we will outline how to adapt the results of Step~1 and Step~2 for $\vec{\omega}_{\ell} (t)$. 

In order to prove a version of the estimate \eqref{eq7.3}, we take a second look at the Cauchy problem in Lemma~\ref{Lem17}. Following the set-up in \eqref{eq7.4}, we define
$V (t, x) = \chi \left( \frac{x}{\tilde{R}_0 + |t|} \right) Z_{\ell} (x)$ for some large $\tilde{R}_0 >0$. Then, $V$ satisfies the assumptions of Lemma~\ref{Lem17} with $I = \mathbb{R}$ and $r_0 = \tilde{R}_0$. Letting $r_1 > \tilde{R}_0$ such that 
\begin{align*}
\norm{\vec{\omega}_{\ell} (0)}_{\mathcal{H} (r \geq r_1 )} \leq \delta_0 r_1^{\beta / (p-1)}
\end{align*}
 we obtain
\begin{align}
\sup_{t \in \mathbb{R}} \Vert \vec{\omega}_{\ell} (t) - S(t) (\omega_{\ell,0}, \omega_{\ell,1}) \Vert_{\mathcal{H}} \leq 
\frac{1}{2^d \cdot 100} \Vert (\omega_{\ell,0}, \omega_{\ell,1})  \Vert_{\mathcal{H}} \period
\label{eq10.6}
\end{align}

Having obtained the estimate \eqref{eq10.6} above, we proceed to adjust the result of Lemma~\ref{Lem4}. 
\begin{Lemma}
\label{Lem18}
There exists $\tilde{R}_0>0$ such that for all $R> \tilde{R}_0$ we have
\begin{align}
\Vert \pi^{\perp}_R \vec{\omega}_{\ell} (t) \Vert^2_{\mathcal{H} (r \geq R)}  
\les \frac{1}{4^d \cdot {10}^{4}} \Vert \pi_R \vec{\omega}_{\ell} (t) \Vert^{2}_{\mathcal{H} (r \geq R)}  
\period
\label{eq10.7}
\end{align} 
\end{Lemma}
We omit the proof of Lemma~\ref{Lem18} since it is identical to the proof of Lemma~\ref{Lem4}. Namely, we follow the same procedure and use 
the estimate \eqref{eq10.6} instead of \eqref{eq7.9}, which leads to the power change on the right hand side of \eqref{eq10.7}.

Let us remind that the orthogonal projections in \eqref{eq10.7} will be of the following form: 
\begin{align}
\begin{split}
\pi_R { \vec{\omega}_{\ell} (t,r) }  
& = \left(     \sum_{i=1}^{\tilde{k}}  \lambda_{\ell , i} (t,R) r^{2i -d} ,~ \sum_{j=1}^{k} \mu_{\ell , j} (t,R) r^{2j -d}  \right)
\\
\pi^{\perp}_R { \vec{\omega}_{\ell} (t, r)} & = \left(  \omega_{\ell} (t,r) - \sum_{i=1}^{\tilde{k}}  \lambda_{\ell, i} (t,R) r^{2i -d}, ~    \partial_t \omega_{\ell} (t, r) 
-  \sum_{j=1}^{k} \mu_{\ell, j} (t,R) r^{2j -d}    \right)
\end{split}
\label{eq10.22}
\end{align}
where  $\lambda_{\ell,i}(t, R)$ and $\mu_{\ell,j} (t, R)$ are defined using the formulas in \eqref{EQ5.2p}--\eqref{EQ5.3p}. Note that  these projection coefficients must be  adapted to $\vec{\omega}_{\ell} (t)$. However, since $\partial_t \omega_{\ell} (t, r)= \partial_t u(t, r)$, the formula \eqref{EQ5.3p} gives us  $\mu_{\ell, j} (t, R) = \mu_j  (t, R)$.  
We refer the reader to Section~\ref{sec: orthogonal proj} for a comparison. 

Recalling the decay rates of $(\omega_{\ell,0} (r), \omega_{\ell,1} (r))$ in \eqref{eq10.3}, we immediately deduce the asymptotics for 
$\lambda_{\ell,j} (r)$ and $\mu_{\ell,j} (r)$. Namely, we get 
\begin{align}
\begin{split}
|\lambda_{\ell,j} (r)| & = O \left( r^{- (d-2)p  + d +2 - 2j}  \right)   \comma 1 \leq j \leq \tilde{k}
\\
| \mu_{\ell, j} (r)  | & = O\left( r^{-(d-2)p + d+1 - 2j}  \right)   \comma 1 \leq j \leq k .
\end{split}
\label{eq10.8}
\end{align}

Next, we apply the exact same arguments in the proof of Lemma~\ref{Lem1} to prove that $\vec{\omega}_{\ell} (0)$ is compactly supported.  
\begin{Lemma}
\label{Lem19}
Let $(\omega_{\ell,0}, \omega_{\ell,1} )$ be as in \eqref{eq10.2}. Then $( \partial_r \omega_{\ell,0},  \omega_{\ell,1})$ is compactly supported. 
\end{Lemma}

\begin{proof}[Proof of Lemma~\ref{Lem19}] First, we rewrite the estimate \eqref{eq10.7} at $t=0$ in terms of $\lambda_{\ell,j} (r)$  and $\mu_{\ell, j} (r)$. Recalling Lemma~\ref{Lem5},
 we get
\begin{align}
\begin{split}
& 
\int_{R}^{\infty}\sum_{i=1}^{\tilde{k}}  \left( \partial_r \lambda_{\ell, i} ( r) r^{2i - \frac{(d+1)}{2}} \right)^2 
+   \sum_{i=1}^k \left( \partial_r \mu_{\ell, i} (r) r^{2i - \frac{(d-1)}{2}}  \right)^2 
~dr
 \\ &
\indeq \indeq \indeq  \les 
\frac{1}{4^d  \cdot 10^4}
\left(
 \sum_{i=1}^{\tilde{k}} \lambda_{\ell, i}^{2} (R) R^{  4i - (d+2) } + \sum_{i=1}^{k} \mu_{\ell, i}^{2} (R) R^{ 4i -  d  }   \right) 
\label{eq10.10}
\end{split}
\end{align}
for all $R > \tilde{R}_0$.
We argue exactly as in the proof of Lemma~\ref{Lem7} to obtain the difference estimates from \eqref{eq10.10}. For all $\tilde{R}_0  \leq r \leq r' \leq 2 r$  we have  
\begin{align}
\begin{split}
\left| \lambda_{\ell, j} (r) - \lambda_{\ell, j} (r') \right|^2 &  \leq
\left( \int_r^{r'}   \left(      \partial_s \lambda_{\ell, j} (s) s^{2 j - \frac{(d+1)}{2}}   \right)^2ds \right)  \left(  \int_r^{r'} s^{d+1 - 4i} ds  \right)
\\
&  \leq 
\frac{(2r)^{d +2- 4j}}{4^d 10^4} 
\left(
 \sum_{i=1}^{\tilde{k}} \lambda_{\ell, i}^{2} (r) r^{  4i - (d+2) } + \sum_{i=1}^{k} \mu_{\ell, i}^{2} (r) r^{ 4i -  d  }   \right) 
 \end{split}
\label{eq10.11}
\end{align}
with $1 \leq j \leq \tilde{k}$.
Similarly,  for all $\tilde{R}_0  \leq r \leq r' \leq 2 r$    we find
\begin{equation}
\left| \mu_{\ell, j} (r) - \mu_{\ell, j} (r') \right|^2 \leq 
\frac{(2r)^{d- 4j}}{4^d \cdot  10^4} \left(  \sum_{i=1}^{\tilde{k}} \lambda_{\ell, i}^{2} (r) r^{  4i - (d+2) } + \sum_{i=1}^{k} \mu_{\ell, i}^{2} (r) r^{ 4i -  d  }  \right) 
\label{eq10.13}
\end{equation}
where  $1 \leq j \leq k$.
Next, we define
\begin{equation}
{H} (r) = \sum_{j =1}^{\tilde{k}} 
\lambda_{\ell, j} (r)   r^{2j - \frac{(d+2)}{2}} + \sum_{j=1}^{k} \mu_{\ell, j} (r) r^{2j -  \frac{d}{2}}.
\end{equation}
Selecting $r_0 > \tilde{R}_0$ we 
combine the inequalities \eqref{eq10.11}--\eqref{eq10.13} to obtain
\begin{align*}
\left| {H} (2^{n+1} r_0) - {H} (2^{n} r_0) \right| \leq \frac{1}{4} \left| {H} (2^n r_0) \right)| 
\period
\end{align*}
This implies that
\begin{align*}
\left| {H} (2^{n+1} r_0)  \right| \geq \frac{3}{4} \left| {H} (2^{n} r_0) \right|. 
\end{align*}
Via iteration on $n$ we deduce
\begin{align}
\left| {H} (2^{n} r_0)  \right| \geq  \left( \frac{3}{4} \right)^n  \left| {H} ( r_0) \right|. 
\label{eq10.14}
\end{align}
However, the asymptotic decay rates in \eqref{eq10.8}  yield
\begin{align}
\left| { H} (2^n r_0) \right| \les (2^n r_0)^{- (d-2) p + \frac{d+2}{2}} . 
\label{eq10.15}
\end{align}
By \eqref{eq10.14}--\eqref{eq10.15} we get
\begin{align*}
 \left| {H} (r_0) \right|   \les   \frac{4^{n}   r_0^{- p (d-2) +\frac{d+2}{2}} }{ 3^n \cdot 2^{n  ( (d-2)p  - \frac{d+2}{2} )  }} .
\end{align*}
Letting $n \to \infty$ we deduce that ${H} (r_0) = 0$.  Going back to \eqref{eq10.10} and using the fact that 
\begin{equation}
\lambda_{\ell, j} (r_0)  = \mu_{\ell, i} (r_0) =0 \comma 1 \leq j \leq \tilde{k}, ~1 \leq i \leq k
\end{equation}
we obtain 
\begin{align*}
\int_{r_0}^{\infty} \left( \partial_r \lambda_{\ell, 1} (r) r^{-2} \right)^2 
+ \left( \partial_r \lambda_{\ell, 2} (r) \right)^2 
+  \left( \partial_r \mu (r) r^{-1} \right)^2~dr = 0. 
\end{align*}
Hence, 
\begin{align*}
\norm{\vec{\omega}_{\ell} (0)}_{\mathcal{H}(r \geq r_0)}^{2} = \Vert \pi^{\perp}_{r_0} \vec{\omega}_{\ell} (0) \Vert^2_{\mathcal{H} (r \geq r_0)}  
       +  \Vert \pi_{r_0} \vec{\omega}_{\ell} (0) \Vert^{2}_{\mathcal{H} (r \geq r_0)}  =0 
\end{align*}
which proves that $( \partial_r \omega_{\ell,0},  \omega_{\ell,1})$ is compactly supported. 
\end{proof}

Finally we proced with the proof of Lemma~\ref{Lem16}. 
\begin{proof}[Proof of Lemma~\ref{Lem16}] We follow the same argument used in the proof of Lemma~\ref{Lem2}. By the way of contradiction, we assume that
$( \partial_r \omega_{\ell,0},  \omega_{\ell,1}) \neq (0,0)$, and define
\begin{align}
\rho_0 := \inf \{\rho: \norm{\vec{\omega}_{\ell} (0) }_{\mathcal{H} (r \geq \rho)} = 0 \}. 
\label{eq10.16}
\end{align} 
By hypothesis, we get  $\rho_0 >0$ and we deduce 
\begin{equation}
\lambda_{\ell, j} (\rho_0) = \mu_{\ell, i} (\rho_0)  =0 \comma 1 \leq j \leq \tilde{k}, ~1 \leq i \leq k.
\label{eq10.21}
\end{equation}
We then take $\rho_1 \in (\frac{\rho_0}{2}, \rho_0)$ such that
\begin{align}
\norm{\vec{\omega}_{\ell} (0) }_{\mathcal{H} (r \geq \rho_1)} < {\delta}_2 < \delta_0 \rho_1^{\beta / (p-1)} . 
\label{eq10.17}
\end{align}
Above, we select ${\delta}_2 >0 $  sufficiently small that  \eqref{eq10.6} holds. 
Thus, the second inequality in \eqref{eq10.17} guarantees that Lemma~\ref{Lem18} holds with $R= \rho_1$. 
Reformulating that in terms of the projection coefficients, we get
\begin{align}
\begin{split}
& 
\int_{\rho_1}^{\infty}\sum_{i=1}^{\tilde{k}}  \left( \partial_r \lambda_{\ell, i} ( r) r^{2i - \frac{(d+1)}{2}} \right)^2 
+   \sum_{i=1}^k \left( \partial_r \mu_{\ell, i} (r) r^{2i - \frac{(d-1)}{2}}  \right)^2 
~dr
 \\ &
\indeq \indeq \indeq  \les 
\frac{1}{4^d \cdot 10^4}
\left(
 \sum_{i=1}^{\tilde{k}} \lambda_{\ell, i}^{2} (\rho_1) {\rho_1}^{  4i - (d+2) } + \sum_{i=1}^{k} \mu_{\ell, i}^{2} (\rho_1) {\rho_1}^{ 4i -  d  }   \right) 
\label{eq10.18}
\end{split}
\end{align}
Once again, we use the fundamental theorem of calculus to control the difference $|\lambda_{\ell, i} (\rho_1) - \lambda_{\ell, i} (\rho_0) |$  and $|\mu_{\ell, i} (\rho_1) - \mu_{\ell, i} (\rho_0)|$ in terms of \eqref{eq10.18}.  Repeating the arguments in \eqref{eq10.11}--\eqref{eq10.13}  we then  get
\begin{align}
\begin{split}
\left| \lambda_{\ell, j} (\rho_1) - \lambda_{\ell, j} (\rho_0 ) \right|^2 &  \leq
\frac{{\rho_0}^{d +2- 4j} - {\rho_1}^{d+2 - 4j}  }{4^d \cdot 10^4} 
\left(
 \sum_{i=1}^{\tilde{k}} \lambda_{\ell, i}^{2} (\rho_1) {\rho_1}^{  4i - (d+2) } + \sum_{i=1}^{k} \mu_{\ell, i}^{2} (\rho_1) {\rho_1}^{ 4i -  d  }   \right) 
\end{split}\label{eq10.19}
\end{align}
 for every $1 \leq j \leq \tilde{k}$.
Likewise, for every $1 \leq j \leq k$
\begin{equation}
\left| \mu_{\ell, j} (\rho_1) - \mu_{\ell, j} (\rho_0) \right|^2 \leq 
\frac{{\rho_0}^{d - 4j} - {\rho_1}^{d - 4j}  }{4^d \cdot  10^4} 
\left(
 \sum_{i=1}^{\tilde{k}} \lambda_{\ell, i}^{2} (\rho_1) {\rho_1}^{  4i - (d+2) } + \sum_{i=1}^{k} \mu_{\ell, i}^{2} (\rho_1) {\rho_1}^{ 4i -  d  }   \right) .
\label{eq10.20}
\end{equation}

Combining \eqref{eq10.19}-\eqref{eq10.20},  and noting \eqref{eq10.21} we estimate
\begin{align*}
\sum_{j=1}^{\tilde{k}} \left| \lambda_{\ell, j}(\rho_1) \right|^2
+
\sum_{j=1}^{k} \left| \mu_{\ell, j} (\rho_1)  \right|^2
\leq 
(\rho_0 - \rho_1) \tilde{C} 
\left(
\sum_{j=1}^{\tilde{k}} \left| \lambda_{\ell, j}(\rho_1) \right|^2
+
\sum_{j=1}^{k} \left| \mu_{\ell, j} (\rho_1)  \right|^2
\right)
\end{align*}
where $\tilde{C}>0 $ depends only on $\rho_0, k , \tilde{k}$ as we have  $ \rho_1 \in (\frac{\rho_0}{2}, \rho_0)$.  Finally, selecting $\rho_1$ so that 
\[
0 <
(\rho_0 -\rho_1)  \leq  \frac{1}{ 2 \tilde{C}}
\]
we arrive at the conclusion that
\[
 \lambda_{\ell, j} (\rho_1)
 =
  \mu_{\ell, i} (\rho_1)
  =0 \comma 
  1\leq j \leq \tilde{k}, ~ 1\leq i \leq k .
\]
By \eqref{eq10.18} and \eqref{eq10.22}, we then have
\[
\norm{\vec{\omega}_{\ell} (0)}_{\mathcal{H} (r \geq \rho_1)} = 0 
\]
which contradicts the definition of $\rho_0$ since $\rho_1 < \rho_0$.  Therefore, 
$(\partial_r \omega_{\ell, 0} ,  \omega_{\ell, 1} ) \equiv (0,0)$. 
Since $\omega_{\ell,0} (r) \to 0$ as $r \to \infty$, we must have $(\omega_{\ell, 0} ,  \omega_{\ell, 1} ) \equiv (0,0)$. 
\end{proof} 

\begin{proof}[Proof of Proposition~\ref{Proposition1}] We may now close the proof of Proposition~\ref{Proposition1}
by tracing our steps in Section~{\ref{sec: rigidity}}. Let $\vec{u} (t)$ be a solution of \eqref{wave equations} as in 
Proposition~\ref{Proposition1}. By Proposition~{\ref{Proposition2}}, there exists $\ell \in \mathbb{R}$ so that
\begin{align*}
& r^{d-2} u_0 (r)  = \ell + O(r^{-(d-2)p +d})   \inon{ as }~ r \to \infty    \\
& \int_r^{\infty} u_1(s) s^{2i -1}~ds  = O(r^{- (d-2)p +2i +1}) \inon{ as }~ r \to \infty \period 
\end{align*}
If $\ell$ is zero, then Lemma~\ref{Lem2} shows that $\vec{u} (0) = (0,0)$ and in turn verifies Proposition~\ref{Proposition1}.
 On the other hand, if $\ell$ is nonzero, by Lemma~\ref{Lem16} we get
$\vec{u} (0) = (Z_{\ell}, 0)$, where $Z_{\ell}$ is the singular stationary solution constructed in Section~\ref{sec:stationary solutions}. 
Finally, this yields the desired contradiction eliminating the case $\ell \neq 0$ 
since $Z_{\ell}$ is a nonzero solution to \eqref{SS_eq1} with $Z_{\ell} \notin \dot{H}^{s_p} (\mathbb{R}^d)$ and $\vec{u} (0) \in \dot{H}^{s_p} \times \dot{H}^{s_{p} -1} (\mathbb{R}^d)$. 
\end{proof}

\section{Appendix: Preliminary Results From Harmonic Analysis}
In this section, we will give an overview of some prelimininary results from the Littlewood-Paley theory, which will yield the proofs of Leibniz rule and the $C^1$ chain rule in Section~\ref{sec:A review of the Cauchy problem}.  

Recalling the definitions of the Littlewood-Paley multipliers from Section~\ref{sec:A review of the Cauchy problem}, we state the Bernstein inequalities. The version stated below is included in the book \cite[Appendix~A]{Tao06}. 
\begin{Lemma}[Bernstein's inequalities]
\label{Lem22}
Let $s \geq 0$ and $1 \leq p \leq q \leq \infty$. For $f: \mathbb{R}^d \to \mathbb{R}$, we have
\begin{align*}
\Vert  P_{\geq N} f   \Vert_{L^p (\mathbb{R}^d) } &  \les N^{-s}  \Vert D^s P_{\geq N} f   \Vert_{L^p (\mathbb{R}^d)} \\
\Vert P_{\leq N}  D^s f \Vert_{L^p (\mathbb{R}^d) } & \les N^s  \Vert P_{\leq N} f  \Vert_{L^p (\mathbb{R}^d)} \\
\Vert P_{ N}  D^{\pm s} f \Vert_{L^p (\mathbb{R}^d )} & \cong N^{\pm s }  \Vert P_{ N} f  \Vert_{L^p (\mathbb{R}^d)} \\
\Vert  P_{\leq N} f  \Vert_{L^q (\mathbb{R}^d)} & \les N^{\frac{d}{p} - \frac{d}{q}}  \Vert P_{\leq N } f \Vert_{L^p (\mathbb{R}^d)} \\
\Vert P_N f \Vert_{L^q (\mathbb{R}^d)} & \les   N^{\frac{d}{p} - \frac{d}{q}} \Vert P_N f  \Vert_{L^p (\mathbb{R}^d)} 
\end{align*}
where the implicit constants depend on $p,s,d$ in the first three inequalities and on $p, q, d$ in the latter two inequalities. 
\end{Lemma}

Note that the definition of $P_N$ leads to the telescoping identities 
\begin{equation*}
P_{\leq N} f = \sum_{M \leq N} P_M f   \comma    P_{>N} f  = \sum_{M > N} P_M f  \comma           f  = \sum_{N} P_N f
\end{equation*}
for every $f$ belonging to the Schwarz class, which in turn implies that 
\begin{equation}
\norm{f}_{L^q (\mathbb{R}^d)} \leq \sum_N \norm{P_N f}_{L^q (\mathbb{R}^d)}
\label{A.5} 
\end{equation}
by the triangle inequality. Also, since $P_{\leq N}$ is a convolution operator, we have
\begin{equation}
\Vert P_{\leq N} f \Vert_{L^q (\mathbb{R}^d) } \les \Vert N^d \check{\varphi} (N \cdot) \Vert_{L^1 (\mathbb{R}^d)} 
\Vert f \Vert_{L^q (\mathbb{R}^d)}
\les \Vert f \Vert_{L^q (\mathbb{R}^d)}
\label{A.6}
\end{equation}
by Young's inequality. 
\begin{proof}[Proof of Lemma~\ref{Leibniz Rule}]
In the case $s=0$, the claim follows from H\"older's inequality. In the general case $s >0$, the proof is similar to the proof of Lemma~A$.8$ in \cite[Appendix~A]{Tao06}.
First, recalling the definition of the Besov norm in \eqref{EQ_Besov1} we split
\begin{align}
\begin{split}
\Vert P_N (fg) \Vert_{L^r_x} & \les \Vert P_{N} ((P_{< N/8} f ) g) \Vert_{L^r_x} + \sum_{M > N/8} \Vert P_N (( P_M f ) g ) \Vert_{L^r_x}  \\
& \cong I_1 + I_2. \label{A.7}
\end{split}
\end{align}
Following the arguments in \cite[Lemma~A.8]{Tao06}, we find that
\begin{align*}
I_1 \les \Vert f \Vert_{L^{q_1}_x} \sum_{M= N/4}^{4N} \Vert P_M g \Vert_{L^{q_2}_x} 
\end{align*}
which then leads to
\begin{align}
\sum_{N} N^{2s} I_1^2 &  \les \Vert f \Vert^2_{L^{q_1}_x} \left( \sum_{N} N^{2s} \sum_{M= N/4}^{4N} \Vert P_M g \Vert^2_{L^{q_2}_x} \right) \\
& \les \Vert f \Vert^2_{L^{q_1}_x} \Vert g \Vert^2_{\dot{B}^{s}_{ q_2, 2}}
\label{A.8}
\end{align}
by the definition of the space $\dot{B}^{s}_{ q_{2} , 2}$. 

Next, we estimate $I_2$. By \eqref{A.6} and H\"older, we have
\begin{align*}
I_2 &  \les \sum_{M > N/8} \Vert (P_M f) g \Vert_{L^r_x} \\
& \les \sum_{M > N/8} \Vert g \Vert_{L^{p_2}_x}  \Vert P_M f \Vert_{L^{p_1}_x} .
\end{align*}
We then apply Cauchy-Schwarz inequality, and use that $N$, $M$ are dyadic integers, which yields
\begin{align*}
N^{2s} I_2^{2} & \les \Vert g \Vert^2_{L^{p_2}_x} N^{2s} \left( \sum_{M > N/8} M^{-s} M^{s} \Vert P_M f \Vert_{L^{p_1}_x} \right)^2 \\
& \les \Vert g \Vert^2_{L^{p_2}_x}  \left( \sum_{M > N/8} \frac{N^s}{M^{s}} \right) \left( N^{s} \sum_{M > N/8} M^{-s} M^{2s} \Vert P_M f  \Vert^2_{L^{p_1}_x} \right) \\
& \les  \Vert g \Vert^2_{L^{p_2}_x} \left( \sum_{M > N/8} \frac{N^s}{M^s} M^{2s} \Vert P_M f  \Vert^2_{L^{p_1}_x} \right) .
\end{align*}
Summing in $N$, we obtain
\begin{align}
\sum_N N^{2s} I_2^2 \les
 \Vert g \Vert^2_{L^{p_2}_x} \sum_N \sum_{M > N/8} \frac{N^s}{M^s}  M^{2s} \Vert P_M f  \Vert^2_{L^{p_1}_x} .
 \label{A.9}
\end{align}
Note that by changing the order of summation we estimate the second factor above as
\begin{align*}
 \sum_{N} \sum_{M > N/8} \frac{N^s}{M^s}  M^{2s} \Vert P_M f  \Vert^2_{L^{p_1}_x} 
& \les \sum_M  M^{2s} \Vert P_M f  \Vert^2_{L^{p_1}_x}  
 + \sum_M \left( M^{2s} \Vert P_M f  \Vert^2_{L^{p_1}_x}   \sum_{N < M/4 } \frac{N^s}{M^s}  \right) \\
 & \les M^{2s} \Vert P_M f  \Vert^2_{L^{p_1}_x}  
\end{align*}
By \eqref{A.7}, combining \eqref{A.8} and \eqref{A.9} we obtain the result.
\end{proof}

Next, we recall Fa\`{a} di Bruno's formula from \cite[Corollary~2.10]{CS96} Let $\beta$ be a multi-index with $|\beta|=k$. If $g \in C^{k} (O, \mathbb{R})$, where $O$ is an open subset of $\mathbb{R}^d$ and 
if $f \in C^{k} (U, \mathbb{R})$ where $U$ is a neighborhood of $g (O)$, then on the set $O$
\begin{equation}
\partial^{\beta}_x f (g (\cdot)) = \sum_{n=1}^k 
\left( f^{(n)} (g (\cdot)) \sum_{I (\beta, n)}  \prod_{\ell =1}^k C_{\beta, n} ( \partial^{\beta_{\ell}}_x g (\cdot)  )^{v_{\ell}} \right)
\label{A11}
\end{equation} 
where $C_{\beta, n}$ is the corresponding factorial coefficient, and 
\begin{equation*}
I (\beta, n)= \left\lbrace  (v_1, \ldots, v_k, \beta_1, \ldots, \beta_k) : v_{\ell} \geq 0, \sum_{\ell =1}^k v_{\ell} = n, \beta_{\ell} \leq \beta~ \mbox{multi-indices}~
 \sum_{\ell =1}^k v_{\ell} \beta_{\ell} = \beta  \right\rbrace . 
\end{equation*}

\begin{proof}[Proof of Lemma~\ref{Chain Rule}]
Let $k =  1+ \left[ {s}\right]$ and assume that $G \in C^k (\mathbb{R})$ is a function of the form $G(x)= |x|^q$ or $G(x) = |x|^{q-1} x$. 
By the mean value theorem, we have 
\begin{align}
G(u) = G(P_{<N} u ) + G' ( \lambda P_{< N} u + (1 - \lambda) P_{\geq N} u ) P_{\geq N} u 
\label{A10}
\end{align}
where $\lambda = \lambda (x, t, u) \in (0,1)$. Since $G'$ is a power type function and $\lambda \in (0,1)$, applying H\"older's inequality we get
\begin{align}
\begin{split}
\Vert P_N G(u) \Vert_{L^{r}_x}  
&\les \Vert P_N G (P_{< N} u) \Vert_{L^{r}_x}   \\
& \indeq + \left( \Vert G' (P_{< N} u) \Vert_{L^{r_1}_x}   + \Vert G' (P_{\geq N} u) \Vert_{L^{r_1}_x}    \right)   \Vert P_{\geq N} u \Vert _{L^{p_2}_x}   \\
& = I_1 + I_2
\end{split}
\label{A10.2}
\end{align}
where $1/r = 1/r_1 + 1/ p_2$. 

We begin with $I_1$. We apply Fa\`{a} di Bruno's formula with $f = G$ and $g = P_{<N} u$. As $P_{<N} u$ is a smooth function on $\mathbb{R}^d$ and $G$ is assumed to be $C^k$, equation \eqref{A11} combined with H\"older  implies that 
\begin{equation}
\Vert \nabla^k G (P_{N} u) \Vert_{L^{r}_x}   \les \sup_{\substack{k_1 + \cdots + k_m =k \\ k_i \in \{0, \ldots, k \}}}
\Vert P_{< N} u \Vert^{q -m}_{L^{\tilde{r}_0}_x}   \Vert \nabla^{k_1} P_{<N} u \Vert_{L^{\tilde{r}_1}_x} \cdots \Vert \nabla^{k_m} P_{<N} u  \Vert_{L^{\tilde{r}_m}_x}
\label{A12}
\end{equation}
where 
\begin{equation*}
\frac{1}{r} = \frac{q-m}{\tilde{r}_0} + \frac{1}{\tilde{r}_1} + \cdots + \frac{1}{\tilde{r}_m}.
\end{equation*}

In the above formula, we denote by
\begin{equation*}
\Vert \nabla^k f \Vert_{L^{r}_x} := \sup_{|\beta| =k}  \Vert \partial_x^{\beta} f  \Vert_{L^{r}_x}
\end{equation*}
where the supremum is taken over all multi-indices of order $k$. Also, the parameters $\tilde{r}_0, \ldots, \tilde{r}_m$ will be fixed below. 

We utilize \eqref{A12} to estimate $I_1$. By Bernstein's inequality, we have
\begin{equation}
\Vert P_N (G (P_{<N} u) ) \Vert_{L^r_x} \les N^{-k} \Vert \nabla^k G( P_{<N} u) \Vert_{L^r_x}.  
\label{A13}
\end{equation}
By \eqref{A12}, we estimate
\begin{align}
\begin{split}
& \Vert \nabla^k G( P_{<N} u) \Vert  \les \sup_{ k_1 + \cdots + k_m =k }
\Vert P_{< N} u \Vert^{q -m}_{L^{\tilde{r}_0}_x}   \Vert \nabla^{k_1} P_{<N} u \Vert_{L^{\tilde{r}_1}_x} \cdots \Vert \nabla^{k_m} P_{<N} u  \Vert_{L^{\tilde{r}_m}_x} \\
& \indeq \indeq  \les \sup_{ k_1 + \cdots + k_m =k }   \Vert P_{< N} u \Vert^{q -m}_{L^{\tilde{r}_0}_x} 
 \left(  \sum_{N_1 < N} \cdots \sum_{N_m < N}  \Vert \nabla^{k_1} P_{N_1} u \Vert_{L^{\tilde{r}_1}_x}  \cdots \Vert \nabla^{k_m} P_{N_m} u  \Vert_{L^{\tilde{r}_m}_x}  \right) \\
&\indeq  \indeq  \les  \sup_{ k_1 + \cdots + k_m =k }   \Vert P_{< N} u \Vert^{q -m}_{L^{\tilde{r}_0}_x} 
 \left(  \sum_{N_1 <  \cdots<  N_m < N}  \Vert \nabla^{k_1} P_{N_1} u \Vert_{L^{\tilde{r}_1}_x}  \cdots \Vert \nabla^{k_m} P_{N_m} u  \Vert_{L^{\tilde{r}_m}_x}  \right) 
\end{split}
\label{A14}
\end{align}
where in the last step we allowed $N_1 < N_2 < \cdots < N_m < N$ by giving up a scaling factor in the order of $r! \cong O (k)$. Even though the terms $\Vert \nabla^{k_1} P_{N_1} u \Vert, \ldots, \Vert \nabla^{k_m} P_{N_m} u \Vert$ are not symmetric in $N_1, \ldots, N_m$ due to the presence of the derivative operators, this does not matter when we take the supremum in $k_1, \ldots, k_m$.  
Then, using Bernstein's inequality once again, the right side of \eqref{A14} is bounded from above by
\begin{equation}
 \sup_{ k_1 + \cdots + k_m =k }   \Vert P_{< N} u \Vert^{q -m}_{L^{\tilde{r}_0}_x} 
  \left(  \sum_{N_1 <  \cdots<  N_m < N}  N_1^{k_1} \cdots N_m^{k_m} \Vert  P_{N_1} u \Vert_{L^{\tilde{r}_1}_x}  \cdots \Vert  P_{N_m} u  \Vert_{L^{\tilde{r}_m}_x}  \right) 
  \label{A15}
\end{equation}
Also, since each $N_i$ is a dyadic integer, we may use
\begin{align}
 \Vert P_{< N} u \Vert_{L^{\tilde{r}_0}_x} & \les \Vert u \Vert_{L^{\tilde{r}_0}_x} \\
\sum_{N_i < N_{i+1}} N_i^{k_i} \Vert  P_{N_i} u  \Vert_{L^{\tilde{r}_i}_x}  & \les \Vert  u  \Vert_{L^{\tilde{r}_i}_x} N_{i+1}^{k_i}
\comma i = 1, \dots, m-1.
\end{align}
Selecting 
\begin{equation*}
\tilde{r}_0 = \cdots \tilde{r}_{m-1} = p_1 \comma \tilde{r}_m= p_2 ,
 \end{equation*}
 we may then estimate \eqref{A15} by
\begin{equation}
  \Vert u \Vert^{q -1}_{L^{p_1}_x} 
\left( \sum_{N_m < N}  N_m^{k}  \Vert P_{N_m} u \Vert_{L^{p_2}_x } \right) 
\end{equation}
which in turn estimates the right side of \eqref{A14}. All in all, combining \eqref{A13}--\eqref{A15} we obtain
\begin{align*}
I_1 \les N^{-k}  \Vert u \Vert^{q -1}_{L^{p_1}_x} 
\left( \sum_{N_m < N}  N_m^{k}  \Vert P_{N_m} u \Vert_{L^{p_2}_x } \right) 
\end{align*}
which leads to
\begin{align}
\sum_{N} N^{2s} I_1^2 \les  \Vert u \Vert^{2 (q -1)}_{L^{p_1}_x} 
\left( \sum_{N} N^{2s - 2k} \left( \sum_{N_m < N}  N_m^{k}  \Vert P_{N_m} u \Vert_{L^{p_2}_x } \right)^2 \right) 
\label{A16}
\end{align}
Note that by Cauchy-Schwarz we get
\begin{equation}
\left( \sum_{N_m < N}  \left( {N_m}/{N}\right)^{k-s}  N_m^{s} \Vert P_{N_m} u \Vert_{L^{p_2}_x}   \right)^2
\les \sum_{N_m < N} \left( {N_m}/{N}\right)^{k-s} N_m^{2s} \Vert P_{N_m} u \Vert_{L^{p_2}_x}^2  .
\end{equation}
We plug this in \eqref{A16} and change the order of summation to find
\begin{align}
\begin{split}
\sum_{N} N^{2s} I_1^2&  \les  \Vert u \Vert^{2 (q -1)}_{L^{p_1}_x} \sum_{N} \sum_{N_m < N} \left( {N_m}/{N}\right)^{k-s} N_m^{2s} \Vert P_{N_m} u \Vert_{L^{p_2}_x}^2  \\
& \les   \Vert u \Vert^{2 (q -1)}_{L^{p_1}_x}   \sum_{N_m} N_m^{2s} \Vert P_{N_m} u \Vert_{L^{p_2}_x}^2.
\label{A17}
\end{split}
\end{align}
Next, we estimate $I_2$. Using our assumptions on $G$ and \eqref{A.6}, we simply have
\begin{align}
\begin{split}
\Vert G' (P_{<N} u) \Vert_{L^{r_1}_x} & \les \Vert P_{<N} u \Vert_{L^{r_1 (q-1)}_x}^{q-1} \\
& \les \Vert u \Vert^{q-1}_{L^{r_1 (q-1)}_x}.
\end{split}
\label{A18}
\end{align}
Similarly, we note that
\begin{equation}
\Vert G' (P_{ \geq N} u) \Vert_{L^{r_1}_x}  \les \Vert u \Vert^{q-1}_{L^{r_1 (q-1)}_x}. \label{A19}
\end{equation}
As a result, \eqref{A18} and \eqref{A19} combined with the definition of $I_2$ in \eqref{A10.2} yields
\begin{equation}
\sum_{N} N^{2s} I_2^2 \les \Vert u \Vert^{2 (q-1)}_{L^{r_1 (q-1)}_x} \left( N^{2s} \Vert P_{\geq N} u \Vert^{2}_{L^{p_2}_x} \right) .
\label{A20}
\end{equation}
Taking a closer look at the sum on the right hand side above, we expand $P_{\geq N}$ with \eqref{A.5}, and change the order of summation to estimate
\begin{align}
\begin{split}
\sum_{N} N^{2s} \Vert P_{\geq N} u \Vert^2_{L^{p_2}_x} & \les 
\sum_{N}  \left(  \sum_{M \geq N}  (N/ M)^{s}  M^{s} \Vert  P_{M} u \Vert_{L^{p_2}_x} \right)^2 \\
&\les \Vert  u \Vert^2_{\dot{B}^s_{p_2,2}}.
\end{split}
\end{align}
In conclusion, collecting the estimates \eqref{A17} and \eqref{A20} for $I_1$ and $I_2$ respectively,  and setting $p_1 = r_1 (q-1)$ we obtain
\begin{align*}
\Vert G(u) \Vert^2_{\dot{B}^s_{r,2}} & \les 
\sum_{N}  N^{2s} I_1^2 + \sum_{N} N^{2s} I_2^2 \\
& \les  \Vert u \Vert^{2(q-1)}_{L^{p_1}_x}   \Vert  u \Vert^2_{\dot{B}^s_{p_2,2}} 
\end{align*}
where $1/r = (q-1)/p_1 + 1/ p_2$. 
\end{proof}

We will also need the following version of the Leibniz rule in  section~\ref{sec:Decay}.   

\begin{Lemma}
\label{Com_Leibniz Rule}
Let $s \geq 0$ and $m \in \mathbb{N}$. Then, 
\begin{align}
\begin{split}
\norm{D^{m+s} (fg)  - f D^{m+s} g }_{L^2_x} & \les 
\sum_{k =1}^{m-1} \sup_{{I}_k} \norm{( \partial_x^{\beta_{k}} f)~ D^{m+s -k - j} (\partial_x^{\beta_{j}} g  ) }_{L^2_x} \\
& \indeq + \sup_{{I}_m} \norm{(\partial_x^{\beta_m} f) ~g }_{\dot{H}^s_x}
\end{split}
\label{A21}
\end{align}
where
\begin{equation*}
{I}_k= \{ (\beta_{k}, \beta_{j}):  \beta_{k}, \beta_{j} ~ \mbox{ are multi-indices with }~|\beta_{k} |= k, |\beta_{j}| = j \leq k,~  m+s-k-j >0  \}
\end{equation*}
and 
\begin{equation*}
{I}_m = \{ \beta_m: \beta_m ~\mbox{is a multi-index with}~ |\beta_m| = m \}.
\end{equation*}
\end{Lemma}

\begin{proof}
We begin with the Littlewood-Paley expansion of the left hand side in \eqref{A21}.  Recalling \eqref{eq1.7p}, we let 
\begin{equation*}
\psi_N (\xi) = \varphi (\xi / N) - \varphi (2 \xi /N). 
\end{equation*}
By the Plancherel theorem, the expansion becomes
\begin{align*}
\begin{split}
\norm{ D^{m+s} (fg) - f D^{m+s} g }^2_{L^2_x} &  = \sum_{N} \norm{ P_N (  D^{m+s} (fg)  )   - P_N ( f D^{m+s} g )  }^2_{L^2_x}
\\ & = 
\sum_{N} \norm{ \psi_N (\xi) |\xi|^{m+s}  \mathcal{F} (fg) (\xi)   - \psi_N (\xi) \mathcal{F} ( f D^{m+s} g)  (\xi)  }^2_{L^2_{\xi}}
\\ & = 
\sum_N I_N .
\end{split}
\end{align*}
Note that
\begin{align}
\begin{split}
I_N &  = \int | \psi (\xi)  |^2   \left( |\xi|^{m+s} \int \hat{f} (\eta) \hat{g} (\xi - \eta )~d\eta ~-~\int \hat{f} (\eta) | \xi- \eta |^{m+s} \hat{g} (\xi - \eta)~d\eta       \right)^2 d\xi \\
& =  \int | \psi (\xi)  |^2 \left( \int (| \xi |^{m+s}    - |\xi - \eta|^{m+s} ) \hat{f} (\eta) \hat{g} (\xi - \eta) ~d\eta   \right)^2 d\xi .
\end{split}
\label{A22}
 \end{align}
Next, we apply Taylor's expansion on the difference $| \xi |^{m+s}    - |\xi - \eta|^{m+s} $. Letting $F(\xi) = |\xi|^{m+s}$, we consider the case $s>0$, which implies that
$F \in C^{m}$. 
By Taylor's expansion we have
\begin{align}
\begin{split}
F(\xi) - F(\xi - \eta) & = \sum_{1 \leq |\beta_k|  \leq m-1} C_{\beta_k} \partial^{\beta_k} F(\xi - \eta)~ \eta^{\beta_k} \\
& \indeq + \sum_{|\beta_m| = m}  C_{\beta_m} \eta^{\beta_m} \int_0^1 (1 - \lambda)^{m-1} \partial^{\beta_m} F(\lambda \xi +(1- \lambda) (\xi - \eta))~d\lambda
.
\end{split}
\label{A23}
\end{align}
Additionally, we recall that
\begin{equation}
\partial^{\beta} F(\xi) = \sum_{(r_j, \alpha_j) \in I_{\beta}} C_j  |\xi|^{r_j} \xi^{\alpha_j} 
\label{A24}
\end{equation}
where $C_j$ is the corresponding binomial coefficient and the sum is taken over the index set 
\begin{equation*}
I_{\beta} = \{ (r, \alpha)  : r \geq 0,~ \alpha   ~ \mbox{ is a multi-index with }~ \alpha \leq \beta, ~r+ |\alpha| = m+s - |\beta|  \}.
\end{equation*}
In view of \eqref{A24}, we plug \eqref{A23} into $I_N$, and then we combine $\eta^{\beta_k}$ and $|\xi- \eta|^{r_j} (\xi- \eta)^{\alpha_j}$ factors with $\hat{f} (\eta)$ and 
$\hat{g} (\xi - \eta)$ respectively. Applying the Plancherel theorem once again and summing up in $N$ leads to the upper bound on the right hand side of \eqref{A21}.
In the case $s=0$, we basically follow the same steps noting that $F$ is a $C^{m-1}$ function when $m$ is odd which implies that  the  Taylor expansion in \eqref{A23} may be done up to order $m-2$. 
\end{proof}

\section*{Acknowledgments} 
We would like to thank the referees for pointing out some inaccuracies in the previous version of the manuscript. Their comments and suggestions significantly improved the quality of the article. 

The second author  was supported in part by the NSF grant DMS--2153794.

\end{document}